%% file: main.tex
\documentclass[article,12pt,english]{elsarticle}
\usepackage[utf8]{inputenc}
\usepackage{amsmath,amssymb,multirow}
\usepackage{color}
\usepackage{psfrag}
\usepackage{soul}
\usepackage{adjustbox}

\usepackage{multirow}
\usepackage{multicol}
\usepackage{enumerate}
\usepackage{comment}
\usepackage{booktabs}

\usepackage{siunitx} 
\usepackage{longtable}
\usepackage{pdflscape}
\usepackage{algorithm}
\usepackage{algorithmic}
\usepackage{tabularx}
\usepackage{colortbl}

\usepackage{mathabx}
\usepackage{wasysym}

\definecolor{Gray}{gray}{0.95}
\definecolor{LightCyan}{rgb}{0.88,1,1}

\newcolumntype{a}{>{\columncolor{Gray}}c}
\newcolumntype{b}{>{\columncolor{white}}c}

\definecolor{LightRed}{cmyk}{0,0.87,0.68,0.32}

\usepackage{graphicx}

\usepackage{amssymb}

\renewcommand\appendix{\par
  \setcounter{section}{0}
  \setcounter{subsection}{0}
  \setcounter{figure}{0}
  \setcounter{table}{0}
  \renewcommand\thesection{Appendix \Alph{section}}
  \renewcommand\thefigure{\Alph{section}\arabic{figure}}
  \renewcommand\thetable{\Alph{section}\arabic{table}}
}

\newcommand{\FigFormat}{png}

\usepackage[top=0.75in, bottom=0.75in, left=0.5in, right=0.5in]{geometry}

\begin{document}
\graphicspath{{FIG/}}

\begin{frontmatter}



\title{Local time-stepping for adaptive multiresolution using natural extension of Runge--Kutta methods } 


\author[CAP,INPE]{M\"uller Moreira Lopes\corref{cor1}}
\ead{muller.lopes@inpe.br}
\cortext[cor1]{Corresponding author}

\author[LAC,INPE]{Margarete Oliveira Domingues}
\ead{margarete.domingues@inpe.br,margarete.oliveira.domingues@gmail.com}
\author[I2M]{Kai Schneider}
\ead{kai.schneider@univ-amu.fr}
\author[DGE,INPE]{Odim Mendes}
\ead{odim.mendes@inpe.br}
%
\address[CAP]{Graduate program in Applied Computing (CAP)}
\address[LAC]{Associate Laboratory of Applied Computing and Mathematics (LAC), Coordination of the Associated Laboratories (CTE)}
\address[DGE]{Space Geophysics Division (DGE),
Coordination of Space Sciences (CEA)}
\address[INPE]{National Institute for Space Research (INPE), \\
Av. dos Astronautas 1758, 12227-010 São José dos Campos, São Paulo, Brazil}

\address[I2M]{Institut de Math\'ematiques de  Marseille (I2M), Aix-Marseille Universit\'e, CNRS, Centrale Marseille \\ 39 rue F. Joliot--Curie, 13453 Marseille Cedex 13, France}

\begin{abstract}
A space-time fully adaptive multiresolution method for evolutionary non-linear partial differential equations is presented introducing  an improved local time-stepping method.
The space discretisation is based on classical finite volumes, endowed with cell average multiresolution analysis for triggering the dynamical grid adaptation.
The explicit time scheme features a natural extension of Runge--Kutta methods which allow local time-stepping while guaranteeing accuracy. 
The use of a compact Runge--Kutta formulation permits further memory reduction.
The precision and computational efficiency of the scheme regarding CPU time and memory compression are assessed for problems in one, two and three space dimensions.
As application Burgers equation, reaction-diffusion equations and the compressible Euler equations are considered. The numerical results illustrate the efficiency and superiority of the proposed local time-stepping method with respect to the reference computations.
\end{abstract}

\end{frontmatter}




\section{Introduction}

Multiresolution (MR) methods improve the computational performance of numerical solvers of evolutionary partial differential equations when the solution exhibits localised structures, point-wise singularities, boundary layers, shocks, coherent vortices such as encountered in combustion and turbulent flow applications  \cite{Harten:1995,RSTB03,Mueller:2003}.
\par In these methods, the fast wavelet transform is the key ingredient to speed-up computations. 
The wavelet coefficients are employed to measure the local smoothness of the solution. 
Then, a thresholding strategy is used to remove non-significant coefficients, obtaining a grid adapted to the solution. 
This grid is coarser in smooth regions and finer there where structures and steep gradients are present. 
The locally refined grid can be rebuilt interactively to a regular grid  with an expected error directly related to the chosen threshold. 
With the use of this adaptive grid, the number of interface flux computations during the time evolution can be significantly reduced, while controlling the error. 
The adaptive grid is checked before each time step to guarantee that it is sufficiently refined to represent possible new structures in the solution. 
Therefore, the adaptive grid is dynamically adapted to track the solution in scale and space.
\par The next step to improve the computational performance associated with these adaptive grids is to use a local time-stepping (LT) approach, especially when explicit time schemes are used. 
In this work, the combination of MR methods with the local time-stepping approach is denoted as MRLT. 
LT consists in performing the time evolution of each cell on the adaptive grid independently according to its required time step. 
Therefore, each cell must have a time step proportional to its refinement and consequently, a larger cell performs a larger time step. 
\par Local time-stepping methods for space adaptive discretisations of partial differential equations have a long tradition, going back to the early work of Osher and Sanders \cite{OsherSanders:1983}. 
Related to multiresolution, M\"uller and Stiriba \cite{MullerStriba:2007} presented some general MRLT schemes that could be applied either to an explicit time scheme, based on Lagrange projection, or an implicit time scheme, for a reference finite-volume method in space. In \cite{MullerStriba:2007} they applied LT to one-dimensional scalar conservation laws. Following this work, Coquel \textit{et al.} \cite{Coqueletal:2010} presented MRLT methods with both explicit and implicit Lagrange-projection schemes, the latter is the novelty concerning \cite{MullerStriba:2007}. Hejazialhosseini \textit{et al.}  combined in \cite{Hejazialhosseinietal:2010} first and second order RK schemes to propose an LT approach for MR in blocks with finite volumes for multi-phase compressible flows implemented on multi-core architectures.


\par More recently, LT methods have been proposed using different approaches for time interpolation required in the algorithms for providing values at intermediate time steps. In \cite{Qi:2018}, the authors use a high-order Taylor type integrator, while a discontinuous Galerkin method with spectral elements is used for spatial discretization. In the context of finite elements methods, \cite{Rietmann:2017} proposed an energy conserving LT algorithm based on a second-order leap-frog scheme.
 Discussions on the proper choice of the time-step in LT schemes are still an important topic of discussion. To this end \textit{a posteriori} error estimators are used in \cite{Auzinger:2018,Mayr:2018}, while Gnedin {\it et al.} (2018) \cite{Gnedin:2018} investigate their effect by enforcing the CFL condition locally over every cell.
 In the context of finite volumes methods, the use of high-order schemes for local time-stepping on adaptive mesh refinement grids were previously discussed in~\cite{Dumbser:2013}. In that work, the accuracy order in space is obtained through a WENO reconstruction, while the accuracy in the time discretization is achieved by a local space-time discontinuous Galerkin predictor method. This approach yields up to fourth order for compressible Euler equations in 2D. Similar and related works in this context were carried out in \cite{Krivodonova:2010} and \cite{Boscheri:2015}.

A detailed discussion on the stability of those MRLT schemes is presented in  \cite{Hovhannisyanetal:2010}. 
In the context of adaptive numerical methods for partial differential equations, the derivation of explicit LT methods based on standard RK schemes typically stays at orders smaller or equal to two \cite{Domingues20083758,DiazGrote:2009,DiazGrote:2015}
The reason why LT methods are limited to low order in time is because a time synchronization is required; for a discussion we refer to \cite{Domingues20083758}. 
 Recently, higher-order LT schemes have been proposed in the context of discontinuous Galerkin methods. 
The works of Gassner \textit{et al.} \cite{gassner2011runge,Gassneretal:20114232} in this context are based on natural continuous extensions for Runge--Kutta methods (NERK). 
Such schemes have been introduced initially by Zennaro in the late $1980$'s for solving general ODEs, with application to delay equations \cite{zennaro1986natural}. 
The idea is to interpolate the intermediate stages of the Runge--Kutta scheme to obtain the values at the requested intermediate time instants required for the time synchronisation in LT schemes.
This method is also called Natural Continuous Extension  in \cite{zennaro1986natural}, Continuous Extension Runge--Kutta  in \cite{OwrenZennaro:1992}, and in a general way Continuous Runge--Kutta, as discussed in \cite{VERMIGLIO1993521}.
\par Recently, an alternative for higher-order LT methods, again in the context of discontinuous Galerkin spectral element schemes, has been proposed by Winters and Kopriva \cite{winters2014high}. 
The underlying ideas are Adams--Bashforth multi-step schemes. 

Moreover, Gassner {\it et al.} (2011) considered a family of explicit one-step time discretisations for finite volume (FV) and discontinuous Galerkin schemes, which are based on a predictor-corrector formulation \cite{Gassneretal:20114232}.

\par The aim here is to use the idea of Gassner \textit{ et al.} \cite{gassner2011runge} using NERK for the first time in the context of MRLT methods to perform the synchronisation. 
Thus, an improvement in time accuracy can be obtained. 
In the current work, only second and third-order schemes are used, but the extension to higher-order is in principle possible. 
Gassner \textit{ et al.} \cite{gassner2011runge} discussed that NERK is $10\%$ slower than the 
Cauchy–-Kovalevskaya scheme. 
However, to work with the latter is complicated as the analytic solution of some nonlinear PDEs is required. 
The goal of the proposed approach is to perform simulations using NERK schemes. 
Therewith, the synchronisation required for the LT approach is possible, and this new class of MRLT methods, named MRLT/NERK, is created. 
This work presents the application and implementation of the MRLT/NERK method for the second (RK2) and third (RK3) order time evolution, named MRLT/NERK2 and MRLT/NERK3, respectively. 
The method proposed in this work is applied for solving the two-dimensional Burgers equation, one and three-dimensional reaction-diffusion problems and the two-dimensional compressible Euler equations. 
The obtained results and CPU times are compared with the MR method using classical RK2 and RK3 time evolution, named MR/RK2 and MR/RK3, respectively. The results are also compared with the MRLT approach given in \cite{Domingues20083758} based on RK2 time evolution, named MRLT/RK2.  

In Section \ref{MR}, we summarise the adaptive multiresolution method proposed in \cite{Harten:1995}. Then, in Section~\ref{Runge-kutta}, we discuss the Runge--Kutta methods and the NERK methods used in the current work to perform the proposed MRLT/NERK approach given in Section~\ref{LT}.
A convergence analysis is conducted in Section~\ref{sec:stability}.
Performance comparisons, considering CPU time and errors, among the FV, MR and MRLT approaches given in \cite{Domingues20083758} and the MRLT/NERK approach introduced here, are presented in Section~\ref{results}.      
Conclusions are drawn in Section~\ref{sec:conclusions}.

\section{Adaptive multiresolution methods using finite volumes}\label{MR}
The following initial value problem for a vector-valued function $\mathbf{Q}(\mathbf{x},t)$ written in divergence form is considered:
\begin{equation} 
\frac{\partial \mathbf{Q}}{\partial t} =
- \nabla \cdot \mathbf{F}(\mathbf{Q}, \nabla \mathbf{Q}) + \mathbf{S}(\mathbf{Q}),  
\quad \text{for} \quad (\mathbf{x},t) \in \Omega \times [0, +\infty ),\, \Omega \subset \mathbb{R}^{\mathit{d}}.
\label{eqn:pde}
\end{equation}
This problem is given in $\mathit{d}$ space dimensions, completed with initial conditions $\mathbf{Q}(\mathbf{x},t=0)  =  \mathbf{Q}_0 (\mathbf{x})$ and appropriate  boundary conditions. The terms $\nabla \cdot \mathbf{F}(\mathbf{Q}, \nabla \mathbf{Q})$ and $\mathbf{S}(\mathbf{Q})$ denote the divergence and source term, respectively. The flux $\mathbf{F}$ can be decomposed into advective and diffusive contributions, \textit{i.e.} $\mathbf{F}(\mathbf{Q}, \nabla \mathbf{Q}) = \mathbf{f}(\mathbf{Q}) - \nu \nabla \mathbf{Q}$, where the diffusion coefficient $\nu$ is positive and assumed to be constant.

\par To discretise Equation (\ref{eqn:pde}) in space, we use a classical finite volume formulation written in standard form. The domain $\Omega$ corresponds to a rectangular parallelepiped in $\mathit{d}=1$, $2$ or $3$ dimensions in Cartesian geometry. It is partitioned into cells $(\Omega_i)_{i \in \Lambda}$, $\Lambda = \{ 0, \ldots, i_{\max}\}$ with $\Omega = \bigcup_{i} \Omega_i$.
\par Defining the volume of the cell by $|\Omega_i| = \int_{\Omega_i} d\mathbf{x}$, we compute the cell-average $\bar{\mathbf{q}}_i(t)$ of a given quantity $\mathbf{Q}$ on $\Omega_i$ at time instant $t$ by,
\[
\bar{\mathbf{q}}_i (t) = \frac{1}{|\Omega_i|} \int_{\Omega_i} \mathbf{Q}(\mathbf{x}, t) \, d\mathbf{x}.
\]
Considering the one-dimensional case, $\Omega_i$ is an interval $[x_{i-\frac{1}{2}}, x_{i+\frac{1}{2}}]$ of length $\Delta x_i = x_{i+\frac{1}{2}} - x_{i-\frac{1}{2}}$. Integrating Equation~(\ref{eqn:pde}) on $\Omega_i$ then yields:
\begin{equation} \label{eqn:int}
\frac{d\bar{\mathbf{q}}_i}{dt} (t) = - \frac{1}{\Delta x_i} \left( \bar{\mathbf{F}}_{i+\frac{1}{2}} - \bar{\mathbf{F}}_{i-\frac{1}{2}} \right) + \bar{\mathbf{S}}_i,
\end{equation}
where $\bar{\mathbf{F}}$ is the numerical flux and $\bar{\mathbf{S}}_i$ is the source term of the cell $\Omega_i$. This formulation can be extended to two- and three dimensional problems.

\par The source term is approximated by $\bar{\mathbf{S}}_i \approx \mathbf{S}(\bar{\mathbf{q}}_i)$, which yields also second-order accuracy in case of a linear source term.

\par
\begin{figure}[htb]
\begin{center}
\includegraphics[scale=0.75]{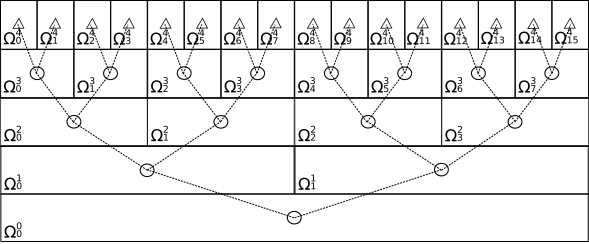} \\
\footnotesize{a) Grid hierarchy for an unidimensional domain.} \\[0.2cm]
\includegraphics[scale=0.75]{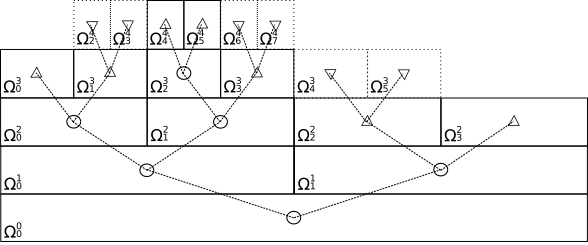} \\
\footnotesize{b) Tree structure for an adapted grid.} \\[0.2cm]
\end{center}
\caption{Dyadic grid hierarchy of MR methods and its implementation in a tree data structure. a) Example of nested unidimensional grids, the circles represent the internal nodes of the tree structure, while the triangles represent the leaves of the tree. b) Example of an adapted grid using the tree structure from a). The virtual leaves are represented by the triangles $\triangledown$. Adapted from \cite{RSTB03}.}
\label{fig:tree23D}
\end{figure}

 The adaptive multiresolution (MR) analysis in the cell average context, proposed by Harten \cite{Harten:1995}, consists in decomposing the cell-averages of the solution into a multilevel representation. This representation, illustrated in Figure~\ref{fig:tree23D}a, consists of a hierarchy of nested grids $\Omega^\ell$, where $\ell$ is the grid refinement level. Each grid $\Omega^\ell$ consists of a regular grid, as defined for the FV formulation, with $2^{\mathit{d}\ell}$ cells. A cell of refinement level $\ell$ and position $i$ is denoted by $\Omega^\ell_i$.

\par Figure \ref{fig:tree23D}a shows the implementation of this structure as a binary tree, where the nodes of level $\ell$ generate the grid $\Omega^\ell$. For the two- and three-dimensional cases, this idea is extended by using a quadtree and an octree, respectively. 

\par This sequence of nested grids corresponds to a scheme where a finer scale represents its subsequent coarser scale plus a sum of details between these levels. These details are the wavelet coefficients. The construction of the adapted grid in the MR method consists in representing the grid using only cells with significant wavelet coefficients. For this, the leaves of the most refined level are checked, if their parent cell has a significant wavelet coefficient, larger than a predetermined threshold $\epsilon$, if not, the leaves are deleted. 
This process is repeated recursively starting from the finest level, going to coarser and coarser levels. 
The procedure to obtain the projection of the solution on the coarser levels and the wavelet coefficients is given in Section~\ref{sec:projection}.

\par The grid generation procedure is mathematically supported by the fact that the magnitude of of the wavelet coefficients are small in regions where the solution is smooth, while they are significant in regions where the solution exhibits steep gradients.
Hence, high compression rates are expected when only a small number of cell-averages is present on the finest scales.  

\par In this work, the data structure for representing the solution is organized as a dynamic graded tree. This tree organization requires that no hole is admitted inside the tree, which means that the connectivity in the tree structure has to be ensured. Moreover, the tree can change in time to track the space-scale evolution of the solution as nodes can be added or removed while guaranteeing its gradedness.
\par This dynamic graded tree organization implies that neighbours of each cell can have a difference of one refinement level at most. This restriction allows the use of virtual leaves for the flux computations between leaves at different refinement levels. These virtual leaves are auxiliary leaves placed as children nodes of the leaves who have an interface with finer leaves. Their values are predicted using the prediction procedure used to compute the wavelet coefficients. A unidimensional adaptive grid represented in a graded tree structure is shown in Figure~\ref{fig:tree23D}b. 
     
\par The flux computations in the MR scheme are performed individually for each leaf. The numerical fluxes are computed between cells at the same refinement level. Using the adaptive grid of Figure~ \ref{fig:tree23D}b, the following interface scenarios for computing the numerical fluxes of a leaf can be identified:  
\begin{itemize}
    \item \textbf{Leaf / Leaf or Virtual leaf:} When both leaves belong to the same refinement level, the flux computation is performed in the same way as in the FV method. The following cases can occur: Leaf/Leaf example: Flux between cells $\Omega^3_0$ and $\Omega^3_1$; Leaf/Virtual leaf example: Flux between cells $\Omega^4_5$ and $\Omega^4_6$. 
    \item \textbf{Leaf / Internal node:} In this scenario, the current leaf has an interface with a finer leaf. To perform the flux computation in this case, the numerical flux is computed using the virtual children of the leaf and its adjacent leaves. Example: Flux between cells $\Omega^2_2$ and $\Omega^2_1$. 
\end{itemize}

\par The construction of the adaptive grid guarantees that the current solution is well represented.  However, after the time evolution process, the new solution should also be well represented on this grid, which a priori cannot be ensured. 
In order to guarantee that the new solution after the time evolution is still well represented on this grid, the leaves with finer neighbours are refined and neighbor cells are added. Further details of the MR scheme and its implementation can be found in \cite{RSTB03}. 
The Algorithms \ref{alg:AdaptiveMesh} and \ref{alg:waveletThres} given in Appendix \ref{App:Alg} describe the adaptive grid creation and its update, respectively.

\subsection*{Projection and prediction operators}\label{sec:projection}

In order to perform the MR method, some operations for projection and prediction are required. For the MR scheme with finite volumes, where the cell values are local averages, a coarser cell $\Omega^\ell_i$ has its value estimated using the finer values and an unique projection operator $P_{\ell+1 \rightarrow \ell} \, : \, \bar{\mathbf{q}}^{\ell+1} \, \mapsto \, \bar{\mathbf{q}}^{\ell}$. In this scheme, the projection operator to obtain the solution of a coarser cell is given by the average value of its children. For the unidimensional case, the projection is performed by:
\begin{equation}
    \bar{\mathbf{q}}^\ell_i = P_{\ell+1 \rightarrow \ell} \left(\bar{\mathbf{q}}^{\ell+1}_{2i}, \bar{\mathbf{q}}^{\ell+1}_{2i+1}\right) = \frac{1}{2}\left( \bar{\mathbf{q}}^{\ell+1}_{2i} + \bar{\mathbf{q}}^{\ell+1}_{2i+1}\right),
\end{equation}
where $\bar{\mathbf{q}}^{\ell}_{i}$ are the average value of the cell $\Omega^\ell_i$. The same idea is extended for the two and three-dimensional cases.

\par The prediction operators are used to perform the opposite path of the projection operators, they allow to obtain the values of the finer cells using the values of the coarser ones. 
For each child cell to be predicted, there is a different prediction operator, represented by $P^i_{\ell \rightarrow \ell+1} \, : \, \bar{\mathbf{q}}^{\ell} \, \mapsto \, \bar{\mathbf{q}}^{\ell+1}$ for the one-dimensional case, $P^{i,j}_{\ell \rightarrow \ell+1} \, : \, \bar{\mathbf{q}}^{\ell} \, \mapsto \, \bar{\mathbf{q}}^{\ell+1}$ for the two-dimensional case and $P^{i,j,k}_{\ell \rightarrow \ell+1} \, : \, \bar{\mathbf{q}}^{\ell} \, \mapsto \, \bar{\mathbf{q}}^{\ell+1}$ for the three-dimensional case. 
These operators yield a non-unique approximation of $\bar{\mathbf{q}}^{\ell+1}_i$ by interpolation.   
We use polynomial interpolation of second degree on the cell-averages, as proposed by Harten \cite{Harten:1995}, which yields third-order accuracy. For the one-dimensional case, it  follows that,
\begin{subequations}
\begin{equation}
\tilde{\mathbf{q}}^{\ell+1}_{2i}  = P^0_{\ell \rightarrow \ell+1} \left(\bar{\mathbf{q}}^{\ell}_{i-1}, \bar{\mathbf{q}}^{\ell}_{i}, \bar{\mathbf{q}}^{\ell}_{i+1}\right) = \bar{\mathbf{q}}^\ell_{i}- \frac{1}{8} (\bar{\mathbf{q}}^\ell_{i+1} - \bar{\mathbf{q}}^\ell_{i-1})
\end{equation}
\begin{equation}
\tilde{\mathbf{q}}^{\ell+1}_{2i+1}  = P^1_{\ell \rightarrow \ell+1} \left(\bar{\mathbf{q}}^{\ell}_{i-1}, \bar{\mathbf{q}}^{\ell}_{i}, \bar{\mathbf{q}}^{\ell}_{i+1}\right) = \bar{\mathbf{q}}^\ell_{i}+ \frac{1}{8} (\bar{\mathbf{q}}^\ell_{i+1} - \bar{\mathbf{q}}^\ell_{i-1}).
\end{equation}
\end{subequations}
where $\tilde{\mathbf{q}}^\ell_{i}$ is an approximation of the value $\bar{\mathbf{q}}^{\ell}_{i}$. Interpolation operators for higher dimensions can be found in \cite{RSTB03}. The operator must satisfy the properties of locality, requiring the interpolation for a child cell to be computed from the cell-averages of its parent and its nearest uncle cells in each direction; and consistency, $P_{\ell+1 \rightarrow \ell} \circ P_{\ell \rightarrow \ell+1} = \mbox{Identity}$. 
\par The prediction operator is used to obtain the wavelet coefficients $\mathbf{d}^\ell_i$ of the finer cells. These coefficients are given by the difference between the cell average $\bar{\mathbf{q}}^{\ell}_i$ and the predicted value $\tilde{\mathbf{q}}^{\ell}_{i}$:
\begin{equation}
    \mathbf{d}^\ell_i = \bar{\mathbf{q}}^{\ell}_{i} - \tilde{\mathbf{q}}^{\ell}_{i}.
    \label{eq:waveletcoeff}
\end{equation}
The values $\mathbf{d}^\ell_i$ are also used for reconstructing the finer levels without interpolation errors. Their norm yields the local approximation error. Moreover, the information of the cell-average value of the two children is equivalent to the knowledge of the cell-average value of the parent and one independent detail. This can be expressed by 
$\left( \bar{\mathbf{q}}^{\ell+1}_{2i}, \bar{\mathbf{q}}^{\ell+1}_{2i+1} \right) \longleftrightarrow  \left( \mathbf{d}^{\ell+1}_{2i} ,\bar{\mathbf{q}}^{\ell}{i} \right).$ 
This procedure can be applied recursively from level $L$ down to the level $0$ creating thus a multiresolution transform of the cell-average values as proposed by Harten~\cite{Harten:1995}. 
Therefore, we have
\begin{equation}
\bar{\mathbf{q}}^{L} \longmapsto (\bar{D}^{L},\, \bar{D}^{L-1},\, \ldots, \bar{D}^1,\, \bar{\mathbf{q}}^0),
\end{equation}
where $\bar{D}^{\ell}$ is the set of wavelet coefficients at level $\ell$.
Accordingly, the information of the cell-average values of all the leaves is equivalent to the knowledge of the cell-average value of the root and the wavelet coefficients of all the other nodes of the tree structure. For two and three dimensions, respectively, the information of the cell-averages of four and eight children is equivalent to the knowledge of three and seven wavelet coefficients in the different directions and the node value \cite{Harten:1995,RSTB03}.


\section{Runge--Kutta methods}\label{Runge-kutta}
After discretising in space the initial value problem given in Eq. (\ref{eqn:pde}), the following system of ordinary differential equations in time is obtained:
\begin{equation}
\frac{d\bar{\mathbf{q}}}{dt} = f(t,\bar{\mathbf{q}})\; ,
\end{equation}
where $\mathbf{Q}(t=0)=\bar{\mathbf{q}}^0$ is the given initial condition. This system yields an equation for each leaf of the grid.
By abuse of or to simplify notation the space discretised solution will be denoted again by $\bar{\mathbf{q}}$. 
The general formulation for an explicit $s$-stage Runge--Kutta (RK) method can be expressed at time $t^{n\!+\!1}$ as:
\begin{equation}
\bar{\mathbf{q}}^{n+1} = \bar{\mathbf{q}}^n + \sum_{i=1}^{s} b_i k_i \;,
\label{eq:RKs}
\end{equation}
with
\begin{equation}
k_{i} = \Delta t^n \, f\left( t^n + c_i\Delta t^n, \bar{\mathbf{q}}^n + \sum_{j=1}^{i-1}a_{ij}k_{j}\right)\;,
\label{eq:RKk}
\end{equation}
where $a_{ij}$, $b_i$ and $c_i$ are the Runge--Kutta coefficients, and $\Delta t^n$ is the time-step used to perform the time evolution from the instant $t^n$ to $t^{n+1}$. The actual convergence order of the RK method depends of the number of stages and a set of well selected RK coefficients. 
\par In this work we consider second order RK methods (RK2) with coefficients $c_{i}$ and $a_{ij}$ in such a way that they are the same coefficients used in the first and second steps of the RK3 method. Namely, the values of these coefficients are 
 $c_1=0$ and $c_2=a_{21}=1$. The other coefficients for RK2 are $b_1=b_2=\dfrac{1}{2}$. To perform RK3, further coefficients are $c_3=\dfrac{1}{2}$, $a_{31}=a_{32}= \dfrac{1}{4}$, $b_1=b_2=\dfrac{1}{6}$ and $b_3=\dfrac{2}{3}$.
 
\subsection*{Natural continuous Extension Runge--Kutta (NERK) method}\label{NERK}

The NERK method, originally introduced by \cite{zennaro1986natural}, produces an approximation of the solution in the time interval $[t^n;t^n+\Delta t^n]$ using the same coefficients $a_{ij}$ and $c_i$ as in the standard Runge--Kutta methods. The difference between NERK and RK is the use of polynomials $\beta_i$ instead of the constant coefficients $b_i$. 
\par The NERK method can be expressed as,
\begin{equation}
 \bar{\mathbf{q}}(t^n + \theta \Delta t^n) = \bar{\mathbf{q}}^n + \sum_{i=1}^{s} \beta _i(\theta) \; k_i , \; \theta \in (0,1], 
 \label{rkcont}
\end{equation}
where the polynomials $\beta_i$ are given as a function of the coefficients $b_i$ of the original RK method. 

\par Using the proposed RK2 coefficients, the following polynomials for the two-stage NERK method are obtained using the methodology given in \cite{OwrenZennaro:1992}:
\begin{equation}
 \beta _1(\theta) = -\frac{1}{2}\theta^2 + \theta , \quad
 \beta _2(\theta) = \frac{1}{2}\theta^2. 
 \end{equation} 

\par In this work, the application of the NERK method consists in producing approximations of a solution at some intermediate time instants inside the interval $\left[t^n, t^n + \Delta t^n\right]$, depending on the choice for RK2 or RK3 time evolution. 

\subsection*{Compact formulation}\label{CompactFormulation}

In order to reduce the memory allocation per cell when performing the time evolution, the compact formulation of the RK methods is of particular interest, used for example in \cite{RSTB03}. 
Based on the standard RK methods, we can obtain the following compact formulation for the two and three stage methods,
 \begin{itemize}
 \item RK2: \centerline{$\bar{\mathbf{q}}^* = \bar{\mathbf{q}}^n +\Delta t^n  f(t^n,\bar{\mathbf{q}}^n)$, $\qquad \bar{\mathbf{q}}^{n+1} = \frac{1}{2}\bar{\mathbf{q}}^n  + \frac{1}{2} \bar{\mathbf{q}}^*  +  \frac{1}{2}\Delta t^n   f\left(t^n + \Delta t^n, \bar{\mathbf{q}}^*\right)$.}
 	\item RK3: \centerline{$\bar{\mathbf{q}}^* = \bar{\mathbf{q}}^n +\Delta t^n f(t^n,\bar{\mathbf{q}}^n)$, $\qquad \bar{\mathbf{q}}^{**} =  \frac{3}{4}\bar{\mathbf{q}}^n + \frac{1}{4}\bar{\mathbf{q}}^* + \frac{1}{4}\Delta t^n f\left(t^n + \Delta t^n, \bar{\mathbf{q}}^*\right)$ ,}\\
\centerline{$\bar{\mathbf{q}}^{n+1} =  \frac{1}{3}\bar{\mathbf{q}}^n+ \frac{2}{3}\bar{\mathbf{q}}^{**}  + \frac{2}{3}\Delta t^nf\Big(t^n+\frac{1}{2}\Delta t^n, \bar{\mathbf{q}}^{**}\Big)$.}
 \end{itemize}

\par When performing the time integration using a local time-stepping approach, \textit{cf.} Section~\ref{LT}, some intermediate values are necessary. In this work, we propose to use the NERK method for that. However, the values $k_i$ are not stored in the RK compact formulation. Hence the NERK solution must be adapted to be compatible with the compact RK method.

\par The following steps produce a NERK/RK2 approximation at the time instant $t^n + \frac{1}{2}\Delta t^n$. Due to the memory management of the numerical code, where the fluxes are stored at the same memory allocation, each one of the following steps is executed immediately after its corresponding compact RK step:
\begin{equation}
\bar{\mathbf{q}}^{*}_{\theta = \frac{1}{2}} = \bar{\mathbf{q}}^n +\frac{3}{8}\Delta t^n f(t^n,\bar{\mathbf{q}}^n), \qquad
\bar{\mathbf{q}}_{\theta = \frac{1}{2}} = \bar{\mathbf{q}}^{*}_{\theta = \frac{1}{2}}+ \frac{1}{8}\Delta t^n f\left(t^n + \Delta t^n, \bar{\mathbf{q}}^*\right).
\end{equation}
where the values $f(t^n,\bar{\mathbf{q}}^n)$ and $f\left(t^n + \Delta t^n, \bar{\mathbf{q}}^*\right)$ are the same as those obtained for the compact RK formulation. In this formulation, the values $\bar{\mathbf{q}}_{\theta = \frac{1}{2}}$ and $\bar{\mathbf{q}}^{*}_{\theta = \frac{1}{2}}$ are stored at the same memory allocation.

\par For the RK3 time evolution, second order approximations at the time instants $t^n + \frac{1}{4}\Delta t^n$ and $t^n + \frac{3}{4}\Delta t^n$ become necessary. These approximations are required to compute the third step of the RK3. However, the NERK/RK3 method only yields these information after the third step. The proposed solution for this problem is to obtain these approximations via NERK/RK2. Then, it is possible to obtain approximations at these desired instants immediately after the second step of the RK method, under the condition that the coefficients $a_{ij}$ and $c_i$ of both, RK2 and RK3, methods are the same. Therefore we use the following approximations:

\begin{itemize}
\item at $t^n + \frac{1}{4}\Delta t^n$: $\qquad \bar{\mathbf{q}}^{*}_{\theta  =  \frac{1}{4}} = \bar{\mathbf{q}}^n +\dfrac{7}{32}\Delta t^n f(t^n,\bar{\mathbf{q}}^n), \qquad \bar{\mathbf{q}}_{\theta = \frac{1}{4}}  =  \bar{\mathbf{q}}^{*}_{\theta = \frac{1}{4}}+ \dfrac{1}{32}\Delta t^n f\left(t^n + \Delta t^n, \bar{\mathbf{q}}^*\right)$.

\item at $t^n + \frac{3}{4}\Delta t^n$: $\qquad \bar{\mathbf{q}}^{*}_{\theta = \frac{3}{4}} =  \bar{\mathbf{q}}^n +\dfrac{15}{32}\Delta t^n f(t^n,\bar{\mathbf{q}}^n), \qquad 
\bar{\mathbf{q}}_{\theta = \frac{3}{4}} =  \bar{\mathbf{q}}^{*}_{\theta = \frac{3}{4}}+ \dfrac{9}{32}\Delta t^n f\left(t^n + \Delta t^n, \bar{\mathbf{q}}^*\right)$.
\end{itemize} 
These approximations can be performed using similar memory management ideas as for the RK2 evolution.


\section{Local time-stepping}\label{LT}

To improve further the computational efficiency of the MR method, a local time-stepping approach was proposed in \cite{Domingues20083758}. This approach consists in using an adapted time-step for each leaf individually. This time step is obtained accordingly to the spatial size of the leaf. Thus, small time steps are only used for fine scale leaves, while large time-steps can be used for coarser leaves. This is possible without violating the stability condition of the explicit time discretisation, as shown in \cite{Domingues20083758}. 

\par Solutions with point-like singularities are well adapted and yield highly efficient multiresolution representations. For those, the MRLT approach is found to be most efficient. Besides, the adaptive grid created for the MRLT scheme is the same graded tree structure used in the MR scheme \cite{RSTB03}. The update procedure for the trees during MRLT schemes is discussed in Section \ref{AdaptTree}.

\par We suppose that the CFL condition implies a time-step $\Delta t$ for the most refined level. In the LT scheme, each cell of level $\ell$ performs its time evolution with a proper time step given by:
\begin{equation}
\Delta t_{\ell}= 2^{L-\ell}\Delta t,
\label{dt}
\end{equation}      
where $L$ is the finest level of the grid. From this point, the notation is adjusted in order to facilitate the understanding. Moreover, considering the Courant number $\sigma$ depending on the ratio between $\Delta t_\ell$ and $\Delta x_\ell$, and using the relations between the cell size and the time-step of different levels, the value $\sigma$ obtained for every cell will thus be the same.

\par
\begin{figure}[htb]
\centering
\includegraphics[scale=0.5]{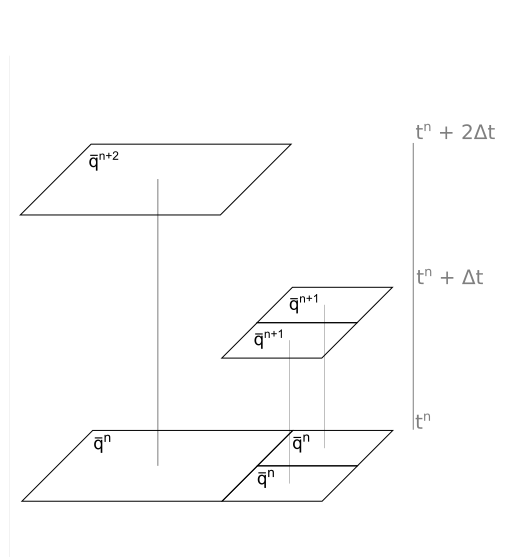}
\caption{LT evolution of adjacent cells at different scales. The finer cells in this representation need a second time evolution to reach the instant $t^{n+2}$.}
\label{fig:NERKSCH}
\end{figure}

 As illustrated in Figure \ref{fig:NERKSCH}, due to the scale dependent time-stepping of the LT scheme, not all leaves of the coarser scales will be evolved during an iteration of the time evolution. We define the coarsest level where the leaves must be evolved in a certain iteration $n$ to be the minimum scale level in which the modulo operator between $n$ and $2^{L-\ell}$ is zero and we denote it as $\ell_{\min}$. 

\par The above approach is restricted to second order time accuracy, because the internal steps of standard Runge--Kutta schemes are not compatible with the dyadic grid size.
The reason is that the intermediate time steps of higher order RK schemes (order larger than two) do not correspond to the time instants of the solution obtained by the RK method when using the dyadic grid size \cite{Domingues20083758}. 
A possible solution to overcome this limitation is to use NERK schemes. 
Their polynomial approximation in time is used to evaluate the solution at the intermediate time instants imposed by the dyadic spatial grid size.

\par The implementation of high order LT schemes leads to three synchronization challenges during the time evolution of a leaf with a neighbour at a different refinement level. Those challenges are discussed in the following subsections. The whole LT method as proposed in this work is performed as in Algorithm \ref{alg-LTS}.  

\subsection{Synchronization during the Runge--Kutta iteration}\label{LTSync1}
Considering the restriction that during an iteration of LT schemes only leaves of refinement level greater or equal $\ell_{\min}$ are evolved, we know that all of these leaves have their solutions at the same time instant $t^n$, requiring no synchronization to perform the first RK step.
Thus, the first RK step can be done with the MR method using the proper $\Delta t_{\ell}$ value. 
After performing the first step of the Runge--Kutta method, the leaves at each level are evolved with their own time step, obtaining a first order solution at time instant $t^n + \Delta t_{\ell}$. 
Using the fluxes obtained in this step, the values $\bar{\mathbf{q}}^{*}_{\theta = \frac{1}{2}}$, 
for RK2, or 
$\bar{\mathbf{q}}^{*}_{\theta = \frac{1}{4}}$ and $\bar{\mathbf{q}}^{*}_{\theta = \frac{3}{4}}$ 
for RK3, are computed at the leaves of every level evolved in this iteration.
\par However, for the next RK steps, due to the different time step size, the solution values after the RK step are given in a different time instant for each refinement level. This implies that some synchronization has to be performed.

\subsubsection{Second Runge--Kutta step synchronization}\label{RKstep2}
\par
\begin{figure}[htb]
    \centering
    \begin{tabular}{cc}
    \includegraphics[width=0.35\linewidth]{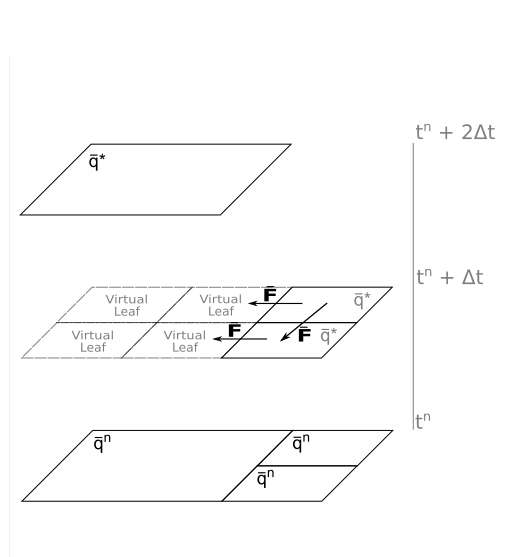} & \includegraphics[width=0.35\linewidth]{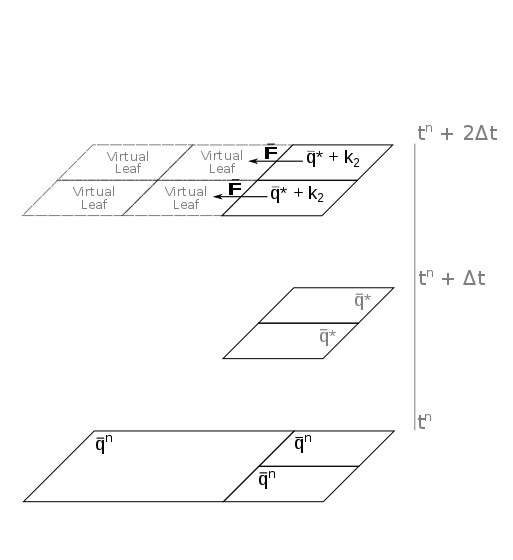}\\
   a) Evaluation of the fine leaves.  & b) Evaluation of the coarse leaves.  \\
    \end{tabular}
    \caption{Schematic view of the synchronization issues for the second RK step. 
    a) presents  the  evaluation of flux computations $\bar{\mathbf{F}}$ at instant $t^n + \Delta t$. 
    For adjacent leaves at the same level we can compute the fluxes directly. 
    However, to compute the fluxes between different levels, we need the virtual leaves of the coarser leaf at the same time instant of the finer leaves.
    For that, first, we interpolate the coarser leaf $\bar{\mathbf{q}}$ at instance  $t^n + \Delta t$  from its respective values at $t^n$ and $t^n + 2 \Delta t$. 
    Then, we predict its value at the finer level and obtain its virtual leaves. 
    After that, we evaluate the flux, and then, the RK step for the finer leaves. 
    b) presents the evolution of the flux at  $t^n + 2 \Delta t$ of the coarser leaf. 
    First, we predict the virtual leaves of the coarse leaf, then we compute $\bar{\mathbf{q}}^* + k_2$, \textit{i.e.} a RK1 step of the finer leaves, which yields a first order approximation.
    After that, we can compute the flux and evaluate the leaves.
    }
    \label{fig:Rk1SyncIssue}
\end{figure}

 In order to perform the second step of the RK method, it is necessary to compute the flux $f(t^n + \Delta t_{\ell}, \bar{\mathbf{q}}^*)$. %
This flux can be interpreted as the flux between leaves at the time instant $t^n + \Delta t_\ell$, where the leaf value is a first order approximation $\bar{\mathbf{q}}^* = \bar{\mathbf{q}}^n + a_{21} k_1$. 
The challenge here lies in the fact that when this approximation is obtained at a finer scale, due to the different time step size used, it is located at an earlier time instant with respect to the approximation at a coarser scale. 
This situation implies that the flux computation for each scale has to be performed at a different time instant, as shown in the scheme presented in Figure~\ref{fig:Rk1SyncIssue}. 
The values required for both cases of the second RK step are not available after the first RK step, requiring thus a synchronization procedure to obtain those values for each scale. Then the second RK step in this scale can be performed. 
To obtain the missing values to compute the fluxes in the situation given in Figure~\ref{fig:Rk1SyncIssue}(a), a tree refreshing procedure, described in Section \ref{syncChallengeTree}, must be performed in order to project the solution for scales $\ell > \ell_{\min}$ at time instant $t^n + \Delta t_{\ell}$ onto the scale $\ell-1$ at instant $t^n + \Delta t_{\ell-1}$. 

\par The proposed synchronization methodology to perform the second RK step consists in using a first order approximation $\bar{\mathbf{q}}_\ell^n + \frac{1}{2}a_{21} k_1$ as a solution at the time instant $t^n + \frac{1}{2}\Delta t_\ell$ for every $\ell \ne L$.
This approximation is in the same time instant as the solution in the next finer level. The use of this approximation consists in predicting the values of the virtual leaves at level $\ell+1$, at the proper time instant, before performing the flux computations of the leaves at level $\ell+1$. This approach is illustrated for the situation presented in Figure~\ref{fig:Rk1SyncIssue}(a).

\par In order to simplify the algorithm, the prediction of the virtual leaves at the coarser scales at the proper time instant, necessary to update the next finer level, is performed during the tree refreshing process, explained in detail for every RK step, in Section~\ref{syncChallengeTree}.
In this Section, we focus on the update of the virtual leaves at level $\ell+1$ in order to perform the flux computations on leaves of level $\ell$, as presented in Figure~\ref{fig:Rk1SyncIssue}(b).

\par Initially, the flux computations for the second RK step are performed on the leaves of level $L$. This choice is due to the fact that at this level, there are no interfaces of the type leaf/internal node, avoiding thus the situation shown in Figure~\ref{fig:Rk1SyncIssue}(b). 

\par Having updated the level $L$, the same process is repeated recursively for the levels $L-1$ down to $L=\ell_{\min}$.
However, due to the leaf/internal node scenario, the leaves and virtual leaves of the previously updated level $\ell+1$ must be synchronized at the instant $t^n + \Delta t_\ell$, which is equivalent to the instant $t^n + 2\Delta t_{\ell+1}$, in order to perform the flux computations of the next coarser level.
For that we propose the use of an extrapolated value, resulting in the situation presented in Figure~\ref{fig:Rk1SyncIssue}(b).

\par In this work, we propose to obtain an extrapolation at time instant $t^n + 2\Delta t_{\ell+1}$, which allows to compute the fluxes for the second step of the RK method. 
For this, the value $k_2$ is used as the value $k_1$ in a second RK1 time evolution. The first order approximation of the solution at the instant $t^n + 2\Delta t_{\ell+1}$ is computed as:

\begin{equation}
\bar{\mathbf{q}}_{\ell+1}\left(t^n + 2\Delta t_{\ell+1} \right) = 
\bar{\mathbf{q}}^*_{\ell+1} + 
\Delta t_{\ell+1}f\left(t+ \Delta t_{\ell+1};\; \bar{\mathbf{q}}^*_{\ell+1}\right) = \bar{\mathbf{q}}^n_{\ell+1} + k_1 + k_2
\label{eq:projleaf}
\end{equation}
where $\bar{\mathbf{q}}^*_{\ell+1}$ is the first RK step, and $\Delta t_{\ell+1}f\left(t+ \Delta t_{\ell+1};\; \bar{\mathbf{q}}^*_{\ell+1}\right)$ is the flux of the second RK step.
Note that this approximation is only possible due to the choice of the coefficients $a_{11} = b_1 = 1$.

\par Once the leaves of level $\ell+1$ are extrapolated to the instant $t^n + 2\Delta t_{\ell+1}$, the virtual leaves of this level are also updated. However, in order to predict the value of the virtual leaves at level $\ell+1$, the values of the virtual leaves and internal nodes of level $\ell$ must be synchronized at instant $t^n + \Delta t_{\ell}$ first.

\par The solution of the virtual leaves and internal nodes of level $\ell$ in this time instant is obtained during the tree refreshing process. However, in order to improve the solution of the internal nodes, their values at time instant $t^n + \Delta t_\ell$ are updated after the evolution of the leaves of level $\ell+1$.
This update is performed via projection of the extrapolated solution at instant $t^n + 2\Delta t_{\ell+1}$ from the level $\ell+1$ to $\ell$. 
In contrast to the leaves, which have an approximation at this instant, the internal nodes do not. The values for the internal nodes are obtained via linear extrapolation:
\begin{equation}
\bar{\mathbf{q}}_{\ell+1}\left(t^n + 2\Delta t_{\ell+1} \right) = 2\bar{\mathbf{q}}^{*}_{\ell+1} - \bar{\mathbf{q}}^n_{\ell+1}. 
\label{eq:projnode}
\end{equation}

\par This approximation is used to perform the projection of the solution of the internal node to the next finer level $\ell$. 
The projection procedure is given in Algorithm \ref{alg:projRK2}. 
After the projection procedure, the virtual leaves of level $\ell+1$ have their predicted values. Then, the flux computations and the update of the leaves at level $\ell$ can be performed. The second step of the RK method is given in Algorithm \ref{Steprk2}.

\par During the second RK step, the fluxes are also used for computing the second step of the NERK approximation with second order at the time instants $t^n+\frac{1}{2}\Delta t_\ell$ for RK2 or $t^n+\frac{1}{4}\Delta t_\ell$ and $t^n+\frac{3}{4}\Delta t_\ell$ for RK3 time evolution. 

As stated before, the NERK approximation for the RK2 time evolution is used to solve the synchronization according the time evolution problem given in Section \ref{LTissue2}. For the RK3 time evolution, the NERK approximations are used to perform the third RK step.

\par
\begin{figure}[htb]
\begin{center}
\begin{tabular}{lp{1cm}l}
\includegraphics[scale=0.5]{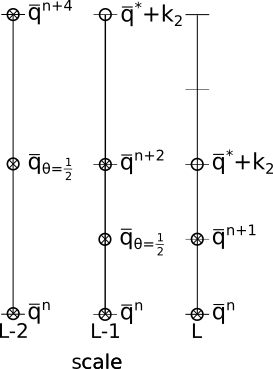} & &  \includegraphics[scale=0.5]{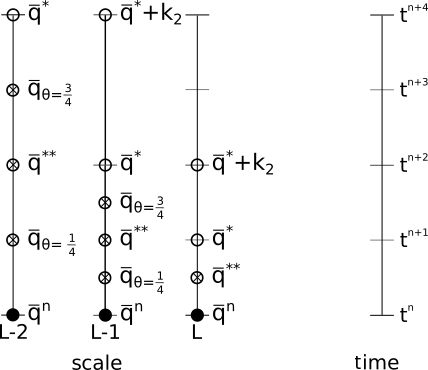}\\
a) RK2 evolution. &  & b) RK3 evolution.
\end{tabular}
\end{center}
\caption{Scheme of the obtained values and its respective instants for cells in each level after performing the second RK step in the LT method. a) for RK2 and b) for the RK3 time evolution. The filled, crossed and clear circles represent an available solution, of first, second and third order, respectively, at the corresponding time instant. }
\label{fig:NERKRK2-end}
\end{figure}

 After the second step of the RK evolution, a solution at the time instant $t^n + \Delta t_{\ell}$ for the RK2 method, and at $t^n + \frac{1}{2}\Delta t_{\ell}$ for the RK3 method is obtained. 
In order to prepare the tree structure for the next RK2 iteration and for the next RK3 step, another tree refreshing procedure must be applied. 
This procedure is given in detail in Section~\ref{syncChallengeTree}. 
After this process, the information obtained for the cells (internal nodes, leaves and virtual leaves) in every refinement level are illustrated in Figure \ref{fig:NERKRK2-end}.  

\subsubsection{Third Runge--Kutta step synchronization}\label{RKstep3}
The third step is done only in the case of the RK3 time evolution.
As for the second step of the RK evolution, in the third step, the evolution must be done first on the finest scale in order to obtain the necessary values to perform the evolution on the next coarser scale. 
However there are, indeed, slightly different synchronization issues in this case which we discuss in the following. 
\par As performed in the second RK step, the virtual leaves and the internal nodes have their values updated during the tree refreshing process, the details are discussed in the next section. 

\par To perform the third RK step, the fluxes of level $L$ are computed in the proper time instant $t^n + \frac{1}{2}\Delta t_L$. Then, the third RK step yields a third order solution at the time instant $t^n + \Delta t_L$.

\par Again as for the second step, this solution is projected onto the next coarser level in order to update the virtual leaves at level $L$ at the instant $t^n + \Delta t_L$, where the fluxes for the third step of the level $L-1$ shall be computed.

\par The projection procedure for the third RK step is, as mentioned above, slightly different from the projection procedure in the second step. 
In this case, the solution at the leaves at level $\ell$ after the third RK step matches with the time instant where the fluxes of level $\ell-1$ are computed. 
Then, the solution is projected from level $\ell$ to level $\ell-1$ by simple averaging, obtaining thus a solution at $t^n + \frac{1}{2}\Delta t_{\ell-1}$. 
This solution is used as an approximation for $\bar{\mathbf{q}}^{**}$. 
Then, the value of the internal nodes at the instant $t^n + \Delta t_{\ell-1}$ should be updated with a second order solution in order to continue the projections recursively during the third RK step. This update is performed using the following relation obtained from the NERK scheme:

\begin{equation}
\bar{\mathbf{q}}_{\ell}(t^n + \Delta t_{\ell}) = \bar{\mathbf{q}}_{\ell}^n +2 \left( \bar{\mathbf{q}}_{\ell, \theta = \frac{3}{4}} - \bar{\mathbf{q}}_{\ell, \theta = \frac{1}{4}} \right).
\label{rk3proj}
\end{equation}

\par The projection procedure inside the third RK step is given in Algorithm \ref{alg:projRK3}, and the execution of the third step of the RK evolution in Algorithm \ref{Steprk3}.

\subsection{Tree refreshing with time synchronization}\label{syncChallengeTree}
In the MR scheme, due to the same time step for every scale, the leaves of every scale store values corresponding to the same time instant. Therefore, projections of the leaves onto internal nodes receive values at the same time instant.
However, in the LT approach, due to fact that leaves of different scales are evolved with different time steps, the tree refreshing may project values at different time instants related to the leaves at the same scale. Accordingly, we need a procedure before each RK step internally to the time cycle, to project the solution to an internal node at its corresponding time instant. 

\par Moreover, these synchronizations are necessary to predict the values of the virtual leaves for the next steps of the RK method, at the required time instant. Only after that, the flux computations can be performed.

\par Before the first RK step, every scale of level $\ell_{\min}$ or greater, which must be evolved, has its solution at the same time instant. 
Therefore, the values of the virtual leaves can be predicted normally, except at level $\ell_{\min}$, which has its virtual leaves predicted using the solution at $t^n + \frac{1}{2}\Delta t_{\ell_{\min}-1}$ from the level $\ell_{\min}-1$.

 After the first step of the RK method, each leaf yields a first order solution at the time instant $t^n + \Delta t_\ell$, which corresponds to the time instant $t^n + \frac{1}{2}\Delta t_{\ell-1}$.
 This implies that the projection of the leaves from level $\ell$ to $\ell-1$ produces an approximation at that time instant. The value of an internal node at this instant is given by the mean value of its children:
\begin{equation}
\bar{\mathbf{q}}_{\ell-1}\left(t^n + \frac{1}{2}\Delta t_{\ell-1}\right) = \bar{\mathbf{q}}_{\ell-1}^n + \frac{1}{2}\Delta t_{\ell-1}f\left(t^n,\bar{\mathbf{q}}_{\ell-1}^n\right)= \frac{1}{2^\mathit{d}}\sum_{i=1}^{2^\mathit{d}}\left[\bar{\mathbf{q}}_{\ell, \; i}^n + \Delta t_{\ell}f\left(t^n,\bar{\mathbf{q}}_{\ell, \; i}^n\right) \right]
\end{equation}
where $\mathit{d}$ is the number of dimensions of the problem and $\bar{\mathbf{q}}_{\ell, \; i}^n$ is the solution of the children cell $i$. This value is stored in order to predict the values of the virtual leaves of level $\ell$ for the second RK step. 
Furthermore, this value must be extrapolated to the instant $t^n + \Delta t_{\ell-1}$, obtaining the value $\bar{\mathbf{q}}^{*}_{\ell-1}$. This value, used to update the virtual leaves of level $\ell$ during the second RK step, is obtained by linear extrapolation: 
\begin{equation}
\bar{\mathbf{q}}_{\ell-1}^{*} = 2\bar{\mathbf{q}}_{\ell-1}\left(t^n + \frac{1}{2}\Delta t_{\ell-1}\right) - \bar{\mathbf{q}}^{n}_{\ell-1}, 
\label{refreshRk2}
\end{equation}
where the value $\bar{\mathbf{q}}^{n}_{\ell-1}$ is the solution before the time evolution. 
This value should be stored for every node before the RK procedure. 
The tree refreshing procedure after the first RK step is given in Algorithm \ref{alg:refRK2}.

\par When performing the RK3 time evolution, before the third RK step, besides the solution, the values $\bar{\mathbf{q}}_{\theta = \frac{1}{4}}$ and $\bar{\mathbf{q}}_{\theta = \frac{3}{4}}$ must be obtained for the internal nodes. The solution in those instants are required to predict the solution of the virtual leaves in the next finer level at the instant required for this RK step.  

\par Those values are obtained through projections from the level $L$ down to the level $\ell_{\min}$, initially by projecting the solution $\bar{\mathbf{q}}_{\ell}^{**}$ at the time instant $t^n + \frac{1}{2}\Delta t_{\ell}$, obtaining thus a value at $t^n + \frac{1}{4}\Delta t_{\ell-1}$, which is used as an approximation for $\bar{\mathbf{q}}_{\ell-1, \;\theta = \frac{1}{4}}$.

\par Using this solution, the Equation \ref{rkcont} and the value $\bar{\mathbf{q}}_{\ell}^{*}$, computed during the projection before the second RK step, we obtain the following relation to reconstruct the solution inside the interval $\left[t^n; t^n+\Delta t^n\right]$:
\begin{equation}
\bar{\mathbf{q}}_{\ell}\left(t^n + \theta\Delta t_{\ell}\right) =  \left[1-\theta - 12\theta^2 \right]\bar{\mathbf{q}}^{n}_{\ell} +  \left[\theta - 4\theta^2\right]\bar{\mathbf{q}}^{*}_{\ell} + 16\theta^2\bar{\mathbf{q}}_{\ell, \;\theta = \frac{1}{4}}.
\end{equation}
This expression is based on the NERK scheme presented in Section \ref{NERK}. Using the value $\theta = \frac{3}{4}$, we have:
\begin{equation}
\bar{\mathbf{q}}_{\ell, \; \theta = \frac{3}{4}} = -\frac{13}{2}\bar{\mathbf{q}}^{n}_{\ell} -\frac{3}{2} \bar{\mathbf{q}}^{*}_{\ell} + 9\bar{\mathbf{q}}_{\ell, \;\theta = \frac{1}{4}}
\label{ProjRK3equations1}
\end{equation}
Subsequently, using these values, the value $\bar{\mathbf{q}}^{**}$ can be computed as:
\begin{equation}
\bar{\mathbf{q}}^{**}_{\ell} = -\frac{11}{2}\bar{\mathbf{q}}^{n}_{\ell} - \frac{3}{2} \bar{\mathbf{q}}^{*}_{\ell} + 8\bar{\mathbf{q}}_{\ell, \;\theta = \frac{1}{4}}
\label{ProjRK3equations2}
\end{equation}

\par The solution $\bar{\mathbf{q}}_\ell^{**}$ allows to continue the projection procedure recursively by obtaining the value $\bar{\mathbf{q}}_{\ell-1, \;\theta = \frac{1}{4}}$ at the next coarser level, while the solution $\bar{\mathbf{q}}_{\ell, \;\theta = \frac{3}{4}}$ is used to compute the value of the virtual leaves at level $\ell+1$ on the following iterations, as explained in Section \ref{LTSync1}. 

\par
\begin{figure}[htb]
\begin{center}
\includegraphics[scale=0.5]{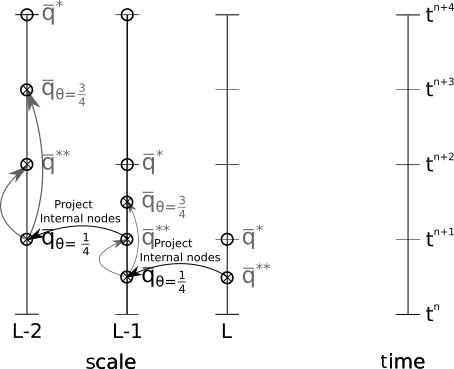}
\end{center}
\caption{Projection scheme before the third RK step.}
\label{projRK2}
\end{figure}

 The tree refreshing procedure, illustrated in Figure \ref{projRK2}, consists in projecting the solution $\bar{\mathbf{q}}^{**}_{L}$ from the level $L$ onto the level $L-1$. This solution is approximated as $\bar{\mathbf{q}}_{L-1, \;\theta = \frac{1}{4}}$. Then, using the previously obtained value $\bar{\mathbf{q}}^{*}_{L-1}$, the values $\bar{\mathbf{q}}_{L-1, \; \theta = \frac{3}{4}}$ and $\bar{\mathbf{q}}^{**}_{L-1}$ can be obtained using Equations (\ref{ProjRK3equations1}) and (\ref{ProjRK3equations2}). This procedure is recursively repeated down to the coarsest scale $\ell_{\min}$.

\par In the third RK step, the fluxes are computed at the intermediary time instant $t^n + \frac{1}{2}\Delta t_\ell$. 
For an adjacent coarser leaf, this instant is $t^n + \frac{1}{4}\Delta t_{\ell-1}$ or $t^n + \frac{3}{4}\Delta t_{\ell-1}$, depending on the current iteration number. 
Basically, if the current scale is $\ell_{\min}$ then the NERK value used to compute the virtual leaves is available at time instant $t^n + \frac{3}{4}\Delta t_{\ell_{\min}-1}$. 
Otherwise, the value at the time instant $t^n + \frac{1}{4}\Delta t_{\ell-1}$ is used.

\par
\begin{figure}[htb]
\begin{center}
\includegraphics[scale=0.5]{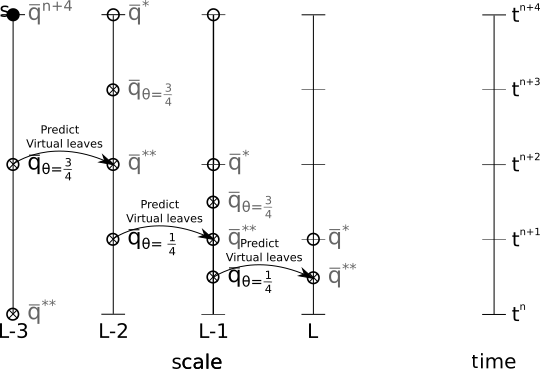}
\end{center}
\caption{Scheme of the choice of the NERK solution to predict the virtual leaves for the third RK step in a LT iteration with $\ell_{\min} = L-2$. The solution $\bar{\mathbf{q}}_{\theta = \frac{3}{4}}$ is used to predict the $\bar{\mathbf{q}}^{**}$ values at level $\ell_{\min}$. For the other scales, the solution $\bar{\mathbf{q}}^{**}$ is predicted with the values $\bar{\mathbf{q}}_{\theta = \frac{1}{4}}$.}
\label{fig:NERKRK3}
\end{figure}

 The choice which value of the NERK approximation shall be used in the flux computations is shown in Figure \ref{fig:NERKRK3}, where we choose $\ell_{\min}=L\!-\!2$.

\par Once the projections have been performed from the level $L$ until the level $\ell_{\min}$, the virtual leaves have their values $\bar{\mathbf{q}^{**}}_\ell$ predicted using the values $\bar{\mathbf{q}}_{\ell-1, \;\theta = \frac{3}{4}}$, in case of $\ell=\ell_{\min}$; or using $\bar{\mathbf{q}}_{\ell, \;\theta = \frac{1}{4}}$, otherwise.

\par Then, the values $\bar{\mathbf{q}}_{\ell, \;\theta = \frac{1}{4}}$ and $\bar{\mathbf{q}}_{\ell, \;\theta = \frac{3}{4}}$ for these virtual leaves are obtained using the Equations~(\ref{ProjRK3equations1}) and (\ref{ProjRK3equations2}).
Those values are required here to predict the values of the virtual leaves at the next finer scale.

\par The tree refreshing procedure before the third RK step is detailed in Algorithm \ref{alg:refRK22}.

\subsection{Synchronization after the time evolution}\label{LTissue2}

After finishing the RK steps, the leaves at level $\ell$ perform a complete time evolution with a time step that is twice the time step of the leaves at level $\ell+1$. In order to perform the next time evolution at scale $\ell$, the finer scale leaves must reach the same time instant. For that, a second time evolution is required for this finer scale leaf. 
\par Moreover, in the next iteration of the time evolution, a new value for $\ell_{\min}$ is obtained. To predict the virtual leaves at the new scale $\ell_{\min}$ at the proper time instant, the values $\bar{\mathbf{q}}_{\ell_{\min}-1,\; \theta = \frac{1}{2}}$ for RK2 or $\bar{\mathbf{q}}_{\ell_{\min}-1}^{**}$ for RK3 are used. For the scales finer than $\ell_{\min}$, the virtual leaves are predicted with the MR method.

\par To solve this synchronization challenge and to maintain the approximation order, a high order approximation at time step $t^n + \frac{1}{2}\Delta t_\ell$, where $\ell$ is the refinement level of the leaf, is required. 
In this case, in the tree terminology, it means that for the coarser leaves, an approximation of its virtual children at the instant $t^n + \frac{1}{2}\Delta t_\ell$ is required. 
\par To obtain these values, it is necessary to obtain approximations at the instant $t^n + \frac{1}{2}\Delta t_\ell$ during the Runge--Kutta evolution for every cell of refinement level $\ell \ne L$ and its neighbours, including internal nodes. Then, the values of the virtual children are obtained by the prediction procedure.
\par For the RK2 time evolution, we use the NERK approach to obtain the approximation with $\theta = \frac{1}{2}$. However, the implementation of the MRLT/NERK approach for the RK3 time evolution requires extra memory. Hence we use the approximation $\bar{\mathbf{q}}^{**}$, obtained after the second step of the compact RK, which is at the instant $t^n + \frac{1}{2}\Delta t_\ell$, instead of the NERK approximation with $\theta = \frac{1}{2}$ in order to reduce the number of variables to be stored in the problem.

\par Besides those values, the updated solution, at time instant $t^n + \Delta t_{\ell}$, is projected onto the coarser level at the time instant $t^n + \frac{1}{2}\Delta t_{\ell-1}$ by simple averaging. This solution is then used as an approximation for $\bar{\mathbf{q}}_{\ell-1, \; \theta = \frac{1}{2}}$, in the RK2 time evolution, or $\bar{\mathbf{q}}^{**}_{\ell-1}$ for the RK3 time evolution. Using the NERK equations given in Section \ref{CompactFormulation} and the current value of the internal node at $t^n + \Delta t_{\ell-1}$ (RK first step), the solution at instant $t^n + \Delta t_{\ell}$ is obtained by the following relation: 
\begin{equation}
\bar{\mathbf{q}}_{\ell}\left(t^n + \Delta t_{\ell}\right) = -2\bar{\mathbf{q}}^{n}_{\ell} - \bar{\mathbf{q}}^*_{\ell} + 4\bar{\mathbf{q}}_{\ell, \;\theta = \frac{1}{2}}.
\end{equation}
\par This value is used to compute the value of the virtual leaves at level $\ell_{\min}+1$ in the following iterations at the second RK step, as explained in Section \ref{RKstep2}.

\subsection*{Adapting the grid during local time-stepping}\label{AdaptTree}
During the LT scheme, in most parts of the numerical simulation, some scales have results at different time instants. Therefore, in order to avoid errors to be caused by converting a cell with a solution value at a determined time instant to another scale where its cells are in a different time instant, there is a need for a criterion when adapting the solution during the LT scheme.
\par This criterion consists in combining or splitting only leaves from a scale to another scale whose leaves are at the same time instant. These scales can be easily detected due to the fact that the time evolution procedure is applied only in scales at the same time instant. In other words, to split a cell into more refined ones, both scales must be included into the current time evolution iteration. The same is valid to combine cells into a coarser one.
\par Moreover, beyond this restriction, the remeshing process is identical to the remeshing process of the multiresolution method.


\subsection{Some remarks on the convergence and stability of MRLT/NERK schemes}
\label{sec:stability}
In this section, we perform a stability analysis to check the convergence order of the MRLT/NERK3 method. For that, we compute the interpolation and extrapolation errors in each approximation required for the method, considering the interface between cell of different levels.

To analyze the stability typically a simple linear ODE,
\begin{equation}
\frac{dq}{dt} \, = \, \lambda u
\end{equation}
where a constant complex-valued coefficient $\lambda \in \mathbb{C}$ with negative or zero real part is considered. The above equation is completed with an initial condition.
Using Fourier analysis of linear PDEs it can be shown that imaginary values of $\lambda$ correspond to pure advection, while real (negative) values correspond to pure diffusion equations.

One-step schemes, including Runge--Kutta methods can then we written in the following form,
\begin{equation}
\hat{\mathbf{q}}^{n+1} \, = \, g(\lambda \Delta t) \, \hat{\mathbf{q}}^n
\end{equation}
where $g$ is a polynomial. In particular, for RK1, RK2 and RK3 we respectively have $g(\lambda \Delta t) = 1 + \lambda \Delta t$, $g(\lambda \Delta t) = 1 + \lambda \Delta t + \frac{1}{2} \lambda^2 \Delta t^2$ and $g(\lambda \Delta t) = 1 + \lambda \Delta t + \frac{1}{2} \lambda^2 \Delta t^2 + \frac{1}{6} \lambda^3 \Delta t^3$, which correspond to the truncated Taylor series of $\exp(\lambda \Delta t)$. A method is called of convergence order $n$ if its polynomial $g(\lambda \Delta t)$ reconstructs the Taylor series until the $\lambda^n \Delta t^n$ term.
 
\par To compute the polynomial $g$ for the MRLT/NERK methods, the analysis is performed by computing and inserting each error $\epsilon_i$ obtained in every approximation required to perform the MRLT method. 
\par The errors obtained in those approximations are inserted into the flux computations of the stability model as follows:
\begin{equation}
f(\bar{\mathbf{q}}) = f(\hat{\mathbf{q}} - \epsilon) \, ,
\end{equation}
where $\bar{\mathbf{q}}$ is an approximation obtained during the MRLT/NERK method and $\hat{\mathbf{q}}$ is the solution that would be obtained by the regular RK method.

\par In the following, we obtain the errors $\epsilon$ for each approximation used for the first RK evolution on both finer and coarser scales. In this evolution, both scales are initially at the same time instant, so there is no approximation error for the first RK step ($k_1$). In other words, it means that $\bar{\mathbf{q}}^n = \hat{\mathbf{q}}^n$. 

\par In this section, we consider the RK3 scheme as given in Section \ref{CompactFormulation}. 

Considering that the solution $\bar{\mathbf{q}}^n$ is at the same instant for every scale to be updated in the time evolution, the first RK stage reads,
\begin{equation}
\bar{\mathbf{q}}^* = \bar{\mathbf{q}}^n +\Delta t^n f(t^n,\bar{\mathbf{q}}^n)
\end{equation}
and yields no interpolation errors due to the LT approach for both scales. That leads to $\bar{\mathbf{q}}^* = \hat{\mathbf{q}}^*$.

\subsubsection*{Approximations after the first RK step}
Here we compute the errors obtained from the predictions and projections after the first RK step, as performed in Section \ref{RKstep2}.
The solution of the intermediary solution $\bar{\mathbf{q}}^*$ is at a different time instant for each refinement level, requiring thus a series of interpolations and extrapolations, which introduces errors $\epsilon_1$ and $\epsilon_2$ into the scheme. 
\begin{itemize}
\item $\epsilon_1$: Prediction error in the finer leaf.
\par Error obtained in the prediction of $\bar{\mathbf{q}}^*$ in the finer leaf using the values of the solution in the lower level with half of the time-step.
\begin{equation}
\epsilon_1 = \hat{\mathbf{q}}^{*}_\ell - P_{\ell-1 \rightarrow \ell}\left[\bar{\mathbf{q}}^n_{\ell-1} - \frac{\Delta t_{\ell-1}}{2}f(\bar{\mathbf{q}^n}_{\ell-1})\right] = 0.
\end{equation}
\par \textbf{Proof:}
Considering $f(\bar{\mathbf{q}^n}_{\ell-1}) = \lambda \bar{\mathbf{q}^n}_{\ell-1}$
\begin{eqnarray}
\epsilon_1 = \hat{\mathbf{q}}^{n}_\ell + \Delta t_\ell \lambda\hat{\mathbf{q}}^{n}_\ell -P_{\ell-1 \rightarrow \ell}\left[ \bar{\mathbf{q}}_{\ell-1} - \frac{\Delta t_{\ell-1}}{2}\lambda\bar{\mathbf{q}}_{\ell-1}\right]\\
\epsilon_1 = \hat{\mathbf{q}}^{n}_\ell + \Delta t_\ell \lambda \hat{\mathbf{q}}^{n}_\ell - \bar{\mathbf{q}}_{\ell} - \frac{\Delta t_{\ell-1}}{2}\lambda\bar{\mathbf{q}}_{\ell}\\
\epsilon_1 = 0
\end{eqnarray}

\item $\epsilon_2$: Projection error in the coarser leaf. 
Error obtained in the projection to obtain $\bar{\mathbf{q}}^*$ in the coarser level. This solution is a first order extrapolation to the time instant $t^n+\Delta t_{\ell-1}$ for the leaves in the finer level.
\begin{equation}
\epsilon_2 = \hat{\mathbf{q}}^{*}_\ell - P_{\ell+1\rightarrow \ell}\left[ \underbrace{\bar{\mathbf{q}}^n_{\ell+1} + \Delta t_{\ell+1} f(\bar{\mathbf{q}}^n_{\ell+1})}_{\text{First RK step in finer leaf}} + \underbrace{\Delta t_{\ell+1} f(\bar{\mathbf{q}}^{*}_{\ell+1})}_{\text{extrapolation to } t^n+\Delta t_{\ell-1}}\right] = -\frac{\Delta t^2_{\ell}\lambda^2 }{4} \hat{\mathbf{q}}^{n}_{\ell}
\end{equation}
\par \textbf{Proof:}
Here we consider the errors caused in the previous approximations in the flux terms. Using $\mathbf{q}^* = \mathbf{q} + \Delta tf(\mathbf{q}) -\epsilon_1$ inside the flux term:
\begin{eqnarray}
\epsilon_2 = \hat{\mathbf{q}}^{n}_\ell + \Delta t_\ell f(\hat{\mathbf{q}}^{n}_\ell) - P_{\ell+1\rightarrow \ell}\left[ \bar{\mathbf{q}}^n_{\ell+1} + \Delta t_{\ell+1} f(\bar{\mathbf{q}}^n_{\ell+1}) + \Delta t_{\ell+1} f(\hat{\mathbf{q}}^{n}_{\ell+1} + \Delta t_{\ell+1} f(\hat{\mathbf{q}}^{n}_{\ell+1}) - \epsilon_1)\right] 
\end{eqnarray}
Considering $f(\bar{\mathbf{q}^n}_{\ell+1}) = \lambda \bar{\mathbf{q}^n}_{\ell+1}$ we have
\begin{eqnarray}
\epsilon_2 = \hat{\mathbf{q}}^{n}_\ell + \Delta t_\ell \lambda\hat{\mathbf{q}}^{n}_\ell - P_{\ell+1\rightarrow \ell}\left[ \hat{\mathbf{q}}^n_{\ell+1} + \Delta t_{\ell+1} \lambda \hat{\mathbf{q}}^n_{\ell+1} + \Delta t_{\ell+1} \lambda( \hat{\mathbf{q}}^{n}_{\ell+1} + \Delta t_{\ell+1} \lambda(\hat{\mathbf{q}}^{n}_{\ell+1}) )\right]
\end{eqnarray}
Applying the projection operator we get
\begin{eqnarray}
\epsilon_2 = \hat{\mathbf{q}}^{n}_\ell + \Delta t_\ell \lambda\hat{\mathbf{q}}^{n}_\ell - \hat{\mathbf{q}}^n_{\ell} - \Delta t_{\ell+1} \lambda \hat{\mathbf{q}}^n_{\ell} - \Delta t_{\ell+1} \lambda\left( \hat{\mathbf{q}}^{n}_{\ell} + \Delta t_{\ell+1} \lambda \hat{\mathbf{q}}^{n}_{\ell} \right)\\
\epsilon_2 = \Delta t_\ell \lambda\hat{\mathbf{q}}^{n}_\ell - \frac{\Delta t_{\ell}}{2} \lambda \hat{\mathbf{q}}^n_{\ell} - \frac{\Delta t_{\ell}}{2} \lambda\left( \hat{\mathbf{q}}^{n}_{\ell} + \frac{\Delta t_{\ell}}{2} \lambda \hat{\mathbf{q}}^{n}_{\ell} \right)\\
\epsilon_2 = -\frac{\Delta t^2_{\ell}\lambda^2 }{4} \hat{\mathbf{q}}^{n}_{\ell}
\end{eqnarray}
\end{itemize}
Once the interpolation errors are obtained, the second Runge--Kutta step is performed for the finer scale:
\begin{eqnarray}
\bar{\mathbf{q}}^{**} =  \frac{3}{4}\bar{\mathbf{q}}^n + \frac{1}{4}\bar{\mathbf{q}}^* + \frac{1}{4}\Delta t^n f\left(t^n + \Delta t^n, \bar{\mathbf{q}}^*\right)\\
\bar{\mathbf{q}}^{**} =  \frac{3}{4}\hat{\mathbf{q}}^n + \frac{1}{4}\hat{\mathbf{q}}^* + \frac{1}{4}\Delta t^n f\left(\hat{\mathbf{q}}^* - \epsilon_1\right)\\
\bar{\mathbf{q}}^{**} =  \frac{3}{4}\hat{\mathbf{q}}^n + \frac{1}{4}\hat{\mathbf{q}}^* + \frac{1}{4}\Delta t^n f\left(\hat{\mathbf{q}}^*\right)\\
\bar{\mathbf{q}}^{**} =  \hat{\mathbf{q}}^{**}
\end{eqnarray}
We observe that the $\bar{\mathbf{q}}^{**}$ solution does not have any interpolation errors due to the LT scheme. Then the second order RK step is performed for the coarser scale:
\begin{eqnarray}
\bar{\mathbf{q}}^{**} =  \frac{3}{4}\bar{\mathbf{q}}^n + \frac{1}{4}\bar{\mathbf{q}}^* + \frac{1}{4}\Delta t^n f\left(t^n + \Delta t^n, \bar{\mathbf{q}}^*\right)\\
\bar{\mathbf{q}}^{**} =  \frac{3}{4}\hat{\mathbf{q}}^n + \frac{1}{4}\hat{\mathbf{q}}^* + \frac{1}{4}\Delta t^n f\left(\hat{\mathbf{q}}^* - \epsilon_2\right)\\
\bar{\mathbf{q}}^{**} =  \frac{3}{4}\hat{\mathbf{q}}^n + \frac{1}{4}\hat{\mathbf{q}}^* + \frac{1}{4}\Delta t^n \lambda \hat{\mathbf{q}}^* - \frac{1}{4}\Delta t^n \lambda\epsilon_2\\
\bar{\mathbf{q}}^{**} =  \frac{3}{4}\hat{\mathbf{q}}^n + \frac{1}{4}\hat{\mathbf{q}}^* + \frac{z}{4} \hat{\mathbf{q}}^* - \frac{\Delta t_\ell \lambda}{4}\epsilon_2\\
\bar{\mathbf{q}}^{**} =  \hat{\mathbf{q}}^{**} -  \frac{\Delta t_\ell \lambda}{4}\epsilon_2 \label{Eq:q**}
\end{eqnarray}
Here, the approximation has an error of $-\frac{\Delta t_\ell \lambda}{4}\epsilon_2$.

\subsubsection*{Approximations after the second RK step}
In the following, we compute the errors obtained from the predictions and projections after the second RK step, as performed in Section \ref{RKstep3}.
The approximations here are the NERK solution at $t^n + \frac{1}{4}\Delta t_{\ell-1}$ to predict the solution $\bar{\mathbf{q}}^{**}_\ell$ on finer scales and the projections of the RK3 evolution of the finer leaves to approximate $\bar{\mathbf{q}}^{**}_{\ell-1}$. The first one is given by:
\begin{itemize}
\item $\epsilon_3$: Prediction error using NERK approximation
\begin{eqnarray}
\epsilon_3 = \hat{\mathbf{q}}^{**}_\ell - P_{\ell-1\rightarrow \ell}\left[ \bar{\mathbf{q}}^n_{\ell-1,\theta = \frac{1}{4}} \right] =  \left(\frac{\Delta t_\ell^2 \lambda^2}{8} - \frac{\Delta t_\ell^3\lambda^3}{16} \right)\hat{\mathbf{q}}^n_\ell 
\end{eqnarray}
\par \textbf{Proof:}
Considering the NERK equation for $\theta = \dfrac{1}{4}$ given in Section \ref{CompactFormulation} we get,
\begin{eqnarray}
\epsilon_3 = \hat{\mathbf{q}}^{**}_\ell - P_{\ell-1\rightarrow \ell}\left[ \bar{\mathbf{q}}^n_{\ell-1} + \frac{7}{32}\Delta t_{\ell-1}f(\bar{\mathbf{q}}^{n}_{\ell-1}) + \frac{1}{32}\Delta t_{\ell-1}f(\bar{\mathbf{q}}^{*}_{\ell-1})\right]
\end{eqnarray}
Then, we consider the error $\epsilon_2$ in the flux computations for the coarser scale:
\begin{eqnarray}
\epsilon_3 = \hat{\mathbf{q}}^n_\ell\left(1 + \frac{\Delta t_\ell \lambda}{2} +  \frac{\Delta t_\ell^2 \lambda^2}{4}\right)  - P_{\ell-1\rightarrow \ell}\left[ \bar{\mathbf{q}}^n_{\ell-1} + \frac{7}{32}\Delta t_{\ell-1}\lambda\bar{\mathbf{q}}^{n}_{\ell-1} + \frac{1}{32}\Delta t_{\ell-1}\lambda\left(\hat{\mathbf{q}}^{*}_{\ell-1}- \epsilon_2\right)\right]\\
\epsilon_3 = \hat{\mathbf{q}}^n_\ell\left(\frac{\Delta t_\ell\lambda}{16} +  \frac{\Delta t_\ell^2 \lambda^2}{4} \right) - \frac{\Delta t_\ell\lambda}{16}P_{\ell-1\rightarrow \ell}\left[\hat{\mathbf{q}}^{*}_{\ell-1}- \epsilon_2\right]\\
\epsilon_3 = \hat{\mathbf{q}}^n_\ell\left(\frac{\Delta t_\ell \lambda}{16} +  \frac{\Delta t_\ell^2\lambda^2}{4}\right) - \frac{\Delta t_\ell \lambda}{16}P_{\ell-1\rightarrow \ell}\left[\hat{\mathbf{q}}^{n}_{\ell-1} +\Delta t_{\ell-1}f(\hat{\mathbf{q}}^{n}_{\ell-1}) + \frac{\Delta t^2_{\ell-1}\lambda^2 }{4} \hat{\mathbf{q}}^{n}_{\ell-1} \right]\\
\epsilon_3 = \frac{\Delta t_\ell^2\lambda^2}{4}\hat{\mathbf{q}}^n_\ell - \frac{\Delta t_\ell^2\lambda^2}{8} \hat{\mathbf{q}}^{n}_\ell - \frac{\Delta t_\ell\lambda}{16}P_{\ell-1\rightarrow \ell}\left[ \frac{\Delta t^2_{\ell-1}\lambda^2 }{4} \hat{\mathbf{q}}^{n}_{\ell-1} \right]\\
\epsilon_3 = \left(\frac{\Delta t_\ell^2\lambda^2}{8} - \frac{\Delta t_\ell^3\lambda^3}{16} \right)\hat{\mathbf{q}}^n_\ell 
\end{eqnarray}
\end{itemize}
Afterwards we perform the third RK step for the finer leaf.
\begin{eqnarray}
\bar{\mathbf{q}}^{n+1} =  \frac{1}{3}\bar{\mathbf{q}}^n+ \frac{2}{3}\bar{\mathbf{q}}^{**}  + \frac{2}{3}\Delta t^nf\Big(t^n+\frac{1}{2}\Delta t^n, \bar{\mathbf{q}}^{**}\Big)\\
\bar{\mathbf{q}}^{n+1} =  \frac{1}{3}\hat{\mathbf{q}}^n+ \frac{2}{3}\hat{\mathbf{q}}^{**}  + \frac{2z}{3}\Big(\hat{\mathbf{q}}^{**}-\epsilon_3\Big)\\
\bar{\mathbf{q}}^{n+1} =  \frac{1}{3}\hat{\mathbf{q}}^n+ \frac{2}{3}\left(1 + \frac{\Delta t_\ell \lambda}{2} +  \frac{\Delta t_\ell^2\lambda^2}{4}\right)\hat{\mathbf{q}}^n  + \frac{2\Delta t_\ell\lambda}{3}\left[\left(1 + \frac{\Delta t_\ell\lambda}{2} +  \frac{\Delta t_\ell^2\lambda^2}{4}\right)\hat{\mathbf{q}}^n-\epsilon_3\right]\\
\bar{\mathbf{q}}^{n+1} =  \left(1 + \Delta t_\ell \lambda +  \frac{\Delta t_\ell^2\lambda^2}{2} +  \frac{\Delta t_\ell^3\lambda^3}{6} \right)\hat{\mathbf{q}}^n  -\frac{2\Delta t_\ell \lambda}{3}\left(\frac{\Delta t_\ell^2\lambda^2}{8} - \frac{\Delta t_\ell^3\lambda^3}{16} \right)\hat{\mathbf{q}}^n \\
\bar{\mathbf{q}}^{n+1} =  \left(1 + \Delta t_\ell \lambda +  \frac{\Delta t_\ell^2\lambda^2}{2} +  \frac{\Delta t_\ell^3\lambda^3}{12} +\frac{\Delta t_\ell^4\lambda^4}{24} \right)\hat{\mathbf{q}}^n\label{EvRk3finer}
\end{eqnarray}
This shows that the third order accuracy is lost in the evolution of the finer leaf.

\begin{itemize}
\item $\epsilon_4$: Projection error using finer leaf evolution.
\begin{equation}
\epsilon_4 = \hat{\mathbf{q}}^{**}_\ell - P_{\ell+1\rightarrow \ell}\left[ \bar{\mathbf{q}}^{n+1}_{\ell+1} \right] = \left(\frac{\Delta t_\ell^2 \lambda^2}{8} -  \frac{\Delta t_\ell^3 \lambda^3}{96} -\frac{\Delta t_\ell^4 \lambda^4}{384} \right)\hat{\mathbf{q}}^n_{\ell}
\end{equation}

\par \textbf{Proof:} Using Equation \ref{EvRk3finer} to represent the time evolution for the finer leaf, we have,
\begin{eqnarray}
\epsilon_4 = \hat{\mathbf{q}}^{**}_\ell - P_{\ell+1\rightarrow \ell}\left[ \left(1 + \Delta t_{\ell+1}\lambda +  \frac{\left(\Delta t_{\ell+1}\lambda\right)^2}{2} +  \frac{\left(\Delta t_{\ell+1}\lambda\right)^3}{12} +\frac{\left(\Delta t_{\ell+1}\lambda\right)^4}{24} \right)\hat{\mathbf{q}}^n_{\ell+1} \right]\\
\epsilon_4 = \hat{\mathbf{q}}^{**}_\ell - P_{\ell+1\rightarrow \ell}\left[ \left(1 + \frac{\Delta t_{\ell}\lambda}{2} +  \frac{\left(\Delta t_\ell\lambda\right)^2}{8} +  \frac{\left(\Delta t_\ell\lambda\right)^3}{96} +\frac{\left(\Delta t_\ell\lambda\right)^4}{384} \right)\hat{\mathbf{q}}^n_{\ell+1} \right]\\
\epsilon_4 = \left(1 + \frac{\Delta t_\ell \lambda}{2} +  \frac{\Delta t_\ell^2\lambda^2}{4}\right)\hat{\mathbf{q}}^n_\ell - \left(1 + \frac{\Delta t_\ell \lambda}{2} +  \frac{\Delta t_\ell^2\lambda^2}{8} +  \frac{\Delta t_\ell^3\lambda^3}{96} +\frac{\Delta t_\ell^4\lambda^4}{384} \right)\hat{\mathbf{q}}^n_{\ell}\\
\epsilon_4 = \left(\frac{\Delta t_\ell^2\lambda^2}{8} -  \frac{\Delta t_\ell^3\lambda^3}{96} -\frac{\Delta t_\ell^4\lambda^4}{384} \right)\hat{\mathbf{q}}^n_{\ell}
\end{eqnarray}

\end{itemize}
And finally, performing the third RK step for the coarser leaf, we get
\begin{eqnarray}
\bar{\mathbf{q}}^{n+1} =  \frac{1}{3}\bar{\mathbf{q}}^n+ \frac{2}{3}\bar{\mathbf{q}}^{**}  + \frac{2}{3}\Delta t_\ell f\left[t^n+\frac{1}{2}\Delta t_\ell, \bar{\mathbf{q}}^{**}\right]\\
\bar{\mathbf{q}}^{n+1} =  \frac{1}{3}\hat{\mathbf{q}}^n+ \frac{2}{3}\left(\hat{\mathbf{q}}^{**}-  \frac{\Delta t_\ell \lambda}{4}\epsilon_2 \right)  + \frac{2\Delta t_\ell \lambda}{3}\left[\hat{\mathbf{q}}^{**}-\epsilon_4\right]\\
\bar{\mathbf{q}}^{n+1} = \frac{1}{3}\hat{\mathbf{q}}^n+ \frac{2}{3}\left(\hat{\mathbf{q}}^{**} +  \frac{\Delta t^3_{\ell}\lambda^3 }{16} \hat{\mathbf{q}}^{n}_{\ell} \right)  + \frac{2\Delta t_\ell \lambda}{3}\left[\hat{\mathbf{q}}^{**}-\left(\frac{\Delta t_\ell^2 \lambda^2}{8} -  \frac{\Delta t_\ell^3 \lambda^3}{96} -\frac{\Delta t_\ell^4 \lambda^4}{384} \right)\hat{\mathbf{q}}^n_{\ell}\right]\\
\bar{\mathbf{q}}^{n+1} = \hat{\mathbf{q}}^{n+1} + \frac{\Delta t^3_{\ell}\lambda^3 }{24} \hat{\mathbf{q}}^{n}_{\ell} + \left(-\frac{\Delta t_\ell^3 \lambda^3}{12} +  \frac{\Delta t_\ell^4 \lambda^4}{144} +\frac{\Delta t_\ell^5 \lambda^5}{576} \right)\hat{\mathbf{q}}^n_{\ell}\\
\bar{\mathbf{q}}^{n+1} = \left(1+\Delta t_\ell \lambda+\frac{\Delta t_\ell^2 \lambda^2}{2} + \frac{\Delta t_\ell^3 \lambda^3}{8}+ \frac{\Delta t_\ell^4 \lambda^4}{144}+ \frac{\Delta t_\ell^5 \lambda^5}{576} \right)\hat{\mathbf{q}}^{n}. \label{EvRk3coarser}
\end{eqnarray}

This result also shows a loss of the third order accuracy, since the third order term of the Taylor series $\frac{\Delta t_\ell^3 \lambda^3}{6}$ is not found.

\subsection{Discussion on the numerical convergence}
To study the numerical convergence of local time stepping we consider the advection equation
\begin{equation}
\frac{\partial Q}{\partial t} +  \frac{\partial Q}{\partial x}=0 \qquad x\in[0,1]
\end{equation}
with periodic boundary conditions and a Gaussian initial condition, given by
$Q(x) = \exp{\left(-100\left(x-0.25\right)^2\right)}$. The simulation is performed for one time cycle. 
In order to avoid errors due to the MR scheme, we performed the adaptive simulations using two fixed grids in the domain. The first one in the interval $[0,0.5]$ with 
$\Delta x = 1/512$, and the second in the interval $[0.5,1.0]$ with $\Delta x = 1/256$, are both fixed during the entire simulation.
We also use a centered numerical flux.

The convergence order is computed using a self convergence method, obtaining the rate which the MR and MRLT methods solution converges to a solution as $\Delta t \rightarrow 0$.
For that, we perform simulations using subsequently smaller time steps, each one
having half of the time step used in the previous simulation. 
The convergence rate is obtained from the following ratio:
\begin{equation}
\left\|\frac{\bar{\mathbf{q}}_{\Delta t} - \bar{\mathbf{q}}_{\frac{\Delta t}{2}}}{\bar{\mathbf{q}}_{\frac{\Delta t}{2}} - \bar{\mathbf{q}}_{\frac{\Delta t}{4}}}\right\| = \frac{C\Delta t^p - C\left(\frac{\Delta t}{2}\right)^p + O(\Delta t^{p+1})}{C\left(\frac{\Delta t}{2}\right)^p - C\left(\frac{\Delta t}{4}\right)^p + O(\Delta t^{p+1})} = \frac{1- 2^{-p} + O(\Delta t)}{2^{-p}- 2^{-2p} + O(\Delta t)} = 2^p + O(\Delta t)
\end{equation}
where $\bar{\mathbf{q}}_{\Delta t}$ is the solution of a simulation using a time step $\Delta t$ and $p$ is the order of the method, which is approximated using the logarithm:
\begin{equation}
p \approx \log_2 \left\|\frac{\bar{\mathbf{q}}_{\Delta t} - \bar{\mathbf{q}}_{\frac{\Delta t}{2}}}{\bar{\mathbf{q}}_{\frac{\Delta t}{2}} - \bar{\mathbf{q}}_{\frac{\Delta t}{4}}}\right\|. 
\end{equation}

The convergence order obtained for the FV, MR and MRLT methods are given in Table \ref{convergenceTest}. 
In order to check that $p \rightarrow 0$ as $\Delta t \rightarrow 0$, we perform this test using two different values for the coarsest $\Delta t$. 
We observe that FV/RK2 and FV/RK3 yield second and third order time discretisations, respectively, as expected.
For MRLT/NERK2 and MRLT/NERK3 we obtained the expected second order, in particular for the MRLT/NERK3 the second order is justified by the approximation errors $\epsilon_3$ and $\epsilon_4$ which caused the loss of the third order, as shown in the Equations \ref{EvRk3finer} and \ref{EvRk3coarser}. Another reason for this loss in accuracy was the fact that due to the order barrier for NERK methods \cite{Zennaro91order}, a third order solution at instant $t^n + \frac{\Delta t_\ell}{2}$ could not be produced by a 3 stage method. 
Thus, when a leaf performs its second time evolution inside the LT cycle, it may use the second order solution from a coarser leaf, causing the loss from third to second order.
This issue also justifies the observed loss in accuracy for the MRLT/RK2 method. Here, the coarser leaf produces a first order solution at instant $t^n + \frac{\Delta t_\ell}{2}$ to be used in the second evolution of the finer leaf.  

\begin{table}[H]
\centering
\begin{tabular}{c|ccccccc}
\toprule
$\Delta t$ & \multicolumn{7}{c}{Method}\\
 {\small $(\times 10^{-4})$} & {\small FV/RK2} & {\small MR/RK2} & {\small MRLT/RK2} & {\small MRLT/NERK2} & {\small FV/RK3} & {\small MR/RK3} & {\small MRLT/NERK3}\\ \midrule
1.6 & 1.9991 & 1.9984 & 1.0794 & 2.0154 & 3.0000 & 3.0030 & 1.7895\\
0.8 & 2.0003 & 2.0010 & 0.9433 & 2.0002 & 3.0077 & 3.0017 & 1.9152\\ 
\bottomrule
\end{tabular}
\caption{Convergence order for the FV, MR and MRLT methods.}
\label{convergenceTest}
\end{table}

\section{Numerical experiments}\label{results}

In this section we present some comparative results of the proposed MRLT/NERK method with the MR, MRLT methods given in \cite{Domingues20083758} and also the traditional FV method on a uniform grid. 
These methods are applied to solve the two-dimensional Burgers equation, one and three-dimensional reaction-diffusion equations and finally the two-dimensional compressible Euler equations. 
We use the AUSM+ scheme \cite{Liou:1996} to compute the numerical flux in Burgers and Euler equations. To compute the advective term for the reaction--diffusion equation, we use the McCormack scheme \cite{VanLeer:1977}.
The errors are computed in the discrete $L_1$ norm on the fine grid as:
\begin{equation}
e_{L_1}^{\text{method}} = \frac{1}{2^{L\mathit{d}}}\sum_{i=0}^{2^{L\mathit{d}}-1}\left\| \bar{\mathbf{q}}_i^{\text{ref}} - \bar{\mathbf{q}}_i^{\text{method}} \right\|
\end{equation}
where $\bar{\mathbf{q}}^{\text{ref}}$ is the FV/RK3 reference solution of the corresponding problem and $\bar{\mathbf{q}}^{\text{method}}$ is the solution obtained with the analyzed method. 
\par To compare the performance of two adaptive methodologies in terms of CPU time reduction versus accuracy loss, a cost value $\mu^{\text{method}}$ is defined for each adaptive method as:
\begin{equation}
\mu^{\text{method}} = \frac{e_{L_1}^{\text{method}} \cdot t_{\mathbf{CPU}}^{\text{method}}}{t^{\text{FV}}_{\mathbf{CPU}}},
\end{equation}
where $t_{\mathbf{CPU}}^{\text{method}}$ is the CPU time obtained for the analyzed method and $t_{\mathbf{CPU}}^{\text{FV}}$ is the CPU time of the $FV$ method with the same number of scales $L$ and Runge--Kutta of the same order.
\par The ratio between the cost of different adaptive methods yields the parameter $\lambda$, used to measure the advantage of one method compared with the other. In this work, the parameter $\lambda$ is used to compare the proposed MRLT/NERK methods with the MR and MRLT methods, defined as:
\begin{equation}
\lambda^{\text{method}}_{\text{MRLT/NERK}} = \frac{\mu^{\text{method}}}{\mu^{\text{MRLT/NERK}}} \qquad .
\end{equation}
If the parameter $\lambda$ is larger than $1$, the MRLT/NERK approach is considered to be advantageous over the other method. In case of $\lambda < 1$, the MRLT/NERK approach is considered to be disadvantageous over the other method, and in case of $\lambda = 1$, the methods are considered to be equivalent.  

\par All the simulations are performed using a fixed threshold value  $\epsilon$ in order to simplify the experiments. 
In \cite{Dominguesetal:2012CFL} computations are presented using a MR methodology with $\epsilon$ values which depend of the refinement level.  
\subsection{Two-dimensional Burgers equation}
The Burgers equation is a non-linear PDE which represents a simple model for turbulence and is used in astrophysical applications. The inviscid model, in the two-dimensional case, is given by the following equation:
\begin{equation}
\frac{\partial Q}{\partial t} + 
\frac{1}{2} \left(\frac{\partial (Q^2)}{\partial x} + \frac{\partial (Q^2)}{\partial y}\right) = 0 \qquad (x,y)\in \Omega = [0,1]\times[0,1] .
\end{equation}
\par The initial condition used in this work is $Q_0(x,y) = \sin(2\pi x) \sin(2\pi y)$, with Dirichlet boundary conditions given by
$
Q(x,0) = Q(x,1) = Q(0,y) = Q(1,y) = 0. 
$
\par All simulations are performed with a Courant number $\sigma = 0.5$ and a threshold $\epsilon = 0.01$ until the time instant $t_f = 0.9$. The reference solution for this case is obtained using refinement level $L=12$. 
\par
\begin{figure}[htb]
\begin{center}
\begin{tabular}{ccc} 
\multicolumn{3}{c}{\textbf{Solution}}\\
 \multicolumn{1}{l}{a)\qquad \; \; Reference}  &  \multicolumn{1}{l}{b) \qquad \; \;  MRLT/RK2}  & \multicolumn{1}{l}{c) \qquad \; \; MRLT/NERK2 }  \\

\includegraphics[width=0.3\linewidth]{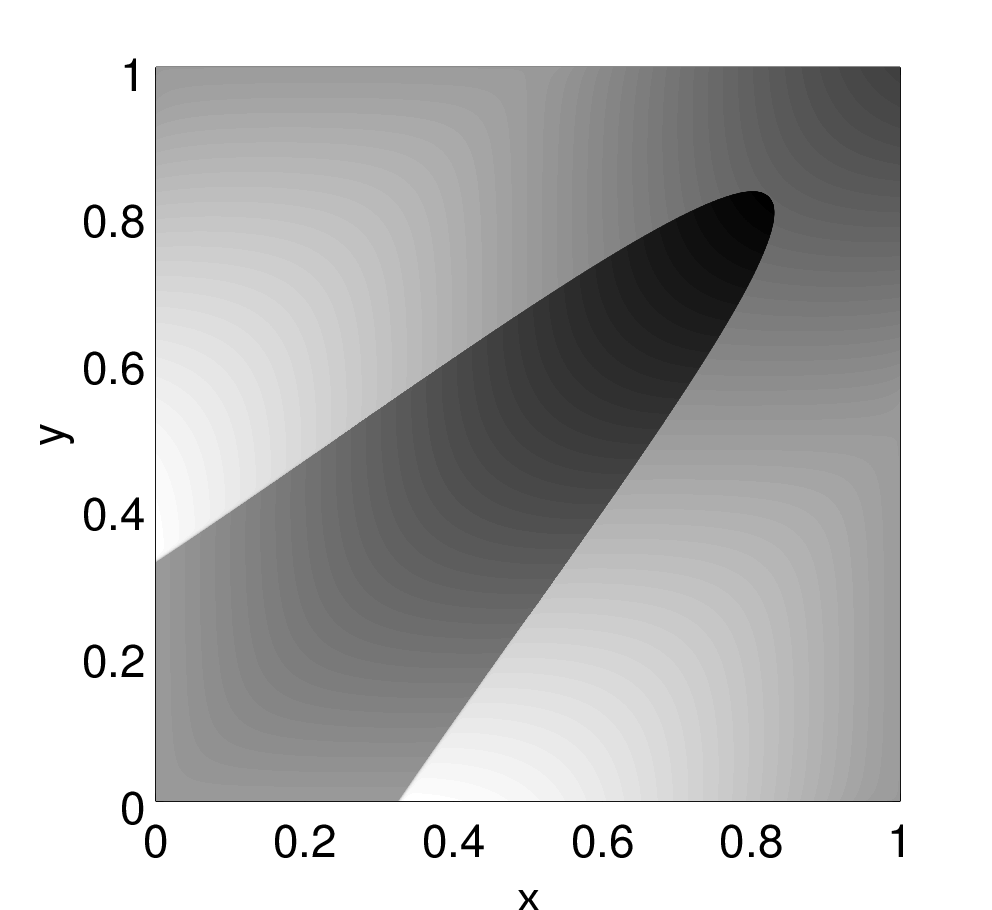}
&
\includegraphics[width=0.3\linewidth]{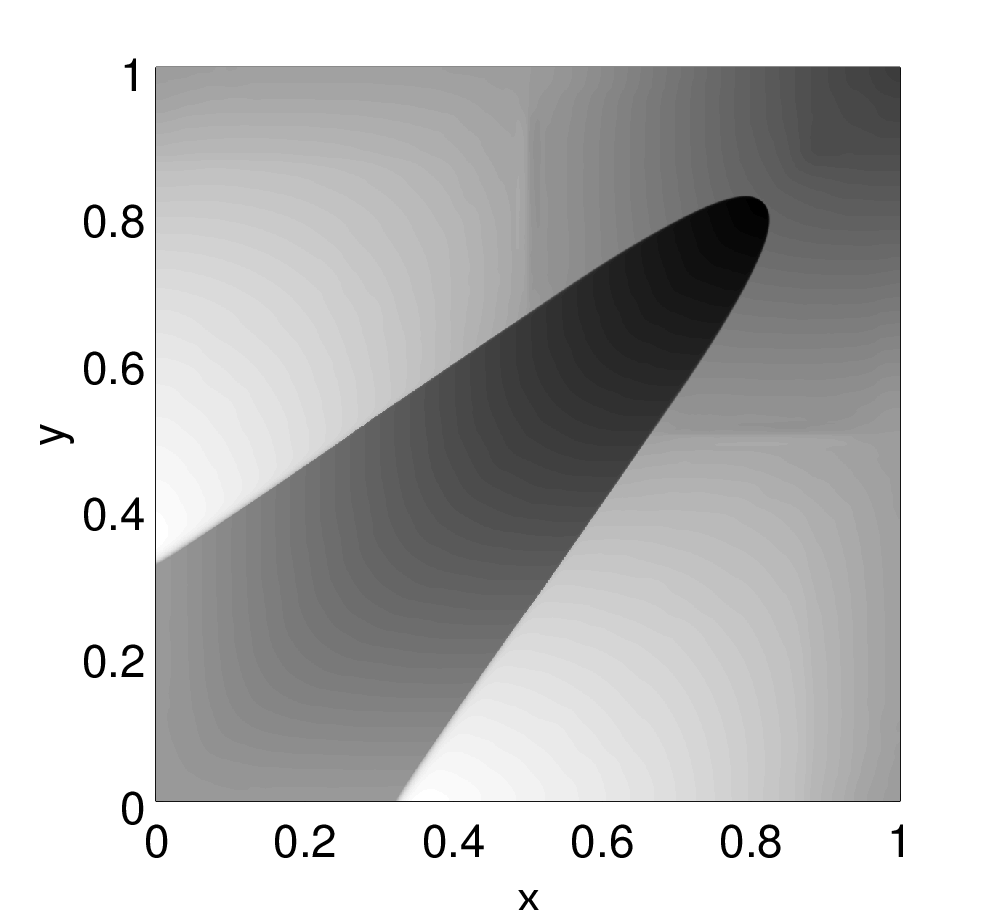}
&
\includegraphics[width=0.3\linewidth]{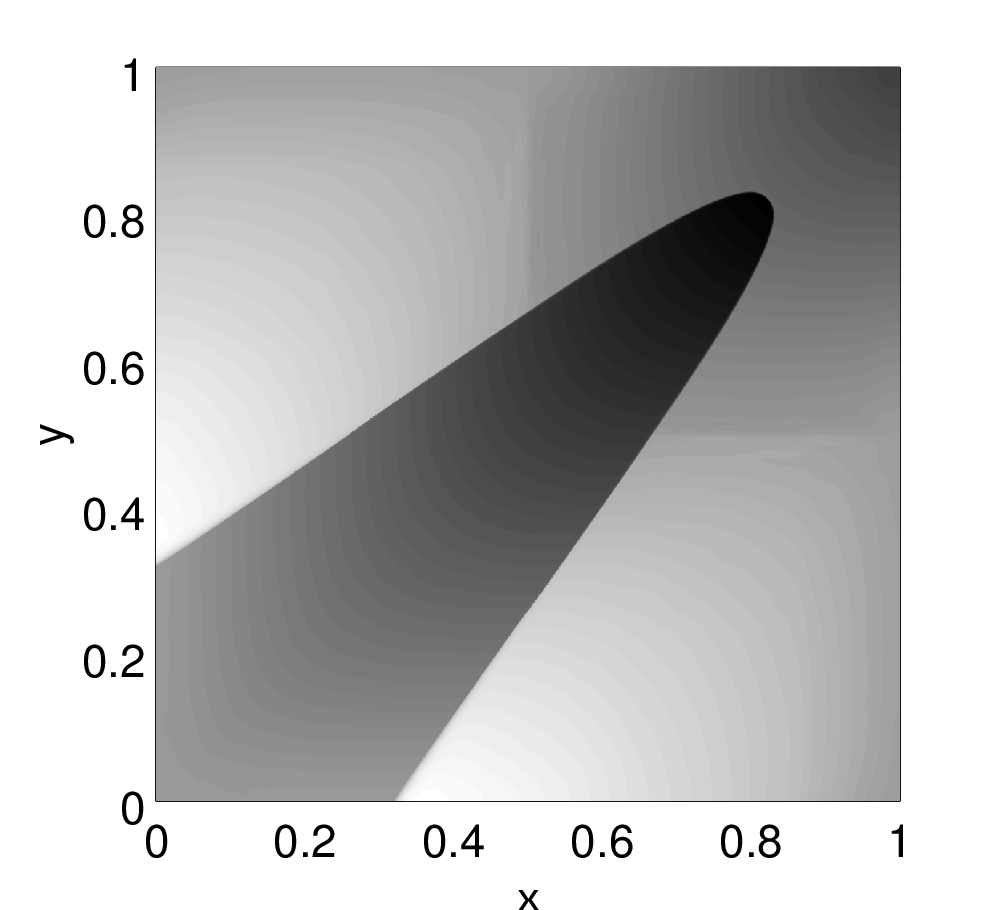}
\\
\multicolumn{3}{c}{\includegraphics[width=0.2\linewidth]{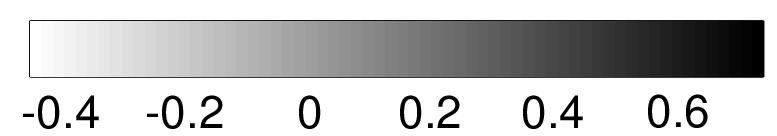}}\\
\\
\multicolumn{3}{c}{\textbf{Error}}\\
 &  \multicolumn{1}{l}{d) \qquad \; \; MRLT/RK2 }  & \multicolumn{1}{l}{e) \qquad \; \; MRLT/NERK2 }  \\
&
\includegraphics[width=0.3\linewidth]{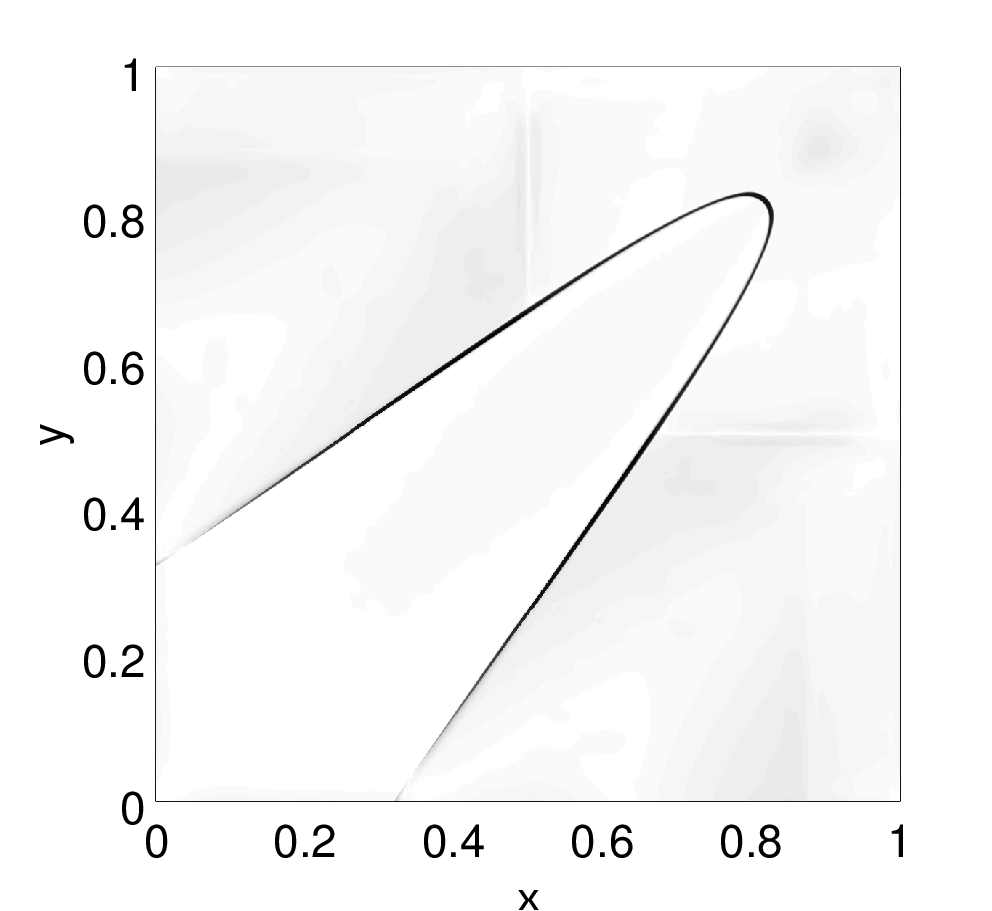}
&
\includegraphics[width=0.3\linewidth]{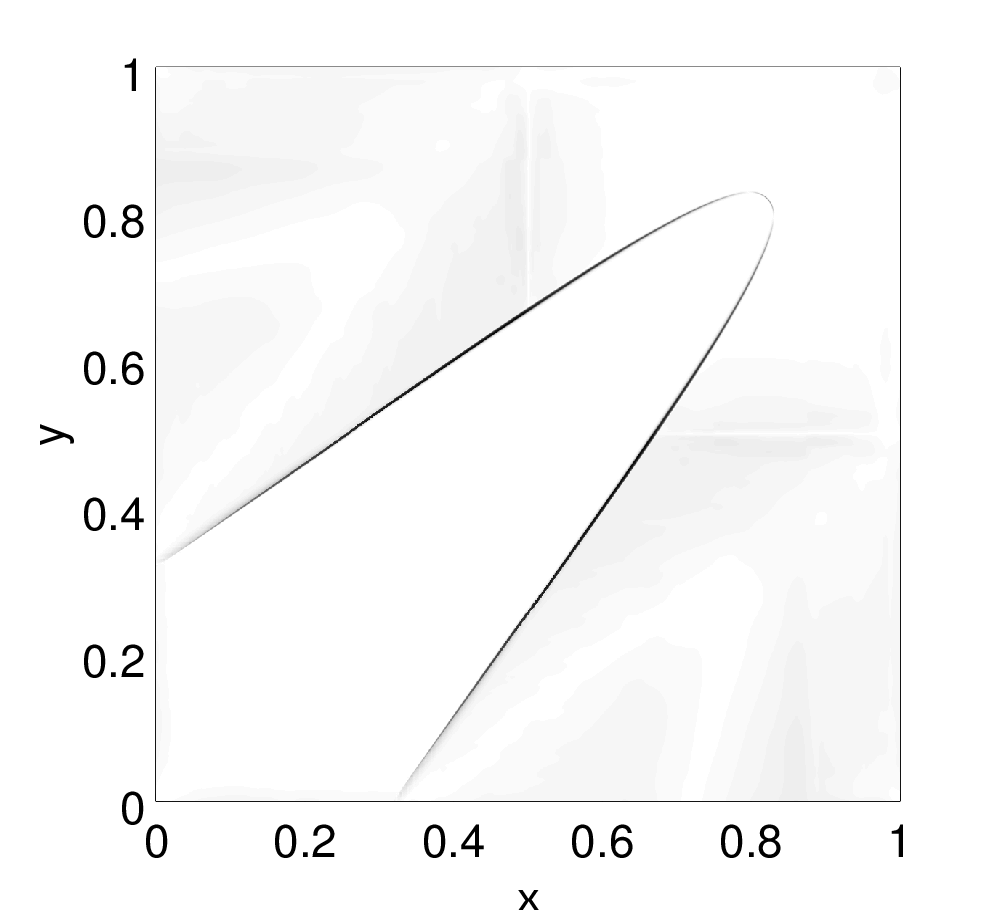}
\\
\multicolumn{3}{c}{\includegraphics[width=0.2\linewidth]{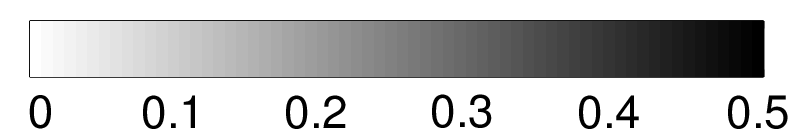}}\\
\\
\multicolumn{3}{c}{\textbf{Adaptive grid}}\\

 &  \multicolumn{1}{l}{f) \qquad \; \; MRLT/RK2}  & \multicolumn{1}{l}{g) \qquad \; \; MRLT/NERK2 }  \\

&
\includegraphics[width=0.3\linewidth]{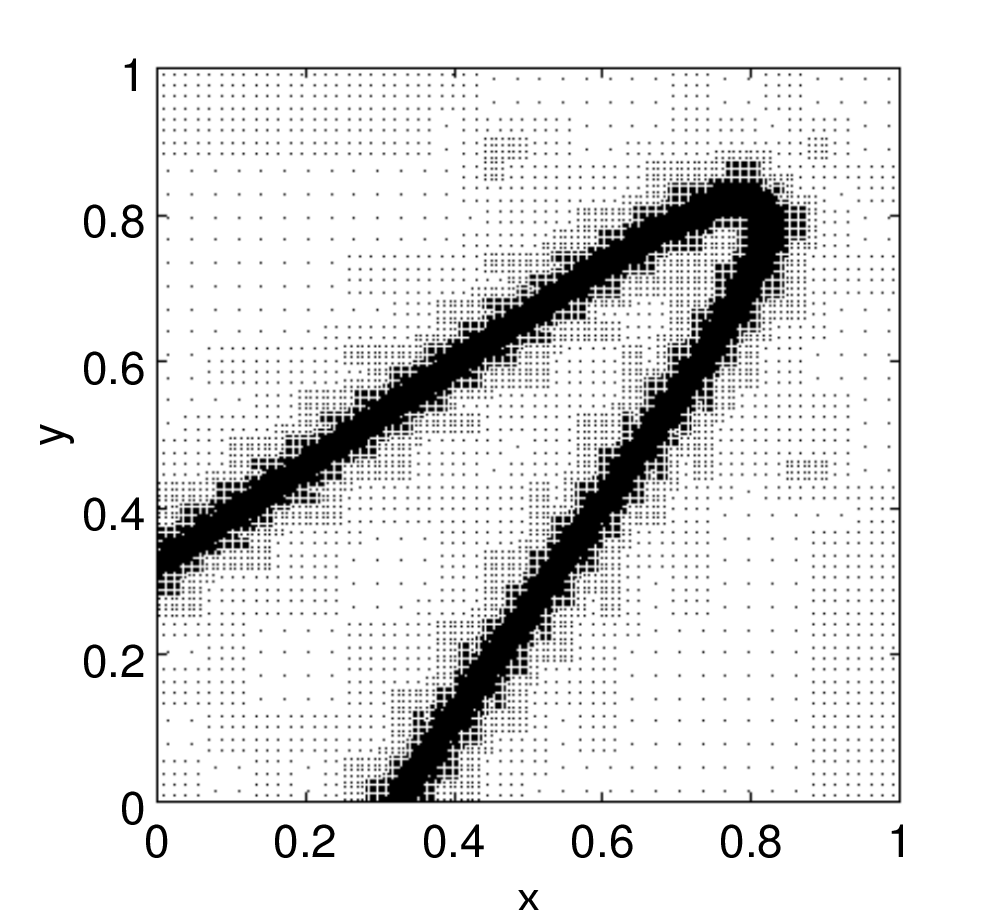}
&
\includegraphics[width=0.3\linewidth]{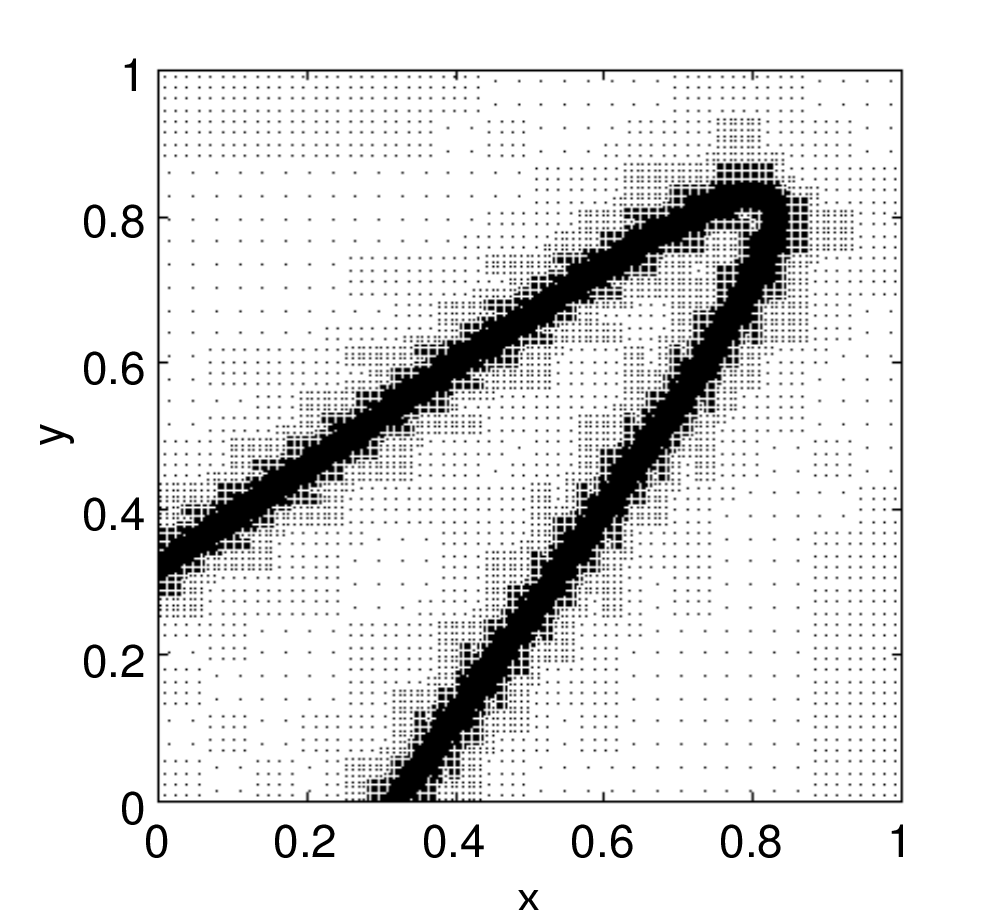}
\\
\end{tabular}
\end{center}
\caption{Reference solution for the two-dimensional Burgers equation a), the solutions obtained by the MRLT/RK2 b) and MRLT/NERK2 c) methods with its respective errors and the corresponding adaptive grids. For all cases we use $L=10$ and the time instant is $t_f= 0.9$.} 
\label{Burgers2Dfigure}

\end{figure}

 Figure \ref{Burgers2Dfigure} shows the reference solution and the solution obtained with the MRLT/RK2 and MRLT/NERK2. The respective difference, in modulus, with respect to the reference solution and adaptive grids at the end of the computation are also shown. The method MRLT/RK2 case exhibits larger errors close to the shocks, especially in the peak of the structure and in its background. The other methods present solutions closer to the one obtained with the MRLT/NERK2 method.

\par The results are compared with the solution in an uniform grid at the same level using a $L_1$ norm, showing perturbation errors. These errors, CPU time and memory compression are summarized in Table \ref{tableBurgers2D}. The CPU time and memory of the adaptive methods are given in percentage of  the number of leaves used in the FV method with the same number of scales and the same Runge--Kutta scheme. For the two-dimensional Burgers equation, the proposed MRLT/NERK methods present a slight gain in precision and a significant gain in CPU time in relation to the other adaptive methodologies. However, the MRLT schemes with NERK time integration require more memory, which decreases when increasing the resolution, \textit{i.e.} for increasing $L$.  
\par The parameters $\lambda$ obtained for the MRLT/NERK methods compared to the MR and MRLT methods are presented in Table \ref{gainBurgers2D}. 
For most of the experiments, the proposed methods yield values of $\lambda$ between $2$ and $3$. This shows that the NERK-based methods are significantly more efficient than the RK-based MR and MRLT methods.

\begin{table}[hbt]
\caption{Two-dimensional Burgers equation: $L_{1}$ errors, CPU time and memory compression.}
\begin{center}
\begin{tabular}{@{}clc*{2}{r}@{}} 
\toprule
Finest scale & Method & \multicolumn{1}{c}{Error ($\times 10^{-2}$)} & CPU Time &Memory \\ 
level & &  $Q$  & \multicolumn{2}{c}{($\%$ FV)}\\ 
\midrule
 & MR/RK2 & $1.5032$ & $68.4$ & $19.5$\\
 & MRLT/RK2 & $1.7701$ & $67.8$ & $19.5$\\
$L = 8$ & MRLT/NERK2 &$1.2958$ & $40.3$ & $24.2$\\
\cmidrule{2-5}
 & MR/RK3 & $1.5045$ & $184.9$ & $19.4$ \\
 & MRLT/NERK3 & $1.2932$ & $94.7$ & $24.2$ \\ 
 \midrule
 & MR/RK2 & $1.3615$ & $34.7$ & $8.9$\\
 & MRLT/RK2 & $1.6645$& $32.5$ & $8.9$\\
$L = 9$ & MRLT/NERK2 & $1.1457$ & $16.7$ & $10.3$\\
\cmidrule{2-5}
 & MR/RK3 & $1.3614$ & $95.4$ & $8.9$ \\
 & MRLT/NERK3 & $1.1423$ &$42.7$ &$10.3$  \\ \hline
 & MR/RK2 & $1.2890$ & $14.7$ & $4.1$\\
 & MRLT/RK2 & $1.6016$ & $13.5$ & $4.1$\\
$L = 10$ & MRLT/NERK2 & $1.0740$& $6.7$ & $4.4$ \\
\cmidrule{2-5}
 & MR/RK3 & $1.2893$ & $15.3$ & $4.1$ \\
 & MRLT/NERK3 & $1.0699$ & $6.4$ &$4.4$  \\
\bottomrule
\end{tabular}
\end{center}
\begin{footnotesize}
Note: All adaptive computations use $\epsilon = 10^{-2}$; the final time is $t_f = 0.9$. 
The computations have been carried out on an Intel Core\texttrademark $i7$ CPU $2.67$GHz. FV/RK2 CPU time: $2.0$ min ($L=8$); $15.7$ min ($L=9$); $2.3$ h ($L=10$). FV/RK3 CPU Time: $1.2$ min ($L=8$); $8.0$ min ($L=9$); $3.3$ h ($L=10$).
\end{footnotesize}
\label{tableBurgers2D}
\end{table}

\begin{table}[hbt]
\caption{Two-dimensional Burgers equation:\newline Computational gain $\lambda$ of the MRLT/NERK methods with respect to the MR and MRLT methods.}
\begin{center}
\begin{tabular}{@{}cl*{3}{c}@{}}  
\toprule
Finest scale &  & MR/RK2 & MRLT/RK2 & MR/RK3 \\ 
level & & & \\
\midrule
\multirow{2}{*}{$L = 8$}
 & MRLT/NERK2 &$1.96$ & $2.29$ & -\\
 & MRLT/NERK3 &  - & - & $2.27$ \\ 
 \midrule
\multirow{2}{*}{$L = 9$}
 & MRLT/NERK2 &$2.46$ & $2.82$ & -\\
 & MRLT/NERK3 &  - & - & $2.66$ \\ 
 \midrule
\multirow{2}{*}{$L = 10$}
 & MRLT/NERK2 &$2.63$ & $3.00$ & -\\
 & MRLT/NERK3 &  - & - & $2.88$ \\
\bottomrule
\end{tabular}
\end{center}
\label{gainBurgers2D}
\end{table}

\clearpage

\subsection{Reaction-diffusion equations} 
\label{cap:flameball}
We consider reaction-diffusion equations in one and three space dimensions and study the formation of a flame front ignited with a spark in an environment with flammable premixed gas.
 This kind of problem is a prototype of non-linear parabolic equations with a non-linearity in the source term, see e.g. \cite{RSTB03,Domingues20083758}.

\subsubsection{One-dimensional case}
In the one-dimensional case, considering equal mass and heat diffusion, this problem can be modeled by the following equation:
\begin{equation}
\frac{\partial T}{\partial t} + v_f \frac{\partial T}{\partial x} = \frac{\partial ^2 T}{\partial x^2} + \omega(T) \quad \textnormal{for} \; x \in (-15, 15)
\end{equation}
where the function $T(x,t)$ is the dimensionless temperature normalized between $0$ (unburned gas) and $1$ (burned gas), $v_f = \int\omega \;dx$ is the flame velocity and $\omega(T)$ is the chemical reaction rate, given by: 

\begin{equation}
\omega(T) = \frac{Ze^2}{2} (1-T) \exp\left(\frac{Ze(1-T)}{\tau (1-T)-1}\right)
\end{equation}
where $Ze$ is a dimensionless activation energy, know as Zeldovich number, and $\tau$ is the burnt-unburnt temperature ratio. In this one-dimensional case, the unburned gas concentration $Y$ is defined as $Y = 1-T$.

\par In the numerical experiments, the following initial condition is used:
\begin{align}
T(x,0) = 
\begin{cases}
1, \qquad & \text{if }x \le 1\\
\exp({1-x}), \qquad &  \text{otherwise}
\end{cases}
\end{align}
with boundary conditions given by: 
\begin{equation}
\frac{\partial T}{\partial x} (-15,t) = 0, \qquad T(15,t) = 0
\end{equation}

\par
\begin{figure}[htb]
\begin{center}
\begin{tabular}{lll} 
\multicolumn{3}{c}{\textbf{Solution}}\\
 \multicolumn{1}{l}{$a)$  \qquad \; \; Reference}  &  \multicolumn{1}{l}{$b)$  \qquad \; \; MRLT/RK2}  & \multicolumn{1}{l}{$c)$  \qquad \; \; MRLT/NERK2}  \\
\includegraphics[width=0.3\linewidth]{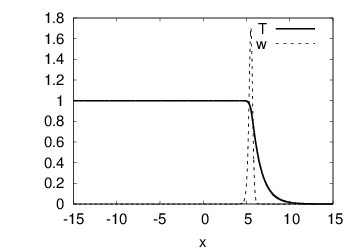}
&
\includegraphics[width=0.3\linewidth]{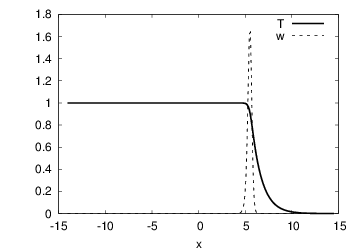}
&
\includegraphics[width=0.3\linewidth]{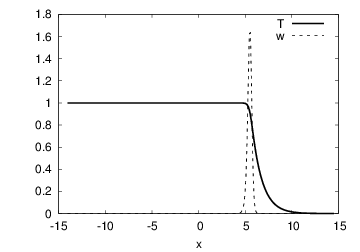}
\\
\\
\multicolumn{3}{c}{\textbf{Error}}\\
 \multicolumn{1}{l}{ }  &  \multicolumn{1}{l}{$d)$ \qquad \; \; MRLT/RK2}  & \multicolumn{1}{l}{$e)$ \qquad \; \; MRLT/NERK2 }  \\
&
\includegraphics[width=0.3\linewidth]{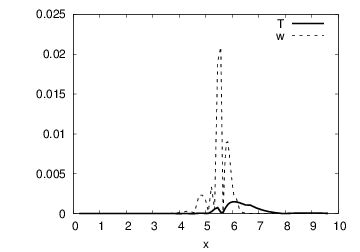}
&
\includegraphics[width=0.3\linewidth]{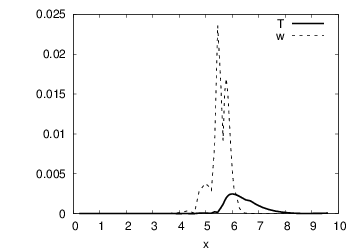}
\\
\\
\multicolumn{3}{c}{\textbf{Adaptive grid}}\\
 \multicolumn{1}{l}{ }  &  \multicolumn{1}{l}{$f)$ \qquad \; \; MRLT/RK2}  & \multicolumn{1}{l}{$g)$ \qquad \; \; MRLT/NERK2 }  \\
&
\includegraphics[width=0.3\linewidth]{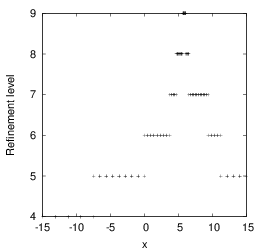}
&
\includegraphics[width=0.3\linewidth]{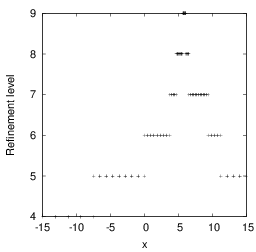}
\\
\end{tabular}
\end{center}
\caption{Reference solution of the one-dimensional reaction-diffusion equations a) and the solutions obtained by the MRLT/RK2 (b) and MRLT/NERK2  $(c)$ methods  with its respective errors $(d,e)$ and final adaptive grids $(f,g)$. For all cases we use $L=13$.}
\label{Flame1Dfigure}
\end{figure}

 The subsequent simulations are performed with the parameters $Ze = 10$, $\tau = 0.8$ and a Courant number $\sigma = 0.5$ until the final time instant $t_{f} = 5.0$. The adaptive simulations are performed with a threshold $\epsilon = 0.01$. The reference solution for this case is obtained using the same Courant number and a refinement level $L=13$. 
Figure \ref{Flame1Dfigure} shows the reference solution and the solution obtained by the MRLT/RK2 and MRLT/NERK2 methods with its respective errors and the corresponding final grids. In this one-dimensional case, the adaptive grid is represented by the position of each cell ($x$-axis) and its refinement level ($y$-axis). 
\par The adaptive methods present a larger error in the flame front region, especially for the variable $\omega$. Among the adaptive methods,  MRLT/RK2 has the smallest errors, while the other methods present very similar errors. However, the MRLT/NERK methods still requires the lowest CPU time. These times and $L_1$ errors are assembled in Table \ref{tableFlame1D}.   

\par The computational gain of the MRLT/NERK methods compared to the MR and MRLT methods are given in Table \ref{gainFlame1D}. In this case, the MRLT/NERK methods yields the most expressive results compared to the MR methods in terms of computational gain. The gain of the MRLT/NERK2 method compared with the MRLT/RK2 method is small for the cases $L=12$ and $13$. The corresponding results are presented in Table \ref{gainFlame1D}.

\begin{table}[hbt]
\caption{One-dimensional reaction-diffusion equations: $L_{1}$ errors, CPU time and memory compression.}
\begin{center}
\begin{tabular}{@{}cl*{3}{c}cc@{}} 
\toprule
Finest scale & Method & \multicolumn{3}{c}{Error ($\times 10^{-4}$)} & CPU Time & Memory \\ 
level & &   $T$ &  $Y$&  $\omega$  & \multicolumn{2}{c}{($\%$FV)}\\ 
\midrule
 & MR/RK2 & $5.045$ & $5.045$ & $24.564$ & $31.9$ & $3.1$\\
 & MRLT/RK2 & $5.666$ & $5.666$ & $30.192$ & $18.0$ & $3.1$\\
$L = 11$ & MRLT/NERK2 & $4.380$ & $4.380$ & $21.308$ & $7.2$ & $3.1$\\
\cmidrule{2-7}
 & MR/RK3 & $5.045$& $5.045$ & $24.566$  & $31.3$ & $3.1$ \\
 & MRLT/NERK3 & $4.543$ & $4.543$ & $22.173$ &$5.3$& $3.1$ \\ \midrule
 & MR/RK2 & $5.043$ & $5.043$ & $24.559$ & $15.9$ & $1.5$\\
 & MRLT/RK2 & $1.780$ & $1.780$  & $14.256$ & $8.0$ & $1.5$\\
$L = 12$ & MRLT/NERK2 & $4.728$ & $4.728$ & $23.107$ & $2.8$ & $1.5$\\
\cmidrule{2-7}
 & MR/RK3 & $5.043$ & $5.043$ & $24.560$  & $16.0$ & $1.5$ \\
 & MRLT/NERK3 & $4.755$  & $4.755$ & $23.261$ &$2.3$& $1.5$ \\ 
\midrule
 & MR/RK2 & $5.054$ & $5.054$ & $24.609$ & $7.9$ & $0.7$\\
 & MRLT/RK2 & $2.904$ & $2.904$ & $15.062$ & $3.9$ & $0.7$\\
$L = 13$ & MRLT/NERK2 & $4.900$ & $4.900$ & $23.903$ & $1.4$ & $0.7$\\
\cmidrule{2-7}
 & MR/RK3 & $5.054$ & $5.054$ & $24.609$ & $7.9$ & $0.7$ \\
 & MRLT/NERK3 & $4.904$ & $4.904$ & $23.923$ &$1.1$& $0.7$ \\ 
\bottomrule
\end{tabular}
\end{center}
\begin{footnotesize}
Note: All adaptive computations use $\epsilon = 10^{-2}$; final time: $t_f = 5.0$. Computed on an Intel Core\texttrademark $i7$ CPU $2.67$GHz. FV/RK2 CPU time: $51.4$ min ($L=11$); $6.9$ h ($L=12$); $54.6$ h ($L=13$). FV/RK3 CPU Time: $51.4$ min ($L=11$); $6.9$ h ($L=12$); $54.6$ h ($L=13$).
\end{footnotesize}
\label{tableFlame1D}
\end{table}

\begin{table}[hbt]
\caption{One-dimensional reaction-diffusion equations: Computational gain, for the variable $T$, of the proposed MRLT/NERK methods compared to the MR and MRLT methods.}
\begin{center}
\begin{tabular}{@{}cl*{3}{c}@{}}  
\toprule
Finest scale &  & MR/RK2 & MRLT/RK2 & MR/RK3 \\ level & & & & \\
\midrule
\multirow{2}{*}{$L = 11$}
 & MRLT/NERK2 & $5.10$ & $3.23$ & -\\
 & MRLT/NERK3 &  - & - & $6.55$ \\ \midrule
\multirow{2}{*}{$L = 12$}
 & MRLT/NERK2 &$6.05$ & $1.07$ & -\\
 & MRLT/NERK3 &  - & - & $7.37$ \\ \midrule
\multirow{2}{*}{$L = 13$}
 & MRLT/NERK2 &$5.82$ & $1.65$ & -\\
 & MRLT/NERK3 &  - & - & $7.40$ \\ 
 \bottomrule
\end{tabular}
\end{center}
\label{gainFlame1D}
\end{table}

\clearpage

\subsubsection{Three-dimensional case}
In the three-dimensional case, the reaction-diffusion equations read:
\begin{subequations}
\begin{equation}
\frac{\partial T}{\partial t} = \nabla ^2 T + \omega - s
\end{equation}
\begin{equation}
\frac{\partial Y}{\partial t} = \frac{1}{Le}\nabla ^2 Y - \omega
\end{equation}
\end{subequations}
where $Le$ denotes the Lewis number, which defines the ratio of mass and heat diffusion and with the chemical reaction rate $\omega$:
\begin{equation}
\omega(T,Y) = \frac{Ze ^2}{2Le} Y \exp\left( \frac{Ze(T-1)}{1 + \tau(T-1)}\right)
\end{equation}
According to the Stefan--Boltzmann law, the heat loss due to radiation $s$ is modeled by: 
\begin{equation}
s(T) = \kappa \left[ (T + \tau ^{-1} -1)^4 - (\tau^{-1} -1)^4 \right]
\end{equation}
where $\kappa$ is a dimensionless radiation coefficient. In this work, we use $\kappa = 0.1$. 

\par The initial condition, described by spherical coordinates, is:
\begin{align}
T(r,0) &= 
\begin{cases}
1, \qquad & \text{if }r \le r_0 \\
\exp\left(1 - \frac{r}{r_0}\right), \qquad &  \text{otherwise}
\end{cases}
\\
Y(r,0) &=
\begin{cases}
0, \qquad & \text{if }r \le r_0 \\
1 - e\left(Le(1 - \frac{r}{r_0})\right), \qquad &  \text{otherwise}
\end{cases}
\end{align}

where $r_0 = 1$ is the initial radius of the ellipsoidal flame ball and $r= \sqrt{\frac{X^2}{a^2} + \frac{Y^2}{b^2} + \frac{Z^2}{c^2}}$ with:

\begin{align}
X &= x\cos (\theta) - y\sin (\theta)\\
Y &= \left[ x\sin(\theta) + y\cos(\theta)\right] \cos(\phi) - z\sin(\phi)\\
Z &= \left[ x\sin(\theta) + y\cos(\theta)\right] \sin(\phi) + z\cos(\phi)
\end{align}

\begin{figure}[htb]
\begin{center}
\begin{tabular}{ccc} 
\multicolumn{3}{c}{\textbf{Solution}}\\
 \multicolumn{1}{l}{a) \; \; Reference}  &  \multicolumn{1}{l}{b) \; \; MRLT/NERK2}  & \multicolumn{1}{l}{c)   \; \; MRLT/NERK3}  \\
\includegraphics[width=0.2\textwidth]{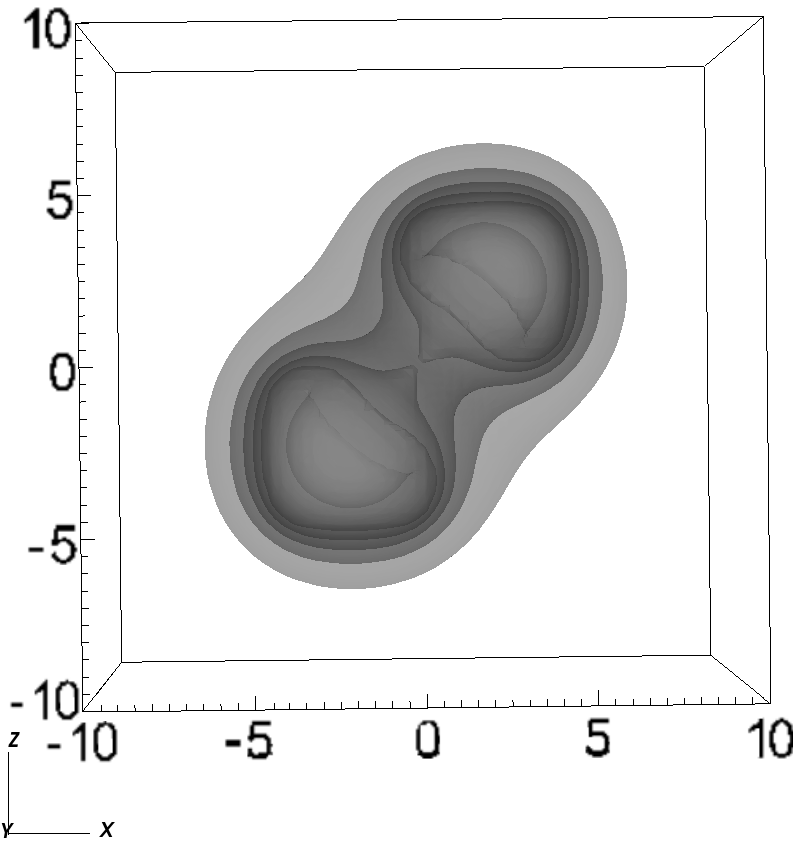}
&
\includegraphics[width=0.2\textwidth]{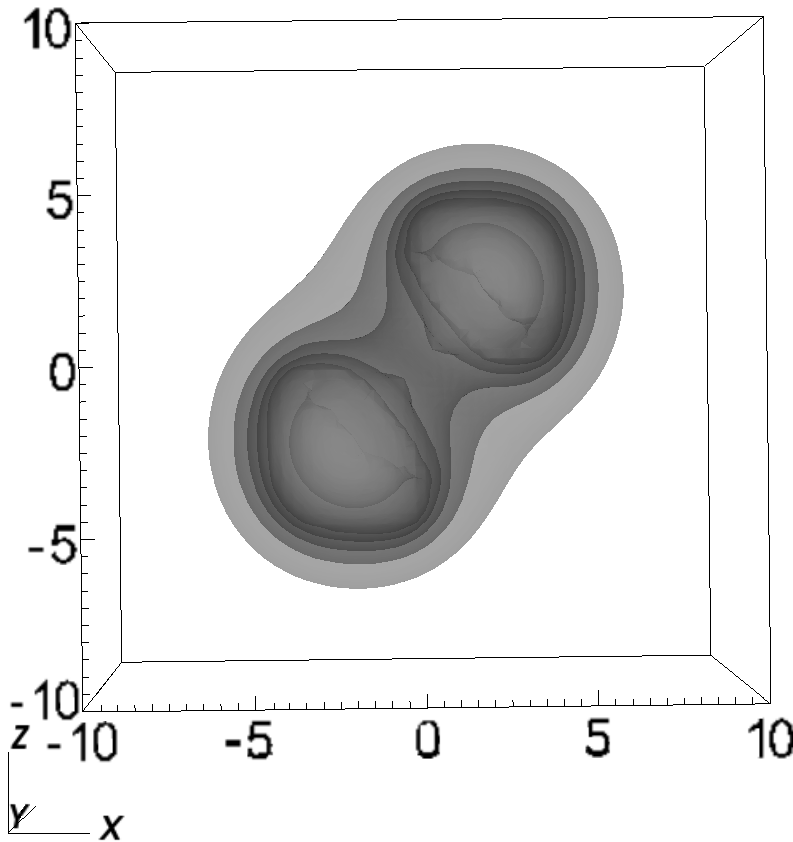}
&
\includegraphics[width=0.2\textwidth]{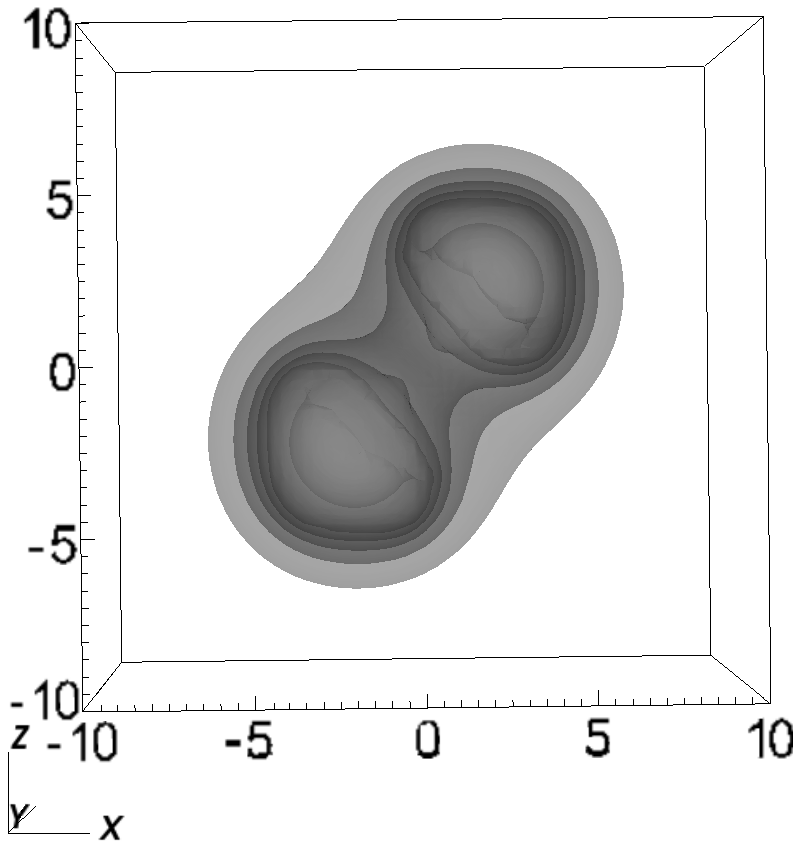}
\\
\multicolumn{3}{c}{\textbf{Solution in the plane $xz$}}\\
 \multicolumn{1}{l}{d)  \; \; Reference}  &  \multicolumn{1}{l}{e)   \; \; MRLT/NERK2}  & \multicolumn{1}{l}{f)   \; \; MRLT/NERK3}  \\
\includegraphics[width=0.2\textwidth]{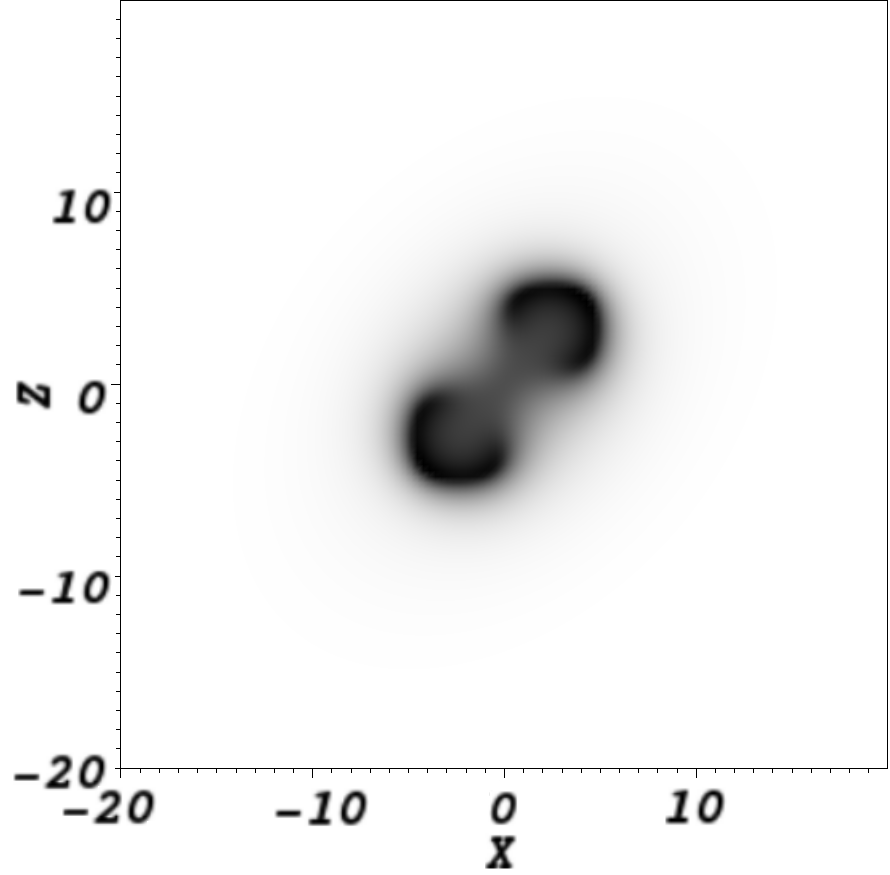}
&
\includegraphics[width=0.2\textwidth]{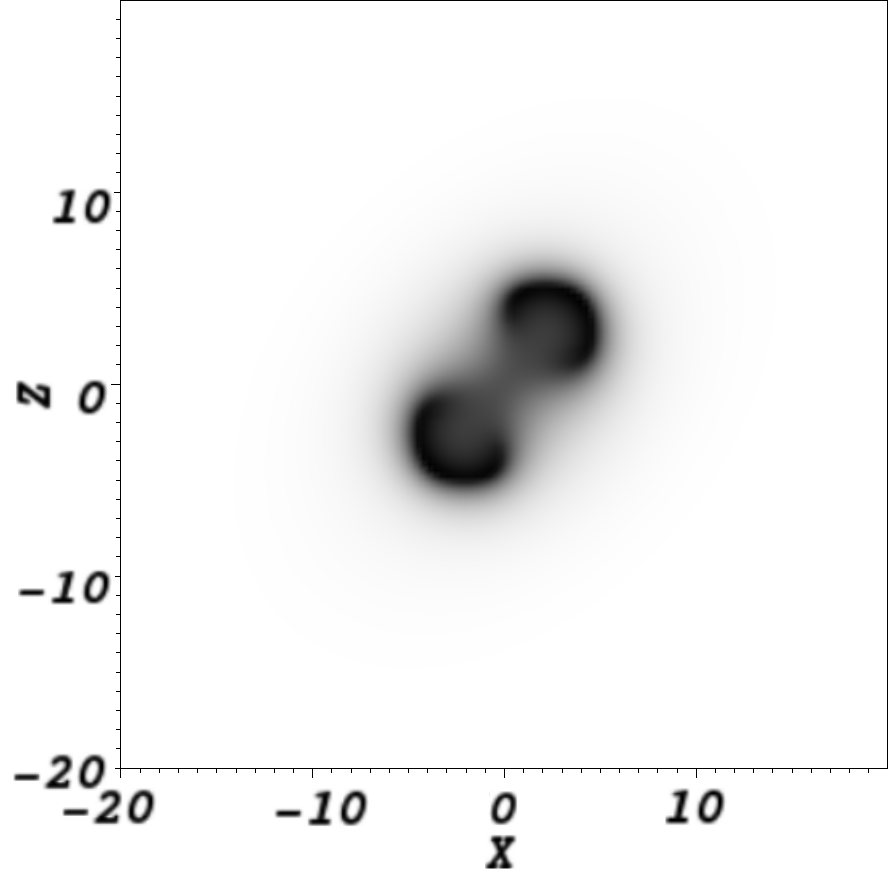}
&
\includegraphics[width=0.2\textwidth]{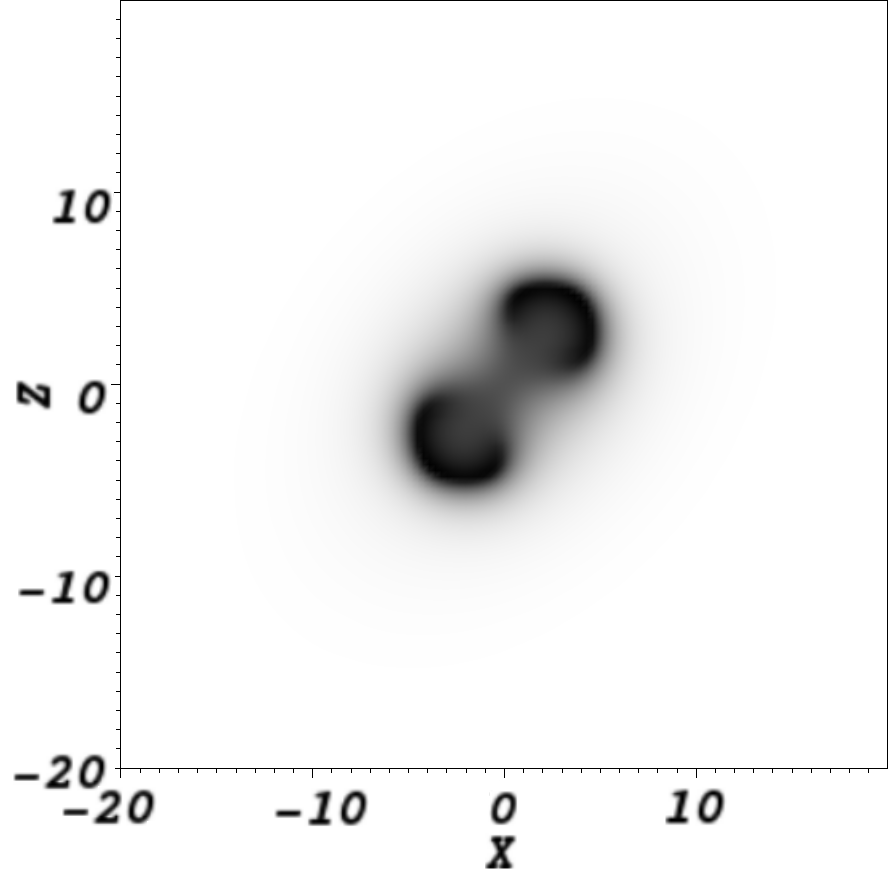}
\\
\multicolumn{3}{c}{\includegraphics[width=0.2\linewidth]{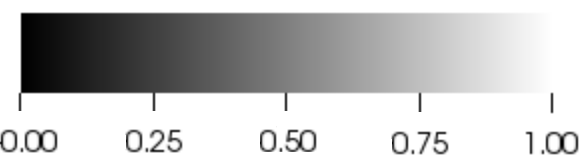}}\\
\multicolumn{3}{c}{\textbf{Error in the plane $xz$}}\\
 &  \multicolumn{1}{l}{g)  \; \; MRLT/NERK2}  & \multicolumn{1}{l}{h)  \; \; MRLT/NERK3}  \\
&
\includegraphics[width=0.2\textwidth]{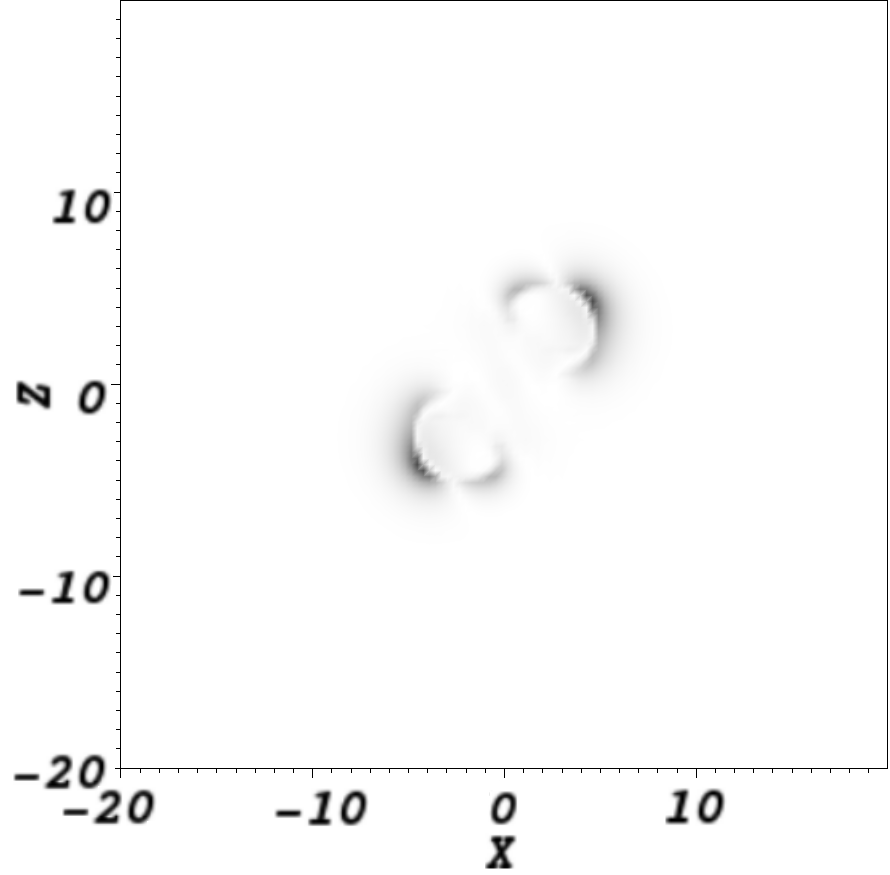}
&
\includegraphics[width=0.2\textwidth]{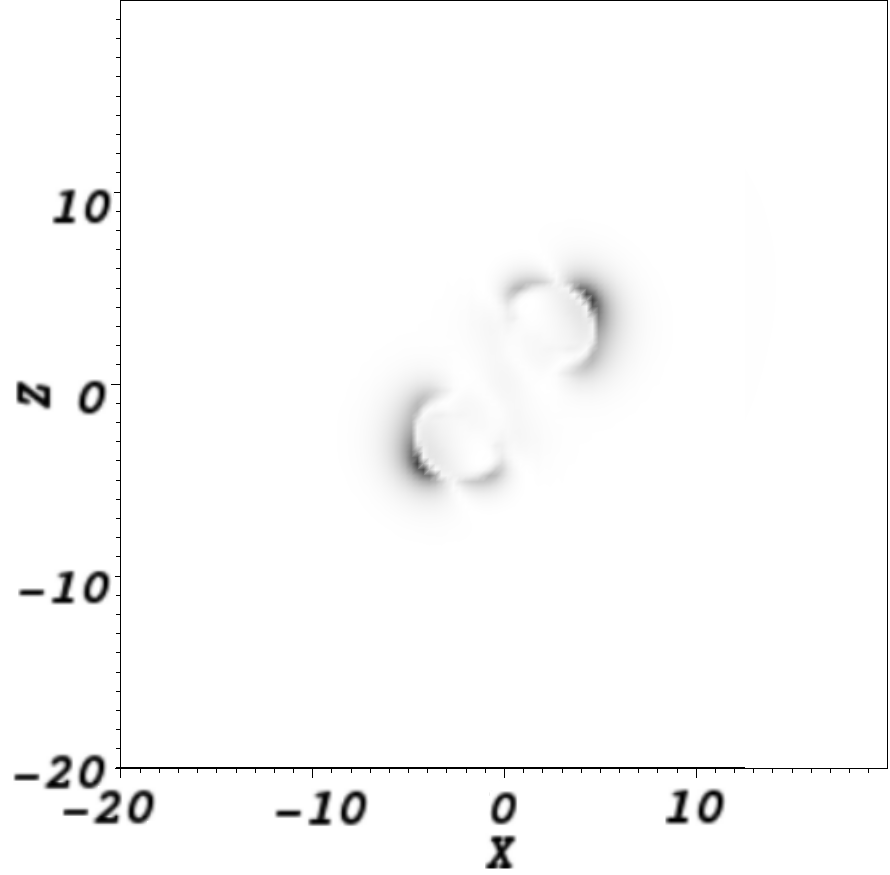}
\\
\multicolumn{3}{c}{\includegraphics[width=0.2\linewidth]{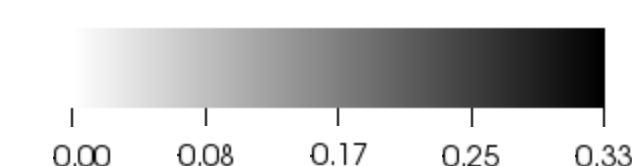}}\\

\multicolumn{3}{c}{\textbf{Grid projection in the plane $xz$}}\\
 &  \multicolumn{1}{l}{i)  \; \; MRLT/NERK2}  & \multicolumn{1}{l}{j)  \; \; MRLT/NERK3}  \\
&
\includegraphics[width=0.25\textwidth]{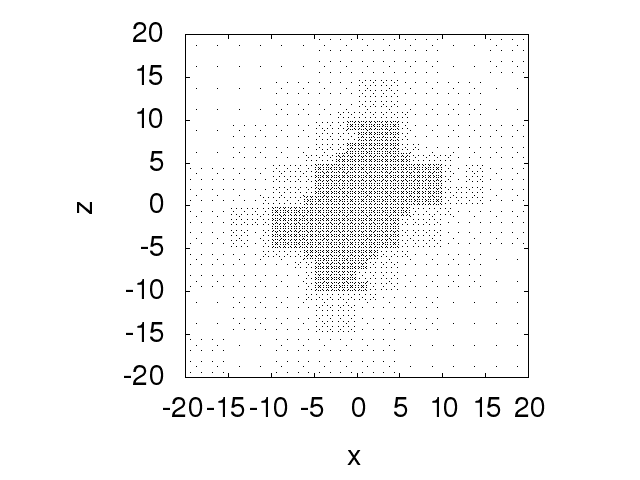}
&
\includegraphics[width=0.25\textwidth]{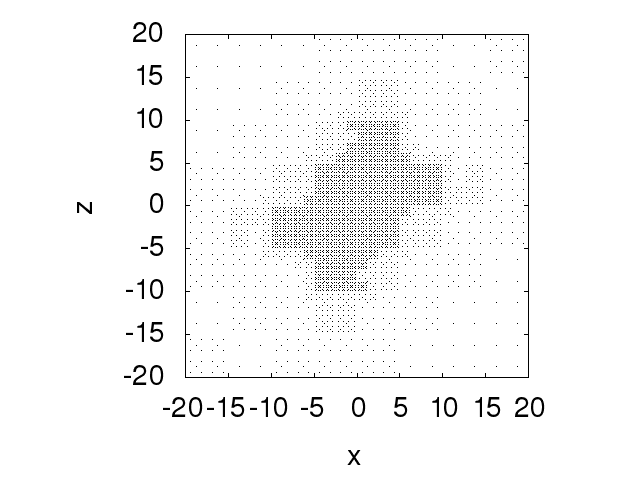}

\end{tabular}
\end{center}
\caption{Isosurface and $xz$ plane reference solution for the variable $T$ in three-dimensions at $t_f = 5.0$. The solutions are obtained with the MRLT/NERK2 and MRLT/NERK3 methods with its respective errors in the $xz$ plane and projections of the center of the adaptive grid cell onto the $xz$ plane. The isosurface plot were build using 5 linearly scaled values.}
\label{Flame3Dtemperature}
\end{figure}

The boundary conditions are of homogeneous Neumann type. 
In these simulations we use the parameters $Ze = 10$, $\tau = 0.64$, $Le = 0.3$, $\theta = \frac{\pi}{3}$, $\phi  = \frac{\pi}{4}$, $a = \frac{3}{2}$, $b=\frac{3}{2}$, $c = 3$, a threshold factor $\epsilon = 0.01$ and a Courant number $\sigma = 0.1$ until the final time instant $t_f = 5.0$. The reference solution for this case is obtained using the same Courant number and the refinement level $L=7$. Figure \ref{Flame3Dtemperature} shows the isosurface of $T$ for the reference and the MRLT/NERK2 and MRLT/NERK3 methods. It also shows the solutions, difference from the reference solution in modulus, and projections of every cell center at the plane $xz$. The adaptive methodologies present higher errors close to the flame ball fronts, which are more rounded in comparison to the FV solutions. The adaptive grids are similar for all adaptive methodologies, with a higher concentration of refined cells in the region of the front and inside the flame ball.   

\par The $L_1$ errors, CPU time and memory compression are summarized in Table \ref{tableFlame3D}. 
In this case, the proposed MRLT/NERK methods present a loss in precision and memory usage with a gain in CPU time in relation to the other adaptive methodologies. However, the parameters $\lambda$ obtained for the MRLT/NERK methods compared to the MR and MRLT methods, presented in Table \ref{gainFlame3D}, still yield favorable values, especially for the third order methods. Thus we find that the MRLT/NERK methods are slightly more efficient than the MR and MRLT methods for the three-dimensional flame ball case, besides a small loss in precision.

\begin{table}[hbt]
\caption{Three-dimensional reaction-diffusion equations: $L_1$ errors, CPU time and memory compression.}
\begin{center}
\begin{tabular}{@{}cl*{3}{c}cc@{}} 
\toprule
Finest scale & Method & \multicolumn{3}{c}{Error ($\times 10^{-4}$)} & CPU Time & Memory \\ 
level & &   $T$ &  $Y$&  $\omega$  & \multicolumn{2}{c}{($\%$FV)}\\ \midrule
 & MR/RK2 & $4.516$ & $7.338$ & $21.831$ & $15.6$ & $2.5$\\
 & MRLT/RK2 & $4.511$& $7.267$ & $21.798$ & $14.9$ & $2.5$\\
$L = 7$ & MRLT/NERK2 & $4.513$ & $7.636$ & $22.038$ & $10.9$ & $7.7$ \\
\cmidrule{2-7}
 & MR/RK3 & $4.516$& $7.338$& $21.831$  & $14.4$ & $2.5$ \\
 & MRLT/NERK3 & $4.526$ & $7.369$ & $21.856$ & $5.5$ & $2.5$ \\
\bottomrule
\end{tabular}
\end{center}
\begin{footnotesize}
Note: All adaptive computations use $\epsilon = 10^{-2}$; final time: $t_f = 5.0$. Computed on a  Quad-Core AMD Opteron\texttrademark CPU $2.4$GHz. FV/RK2 CPU time: $25.8$h ($L=7$). FV/RK3 CPU Time: $38.5$h ($L=7$).
\end{footnotesize}

\label{tableFlame3D}
\end{table}

\begin{table}[hbt]
\caption{Three-dimensional reaction-diffusion equations: \newline Computational gain, for the variable $T$, of the  MRLT/NERK methods compared to the MR and MRLT methods.}

\begin{center}
\begin{tabular}{ll*{3}{c}} 
\hline
Finest scale &  & MR/RK2 & MRLT/RK2 & MR/RK3 \\ 
level & & & & \\ \midrule
\multirow{2}{*}{$L = 7$}
 & MRLT/NERK2 & $1.43$ & $1.36$ & -\\
 & MRLT/NERK3 &  - & - & $2.61$ \\

\hline
\end{tabular}
\end{center}
\label{gainFlame3D}
\end{table}

\clearpage
\subsection{Compressible two-dimensional Euler equations}\label{cap:euler}

The Euler equations, which describe the dynamics of a non-ionized gas, are given by the following relations: 
\begin{subequations}
\begin{equation}
\frac{\partial \rho}{\partial t} + \nabla \cdot (\rho \vec{v}) =0
\end{equation}
\begin{equation}
\frac{\partial \rho \vec{v}}{\partial t} + \nabla \cdot (\rho \vec{v})\vec{v} + \nabla p =0
\end{equation}
\begin{equation}
\frac{\partial E}{\partial t} + \nabla \cdot \bigg(\vec{v}( E+p)\bigg) =0
\end{equation}
\label{EulerEquation}
\end{subequations}
where  $\vec{v}=(v_x,v_y)$ is the velocity,  $\rho$ is the fluid density, $p$ is the pressure, and  the energy per mass unity $E$ is given by $E = \rho e + \frac{\rho \|\vec{v}\| ^2}{2}$. This system is completed by the equation of state of an ideal gas $p = \frac{\rho T}{ \gamma Ma^2}$, where $\gamma$ is the specific heat ratio, $T$ is the temperature and $Ma$ is the Mach number.
 
\par The initial condition used in this numerical test is the classical Lax--Liu configuration \#6  \cite{LaxLiu:1998,Dominguesetal:2012CFL}. 
In this initial condition configuration, the domain is divided into four quadrants, ordered from $1^{\text{st}}$ to $4^{\text{th}}$, and defined by the subdomains 
$\left[0.5;1\right]\times\left[0.5;1\right]$, $\left[0;0.5\right]\times\left[0.5;1\right]$, $\left[0;0.5\right]\times\left[0;0.5\right]$ and  $\left[0.5;1\right]\times\left[0;0.5\right]$ respectively. 
The initial values for each quadrant are given in Table~\ref{InitialCondition}.
This problem is simulated until the final instant $t_f = 0.25$, using homogeneous Neumann boundary conditions.
The numerical parameters are the Courant number $\sigma = 0.5$, $Ma = 1$, and $\gamma = 1.4$ inside the domain $[0;1]\times[0;1]$. 

\begin{table}[htb]
\caption{Initial condition for two-dimensional Euler equations.}
\begin{center}
\begin{tabular}{@{}lrrrr@{}} 
\toprule
\multirow{2}{*}{Variables} & \multicolumn{4}{c}{Quadrant}\\
&  $1^{\text{st}}$ & $2^{\text{nd}}$ & $3^{\text{rd}}$ & $4^{\text{th}}$\\ \midrule
Density, $\rho$ & $1.00$ & $2.00$ & $1.00$ & $3.00$ \\ 
Pressure, $p$ & $1.00$ & $1.00$ & $1.00$ & $1.00$ \\ 
$x$-velocity, $v_{x}$ & $0.75$ & $0.75$ & $-0.75$ & $-0.75$ \\ 
$y$-velocity, $v_{y}$ & $-0.50$ & $0.50$ & $0.50$ & $-0.50$ \\ 
\bottomrule
\end{tabular}
\end{center}
\label{InitialCondition}
\end{table}

\par
\begin{figure}[htb]
\begin{center}
\begin{tabular}{ccc} 
\multicolumn{3}{c}{\textbf{Solution}}\\
\multicolumn{1}{l}{$a)$  \qquad \; \; Reference}  &  \multicolumn{1}{l}{$b)$  \qquad \; \; MRLT/RK2 }  & \multicolumn{1}{l}{$c)$  \qquad \; \; MRLT/NERK2 }  \\

\includegraphics[width=0.3\linewidth]{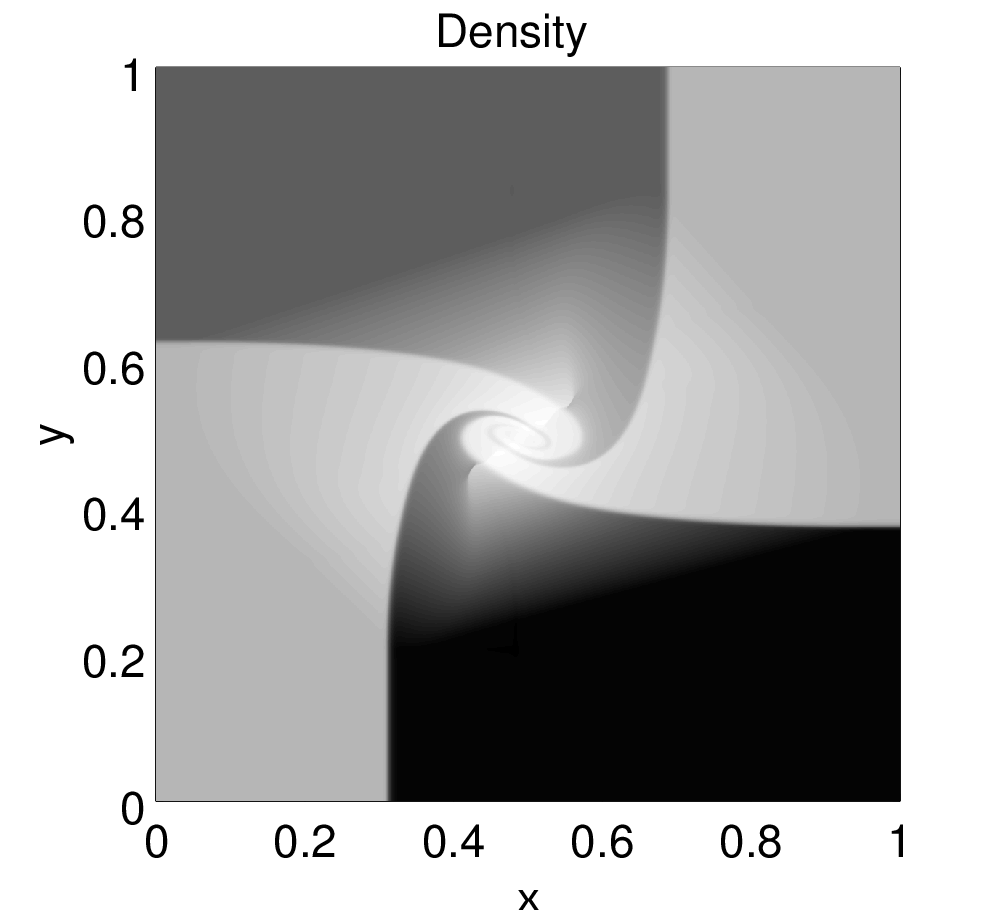}
&
\includegraphics[width=0.3\linewidth]{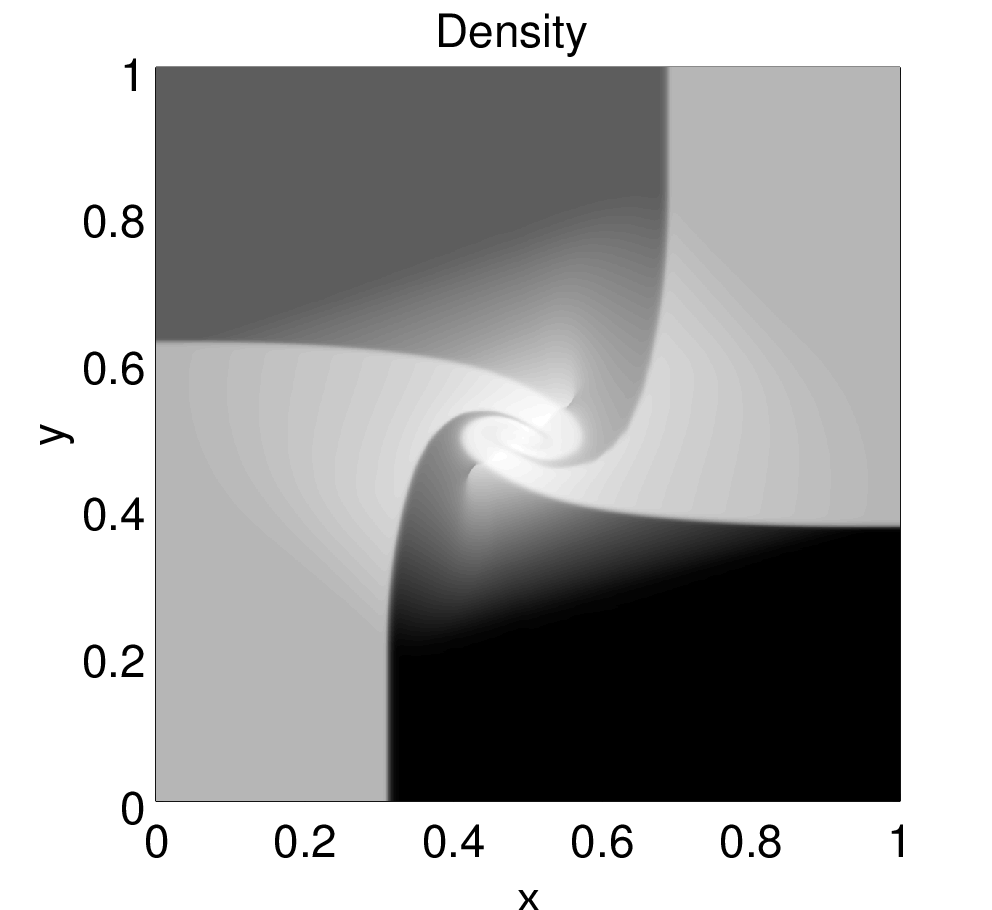}
&
\includegraphics[width=0.3\linewidth]{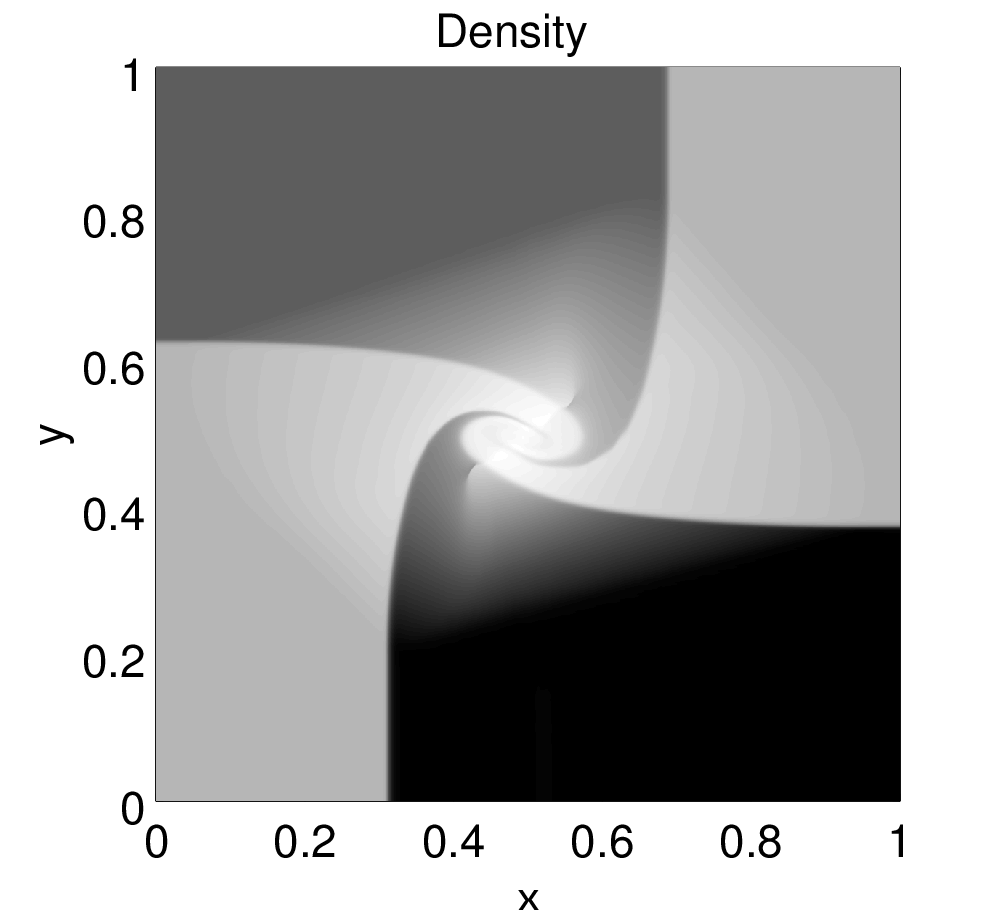}
\\
\multicolumn{3}{c}{\includegraphics[width=0.2\linewidth]{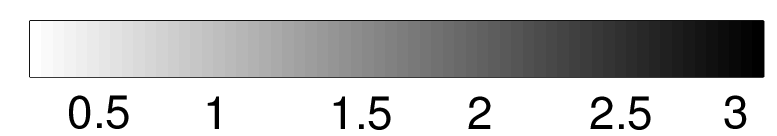}}\\
\\
\multicolumn{3}{c}{\textbf{Error}}\\
   &  \multicolumn{1}{l}{$d)$  \qquad \; \; MRLT/RK2 }  & \multicolumn{1}{l}{$e)$  \qquad \; \; MRLT/NERK2 }  \\
&
\includegraphics[width=0.3\linewidth]{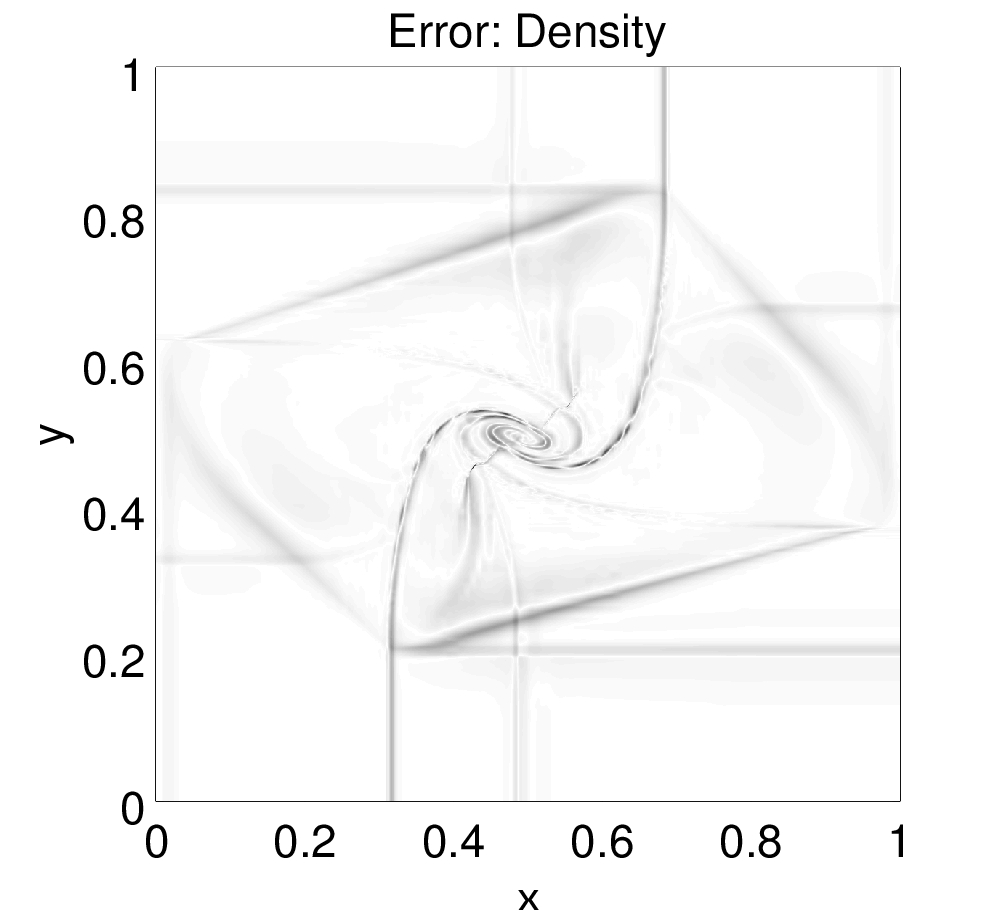}
&
\includegraphics[width=0.3\linewidth]{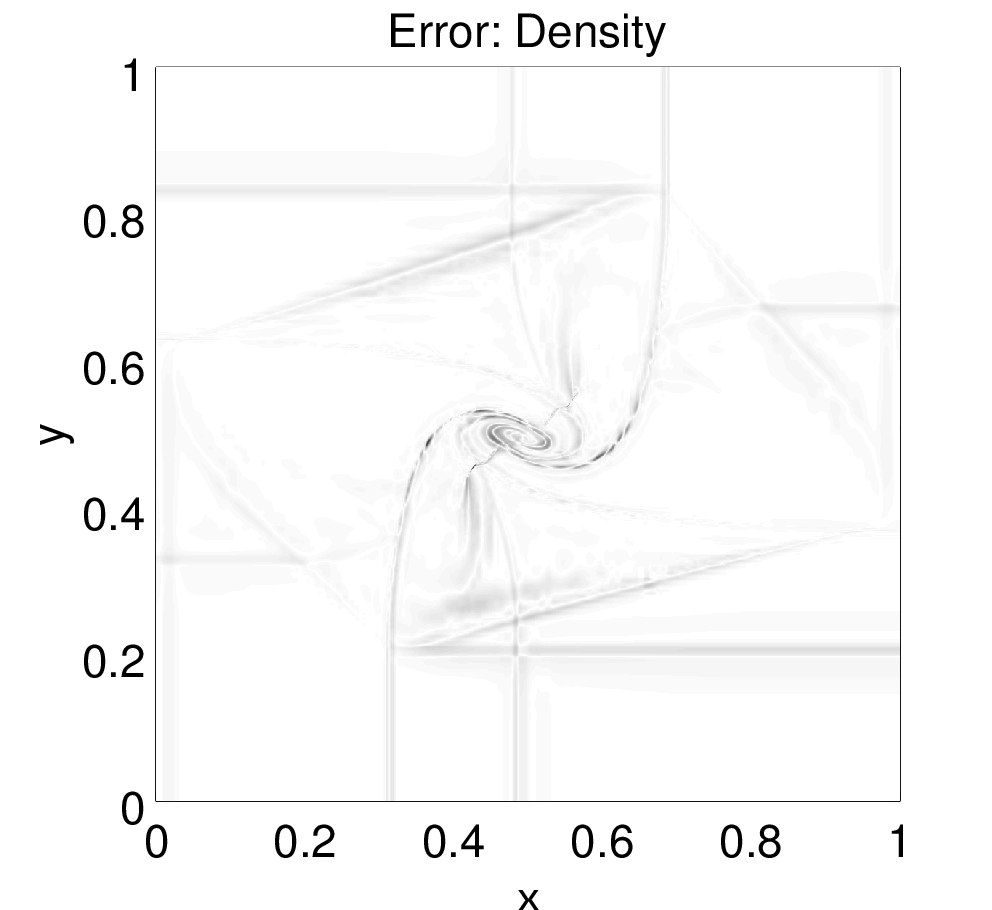}
\\
\multicolumn{3}{c}{\includegraphics[width=0.2\linewidth]{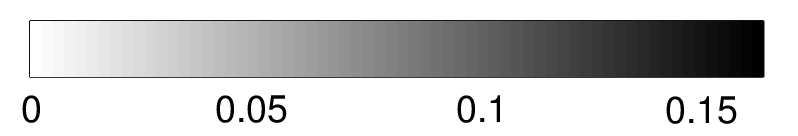}}\\
\\
\multicolumn{3}{c}{\textbf{Adaptive grid}}\\
   &  \multicolumn{1}{l}{$f)$  \qquad \; \; MRLT/RK2 }  & \multicolumn{1}{l}{$g)$  \qquad \; \; MRLT/NERK2 }  \\
&
\includegraphics[width=0.3\linewidth]{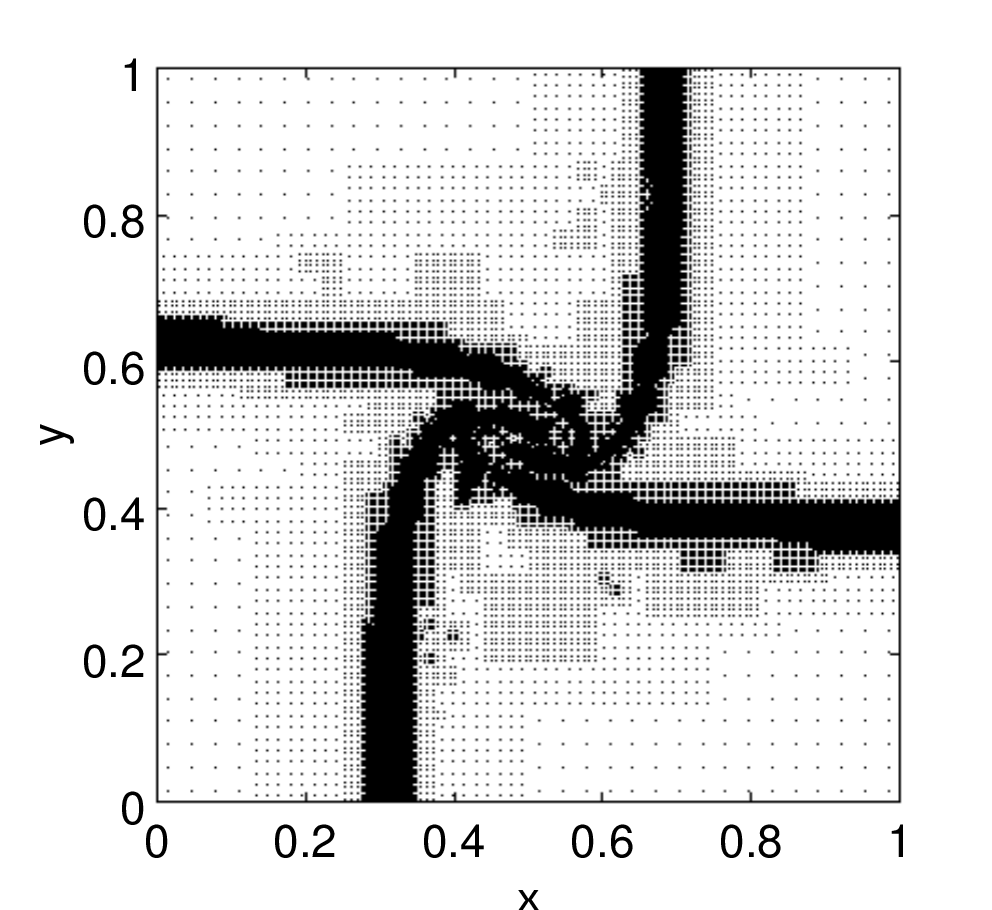}
&
\includegraphics[width=0.3\linewidth]{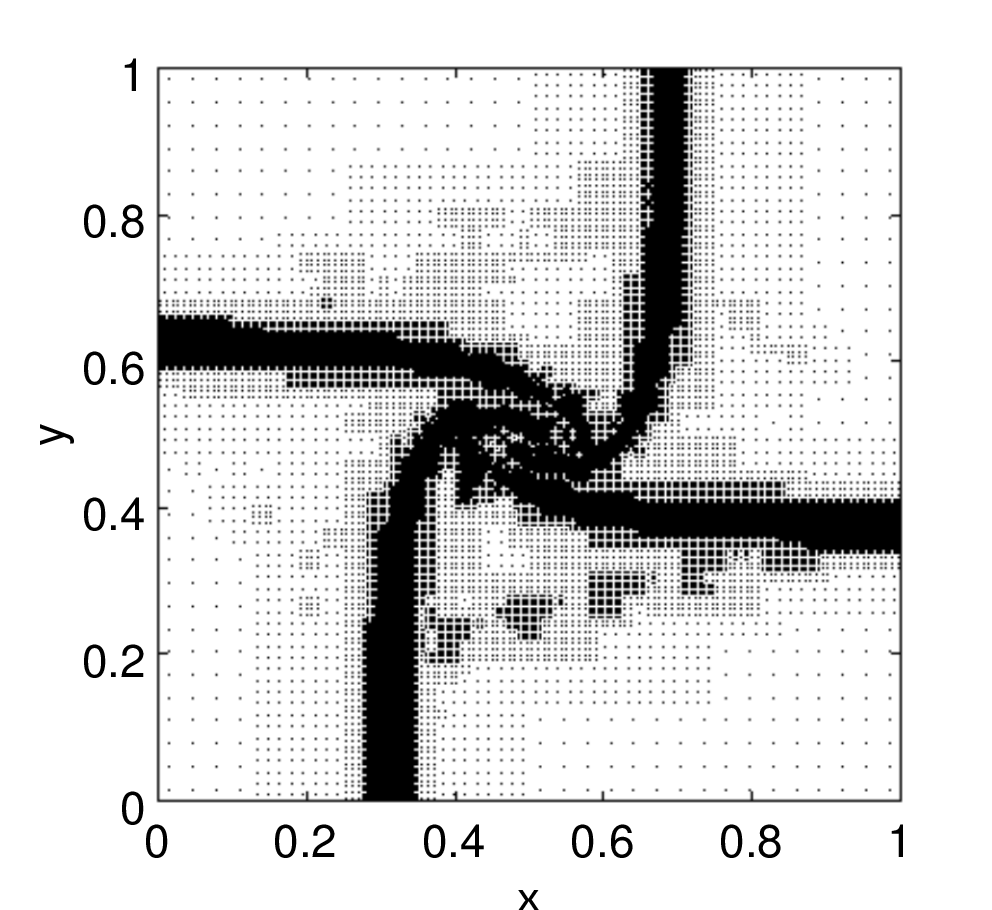}
\\
\end{tabular}
\end{center}
\caption{Reference solution for the two-dimensional Euler equations($a$), the solutions obtained by the MRLT/RK2 ($b$) and MRLT/NERK2 methods, using $L=10$ scales, with its respective errors($d, e$), and corresponding adaptive grids $(f, g)$ at final time $t_f = 0.25$.}
\label{Euler2DDensity}
\end{figure}

 Figure \ref{Euler2DDensity} shows the reference solution and the solution obtained with the MRLT/RK2, MRLT/NERK2 methods and its respective difference, in modulus, from the reference. Furthermore, the corresponding adaptive grid using $L=12$ with CFL $\sigma = 0.5$ is shown.

\par The $L_1$ errors, CPU time and memory compression are presented in Table \ref{tableEuler2D}. 
In this case, the proposed MRLT/NERK methods yield a slightly gain in precision with a significant gain in CPU time in relation to the other adaptive methodologies. 
The memory usage of all adaptive methodologies is quite similar. 
Due to the gain in both precision and CPU time, the parameters $\lambda$ for the MRLT/NERK methods compared with the MR and MRLT methods, presented in Table \ref{gainEuler2D} have expressive values, around $3$ for most of the experiments. 
Thus, using this metric, the proposed methods for the Euler two-dimensional case are significantly more efficient than the MR and MRLT methods.

\clearpage

\begin{table}[hbt]
\caption{Two-dimensional Euler equations: $L_1$ errors, CPU time and memory compression.}
\label{tableEuler2D}
\begin{center}
\begin{small}
\begin{tabular}{@{}cl*{6}{c}cc@{}} 
\toprule
Finest scale & Method & \multicolumn{6}{c}{Error ($\times 10^{-1}$)} & CPU Time & Memory \\ 
level & &  $\rho$ &  $p$&  $T$ & $E$ & $v_x$&$v_y$& \multicolumn{2}{c}{($\%$FV)}\\ \midrule
 & MR/RK2 & $2.2095$ & $0.6075$ & $0.9581$ & $2.0499$ & $1.7521$ & $0.7180$ & $40.7$ & $26.4$\\
 & MRLT/RK2 & $2.2085$ & $0.6049$ & $0.9572$ & $2.0463$ & $1.7510$ & $0.7167$ & $33.7$ & $26.0$\\
$L = 8$ & MRLT/NERK2 & $2.2096$ & $0.6077$ & $0.9582$ & $2.0503$ & $1.7521$ & $0.7183$  & $12.9$ & $26.6$\\
\cmidrule{2-10}
& MR/RK3 & $2.2093$ & $0.6075$ & $0.9580$ & $2.0498$ & $1.7520$ & $0.7180$ & $46.3$ & $26.4$ \\
 & MRLT/NERK3 & $2.2095$ & $0.6078$ & $0.9581$ & $2.0504$ & $1.7522$ & $0.7182$ &$14.6$& $26.6$ \\
\midrule 
 & MR/RK2 &  $1.1365$ & $0.3133$ &  $0.5217$ &  $1.0470$ &  $0.9012$ &  $0.3934$ & $24.4$ & $14.0$\\
 & MRLT/RK2 & $1.1369$ & $0.3133$ & $0.5219$ & $1.0483$ & $0.9012$ & $0.3931$  & $22.6$ & $13.9$\\
$L = 9$ & MRLT/NERK2 & $1.1355$ & $0.3128$ & $0.5216$ & $1.0458$ & $0.9012$ & $0.3931$ & $7.4$ & $14.2$\\
\cmidrule{2-10}
 & MR/RK3 &  $1.1364$ & $0.3133$ & $0.5217$ & $1.0469$ & $0.9012$ & $0.3934$ & $27.3$ & $14.0$ \\
 & MRLT/NERK3 &  $1.1355$ & $0.3130$ & $0.5216$ & $1.0461$ & $0.9013$ & $0.3931$ &$8.0$& $14.2$ \\ 
 \midrule
 & MR/RK2 & $0.5835$ & $0.1604$ & $0.2765$ & $0.5353$ & $0.4611$ & $0.2115$ & $14.7$ & $7.2$\\
 & MRLT/RK2 & $0.5864$ & $0.1643$ & $0.2782$ & $0.5441$ & $0.4622$ & $0.2131$ & $12.8$ & $7.2$ \\
$L = 10$ & MRLT/NERK2 & $0.5821$ & $0.1596$ & $0.2762$ & $0.5335$ & $0.4609$ & $0.2110$ & $4.1$ & $7.3$\\
\cmidrule{2-10}    
 & MR/RK3 & $0.5834$ & $0.1604$ & $0.2765$ & $0.5352$ & $0.4611$ & $0.2115$ & $11.2$ & $7.2$ \\
 & MRLT/NERK3 & $0.5821$ & $0.1598$ & $0.2762$ & $0.5338$ & $0.4610$ & $0.2109$ &$3.6$& $7.3$ \\
\bottomrule
\end{tabular}
\end{small}
\end{center}
\begin{footnotesize}
Note: All adaptive computations use $\epsilon = 10^{-2}$; final time: $t_f = 0.25$. Computed on an Intel Core\texttrademark $i7$ CPU $2.67$GHz. FV/RK2 CPU time: $10.1$ min ($L=8$); $73.3$ min ($L=9$); $8.8$ h ($L=10$). FV/RK3 CPU Time: $11.7$min ($L=8$); $91.7$ min ($L=9$); $13.8$ h ($L=10$).
\end{footnotesize}

\label{TableEuler2D}
\end{table}

\begin{table}[hbt]
\caption{Two-dimensional Euler equations: \newline Computational gain, for the variable $\rho$, of the proposed MRLT/NERK methods compared with the MR and MRLT methods.}
\begin{center}
\begin{tabular}{@{}cl*{3}{c}@{}} 
\toprule
Finest scale &  & MR/RK2 & MRLT/RK2 & MR/RK3 \\ 
level &  &  &  &  \\
\midrule
\multirow{2}{*}{$L = 8$}
 & MRLT/NERK2 & $3.15$ & $2.61$ & -\\
 & MRLT/NERK3 &  - & - & $3.17$ \\ \midrule
\multirow{2}{*}{$L = 9$}
 & MRLT/NERK2 &$3.30$ & $3.05$ & -\\
 & MRLT/NERK3 &  - & - & $3.41$ \\ \midrule
\multirow{2}{*}{$L = 10$}
 & MRLT/NERK2 &$3.59$ & $3.14$ & -\\
 & MRLT/NERK3 &  - & - & $3.11$ \\
\bottomrule
\end{tabular}
\end{center}
\label{gainEuler2D}
\end{table}

\section{Conclusions}\label{sec:conclusions}

In this work, we introduced a new local time-stepping for adaptive multiresolution methods using NERK time integration schemes \cite{zennaro1986natural}.
Interpolating values of the intermediate Runge--Kutta stages yield the required values at intermediate time steps, which are necessary for the time evolution.
Hence the current limitation of local time stepping to second order schemes can be overcome, and the required synchronised solution can be obtained.
The proposed new methodology has been implemented and validated for two and three stage NERK schemes.
In principle the extension of MRLT/NERK schemes to even higher order is possible, but the computational cost would increase due to the order barrier discussed in \cite{Zennaro91order}.
For NERK methods of order larger or equal to three, Owren and Zennaro \cite{OwrenZennaro:1992} have shown that the order of approximation is reduced concerning the underlying RK method.
This order reduction implies that increasing the order of NERK schemes beyond three would become less efficient and thus further research is indeed necessary to obtain well performing local time-stepping schemes with order larger than three.
With the presented numerical experiments we assessed the precision and efficiency of the proposed two and three stage MRLT/NERK approach for different classical nonlinear evolution equations, i.e., for Burgers, reaction-diffusion and the compressible Euler equations considering Cartesian geometries in one, two and three space dimensions.
In all adaptive computations, we observed a significant gain in CPU time in comparison with uniform grid computations where the efficiency is increasing with the grid resolution.
Nevertheless, the precision of the uniform grid computations is controlled in the MRLT/NERK schemes, and the order of convergence is maintained.
Regarding memory consumption, we observed for the MRLT/NERK schemes no necessary increase, compared to MR and classical MRLT methods.

The precision of the MRLT/NERK computations is very reasonable, and in all cases, we found errors about the same order of magnitude as for the MRLT computations using classical RK schemes.

In conclusion, we showed that the MRLT/NERK methods are advantageous compared to the MR and MRLT approaches in all studied cases, obtaining even significant performance gains in some examples.
For the two-dimensional Euler equations, for instance, the MRLT/NERK simulations only required one-third of the CPU time necessary for the MR and MRLT computations.
Most of these gains are due to the significant reduction in CPU time obtained in the MRLT/NERK methods.
 
\section*{Acknowledgements}

The authors are indebted to Prof. Claus-Dieter Munz who motivated the use of NERK methods for local time-stepping. We thank Dr. Olivier Roussel for developing the original Carmen Code and fruitful scientific discussions during the years.
ML thankfully acknowledges financial support from
CAPES and CNPq grant $140626/2014-0$, Brazil.
MD thankfully acknowledges financial support from
Ecole Centrale Marseille, CNPq grant $306038/2015\!-\!3$ and FAPESP grant $2015/25624\!-\!2$.
OM thankfully acknowledges financial support from CNPq grant $307083/2017-9 $ and FINEP grant $01.120527.00$, Brazil.
KS acknowledges financial support from the ANR Grant $15-$CE$40-0019$ (AIFIT), France.

\clearpage
\section{References}

\bibliographystyle{plain}
\bibliography{NERK}
 
\appendix
\input{Appendix_alg}

\end{document}

%% file: Appendix_alg.tex
\newpage
\section{ $\qquad \qquad$ Algorithms}\label{App:Alg}

\begin{algorithm}
\caption{Construction of the adaptive grid.}
\begin{algorithmic}
\label{alg:AdaptiveMesh}
\STATE Do Algorithm~\ref{alg:waveletThres} \COMMENT{Compute the leaves that will belong to the adaptive grid}\\
\STATE -- Refine each leaf with a finer neighbour leaf.
\STATE -- Set these new leaves as virtual leaves;
\end{algorithmic}
\end{algorithm}

\begin{algorithm}
\caption{Grid adaptation.}
\begin{algorithmic}
\label{alg:waveletThres}
\REQUIRE  Finest scale level $L$ of the solution.
\REQUIRE Solution at level $L$.
\REQUIRE  Threshold value $\epsilon$. \\
\COMMENT{Selection of the nodes in the adaptive grid}
\FOR{$\ell = L\!-\!1 \rightarrow 0$, $\ell=\ell\!-\!1$}
\STATE -- Project the solution from the grid $\Omega^{\ell+1}$ to $\Omega^{\ell}$;
\STATE -- Predict the solution of the grid $\Omega^{\ell+1}$ based on the project solution in $\Omega^{\ell}$ ;
\STATE -- Compare the  original solution of the grid $\Omega^{\ell+1}$ with the predicted one, and then obtain the wavelet coefficients $\bar{D}^{\ell+1}$, as the difference with these solutions;\\
\COMMENT{Elimination of the unnecessary nodes and the imposition of  a graded tree}
\FOR{every leaf $\in \Omega^{\ell+1}$}
\IF{$|\mathbf{d}^{\ell+1}| \leq \epsilon$ and adjacent leaves are inside $\Omega^{\ell+1}$ or $\Omega^{\ell}$}
\STATE -- Remove the leaf from the grid $\Omega^{\ell+1}$;
\ENDIF
\ENDFOR
\ENDFOR
\end{algorithmic}
\end{algorithm}

\begin{algorithm}
\caption{Single iteration of the LT scheme.}
\begin{algorithmic}
\label{alg-LTS}
\REQUIRE Coarsest scale $\ell_{\min}$ to be evolved in this time evolution; 
\REQUIRE Current iteration number $n$; 

\STATE Compute $\ell_{\min} = \min\limits_{\ell}\left[\mod\left( n,2^{L-\ell}\right) = 0 \right]$;
\FOR{$\ell = L \to \ell_{\min}$, $\ell\!=\!\ell\!-\!1$}
\FOR{Every internal node $\in \Omega^{\ell}$}
\STATE Obtain the solution $\bar{\mathbf{q}}^n_\ell$ by projecting the solution from its child nodes via simple averaging;
\ENDFOR
\ENDFOR

\IF{$\ell_{\min}$ is not the coarsest scale of the grid}
\FOR{Every internal node $\in\Omega^{\ell_{\min}-1}$}
\STATE Obtain the NERK solution with $\theta=\frac{1}{2}$ by projecting $\bar{\mathbf{q}}^n_{\ell_{\min}}$.
\STATE Extrapolate this solution to the instant $t^{n+2^{L-\ell_{\min}}}$ (Section \ref{LTissue2}).
\ENDFOR
\FOR{Every virtual leaf $\in \Omega^{\ell_{\min}}$}
\STATE Predict $\bar{\mathbf{q}}^n_\ell$ using the NERK (RK2) or $\bar{\mathbf{q}}^{**}_{\ell_{\min}-1}$ (RK3) solutions of the cells $\in \Omega^{\ell_{\min}-1}$ at instant $t^{n}$.
\ENDFOR
\ENDIF

\FOR{$\ell = \ell_{\min}+1 \to L$, $\ell\!=\!\ell\!+\!1$}
\FOR{Every virtual leaf $\in \Omega^{\ell}$}
\STATE Predict $\bar{\mathbf{q}}^n_\ell$ using the values of the cells $\in \Omega^{\ell-1}$ at time instant $t^{n}$. 
\ENDFOR
\ENDFOR

\STATE Remeshing process of cells with refinement level greater or equal than $\ell_{\min}$;   

\FOR{$\ell = L \to \ell_{\min}$, $\ell\!=\!\ell\!-\!1$}
\FOR{Every leaf $\in \Omega^{\ell}$}
\STATE Perform flux computations at instant $t^n$;
\STATE First RK step;
\STATE First order interpolation of the RK evolution at instant $t^n + \frac{1}{2}\Delta t_\ell$;
\IF{we are computing RK2 time evolution}
\STATE Compute $\bar{\mathbf{q}}^{*}_{\theta = \frac{1}{2}}$;
\ELSIF{we are computing RK3 time evolution}
\STATE Compute $\bar{\mathbf{q}}^{*}_{\theta = \frac{1}{4}}$; 
\STATE Compute $\bar{\mathbf{q}}^{*}_{\theta = \frac{3}{4}}$; 
\ENDIF
\ENDFOR
\ENDFOR

\STATE Tree refreshing before the second RK step (Algorithm \ref{alg:refRK2}):
\STATE Second RK step of the time evolution (Algorithm \ref{Steprk2});

\IF{Performing RK3 time evolution}
\STATE Tree refreshing before the third RK step (Section \ref{syncChallengeTree}):
\STATE Third RK step of the time evolution (Algorithm \ref{Steprk3});
\ENDIF

\end{algorithmic}
\end{algorithm}


\begin{algorithm}
\caption{Projection procedure inside the second RK step.}
\begin{algorithmic}
\label{alg:projRK2}
\REQUIRE  Scale $\ell$  to receive the projection.
\REQUIRE Number of dimensions $\mathit{d}$  of the problem.
\FOR{Every internal node $\in \Omega^{\ell}$}
\STATE $\bar{\mathbf{q}}^{*}_{\ell} = 0$; 
\COMMENT{result after RK1, set here to zero to store the projection}
\FOR{Every child cell $\in \Omega^{\ell+1}$} 
\IF{The cell is a leaf}
\STATE $\bar{\mathbf{q}}^{*}_{\ell}  \leftarrow \bar{\mathbf{q}}^{*}_{\ell} + \bar{\mathbf{q}}_{\ell+1} \left(t^n + 2\Delta t_{\ell+1}\right)$; 
\COMMENT{Add the already extrapolated value $\bar{\mathbf{q}}_{\ell+1}$  at $t^n + 2\Delta t_{\ell+1}$ }
\ELSIF{The cell is an internal node}
\STATE  $\bar{\mathbf{q}}^{*}_{\ell} \leftarrow \bar{\mathbf{q}}^{*}_{\ell} + 2\bar{\mathbf{q}}^{*}_{\ell+1} - \bar{\mathbf{q}}_{\ell+1}^n$;
\COMMENT{Add a linear extrapolation of the values $ \bar{\mathbf{q}}_{\ell+1}$ at $t^n + 2\Delta t_{\ell+1}$ }
\ENDIF
\ENDFOR
\STATE $\bar{\mathbf{q}}^{*}_{\ell}  \leftarrow \dfrac{1}{2^\mathit{d}}\bar{\mathbf{q}}^{*}_{\ell}$;
\ENDFOR
\end{algorithmic}
\end{algorithm}

\begin{algorithm}[hb]
\caption{Performing the second RK step in the LT approach.}
\begin{algorithmic}
\label{Steprk2}
\FOR{$\ell = L \to \ell_{\min}$, $\ell=\ell\!-\! 1$}
\IF{$\ell \ne L$}
\STATE Project the leaves and internal nodes from level $\ell+1$ onto level $\ell$ (Algorithm \ref{alg:projRK2});
\STATE Predict the values of the virtual leaves of level $\ell+1$ at instant $t^n + 2\Delta t_{\ell+1}$;\\
\COMMENT{These two steps compute the update in time of the virtual leaves}\\
\ENDIF

\FOR{every leaf $\in \Omega^{\ell}$}
\STATE Flux computation using the values at instant $t^n + \Delta t_\ell$;
\STATE Second step of compact RK;
\STATE Perform the approximation at instant $t^n+2\Delta t_{\ell}$ given in Equation (\ref{eq:projleaf});
\IF{we are computing RK2 time evolution}
\STATE Compute $\bar{\mathbf{q}}_{\theta = \frac{1}{2}}$;
\ELSIF{we are computing  RK3 time evolution}
\STATE Compute $\bar{\mathbf{q}}_{\theta = \frac{1}{4}}$; 
\STATE Compute $\bar{\mathbf{q}}_{\theta = \frac{3}{4}}$; 
\ENDIF
\ENDFOR
\ENDFOR
\end{algorithmic}
\end{algorithm}

\begin{algorithm}[hb]
\caption{Projection procedure inside the third RK step.}
\begin{algorithmic}
\label{alg:projRK3}
\REQUIRE Scale $\ell$ to receive the projection.
\REQUIRE Number of dimensions $\mathit{d}$ of the problem.
\FOR{every internal node $ \in \Omega^{\ell}$}
\STATE $\bar{\mathbf{q}}_{\ell} \left(t^n + \frac{1}{2}\Delta t_{\ell}\right) = 0$;
\COMMENT{result after RK2, set here to zero to store the projection}
\STATE
\FOR{Every child cell $\in \Omega^{\ell+1}$} 
\STATE $\bar{\mathbf{q}}_{\ell} \left(t^n + \frac{1}{2}\Delta t_{\ell}\right) \leftarrow \bar{\mathbf{q}}_{\ell} \left(t^n + \frac{1}{2}\Delta t_{\ell}\right) + \bar{\mathbf{q}}_{\ell+1} \left(t^n + \Delta t_{\ell+1}\right)$;
\COMMENT{Add value $\bar{\mathbf{q}}_{\ell+1}$ from RK3 $3^{\text rd}$ step }
\ENDFOR
\STATE
\STATE $\bar{\mathbf{q}}_{\ell} \left(t^n + \frac{1}{2}\Delta t_{\ell}\right) \leftarrow \dfrac{1}{2^\mathit{d}}\bar{\mathbf{q}}_{\ell} \left(t^n + \frac{1}{2}\Delta t_{\ell}\right)$
\STATE $\bar{\mathbf{q}}^{**}_{\ell} \approx \bar{\mathbf{q}}_{\ell} \left(t^n + \frac{1}{2}\Delta t_{\ell}\right)$;
\STATE
\STATE To obtain a $2^{\text{nd}}$ order approximation for $\bar{\mathbf{q}}_{\ell}(t^n + \Delta t_{\ell})$, use Equation~(\ref{rk3proj}).

\ENDFOR
\end{algorithmic}
\end{algorithm}
\begin{algorithm}[htb]
\caption{Performing the third RK step in the LT approach.}
\begin{algorithmic}
\label{Steprk3}

\FOR{$\ell = L \to \ell_{\min}$, $\ell\!=\!\ell\!-\!1$}
\IF{$\ell \ne L$}

\STATE -- Project the leaves and internal nodes from level $\ell+1$ onto level $\ell$ (Algorithm \ref{alg:projRK3});

\STATE -- Predict the values of the virtual leaves of level $\ell+1$ at instant $t^n + \Delta t_{\ell+1}$;\\
\COMMENT{These two steps compute the update in time of the virtual leaves}\\

\ENDIF

\FOR{Every leaf $\in \Omega^\ell$}
\STATE -- Flux computation using the values at instant $t^n + \frac{1}{2}\Delta t_\ell$;
\STATE -- Third step of compact RK;
\ENDFOR

\ENDFOR

\end{algorithmic}
\end{algorithm}

\begin{algorithm}
\caption{Tree refreshing before the second RK step.}
\begin{algorithmic}
\label{alg:refRK2}
\REQUIRE Scale $\ell$ to receive the projection.  \REQUIRE Number of dimensions $\mathit{d}$ of the problem.
\FOR{$\ell = L-1 \to \ell_{\min}$, $\ell\!=\!\ell\!-\!1$}
\FOR{Every internal node $ \in \Omega^{\ell}$}
\STATE $\bar{\mathbf{q}}_{\ell} \left(t^n + \frac{1}{2}\Delta t_{\ell}\right) = 0$;
\COMMENT{Set the solution equal zero to perform the averaging of its children cells.}\\
\FOR{Every child cell $i$ $\in \Omega^{\ell+1}$} 
\STATE $\bar{\mathbf{q}}_{\ell} \left(t^n + \frac{1}{2}\Delta t_{\ell}\right) \leftarrow \bar{\mathbf{q}}_{\ell} \left(t^n + \frac{1}{2}\Delta t_{\ell}\right) + \bar{\mathbf{q}}_{\ell+1, \; i}^* $ ;
\ENDFOR
\STATE $\bar{\mathbf{q}}_{\ell} \left(t^n + \frac{1}{2}\Delta t_{\ell}\right) \leftarrow \frac{1}{2^\mathit{d}}\bar{\mathbf{q}}_{\ell} \left(t^n + \frac{1}{2}\Delta t_{\ell}\right)$
\STATE Obtain a $1^{\text{st}}$ order approximation for $\bar{\mathbf{q}}_{\ell}^*$ using Equation (\ref{refreshRk2}).
\ENDFOR
\ENDFOR
\FOR{$\ell = \ell_{\min} \to L$, $\ell\!=\!\ell\!+\!1$}
\FOR{Every virtual leaf $ \in \Omega^{\ell}$}
\STATE Use the approximated solution at level $\ell-1$ at time instant $t^n + \frac{1}{2}\Delta t_{\ell-1}$ to predict the solution $\bar{\mathbf{q}}_{\ell}^{*}$.
\STATE Obtain the solution $\bar{\mathbf{q}}_{\ell}\left(t^n + \frac{1}{2}\Delta t_{\ell}\right)$ by linear interpolation.
\COMMENT{These two steps predict the solution of the virtual leaves in the proper time instant to update the level $\ell$.}\\

\ENDFOR
\ENDFOR
\end{algorithmic}
\end{algorithm}

\begin{algorithm}
\caption{Tree refreshing before the third RK step.}
\begin{algorithmic}
\label{alg:refRK22}
\REQUIRE Scale $\ell$ to receive the projection.
\REQUIRE Number of dimensions $\textit{d}$ of the problem.
\FOR{$\ell = L-1 \to \ell_{\min}$, $\ell\!=\!\ell\!-\!1$}
\FOR{Every internal node $ \in \Omega^{\ell}$}
\STATE $\bar{\mathbf{q}}_{\ell, \;\theta = \frac{1}{4}} = 0$;
\FOR{Every child cell $i$ $\in \Omega^{\ell+1}$} 
\STATE $\bar{\mathbf{q}}_{\ell, \;\theta = \frac{1}{4}} \leftarrow \bar{\mathbf{q}}_{\ell, \;\theta = \frac{1}{4}} + \bar{\mathbf{q}}_{\ell+1, \; i}^{**} $ ;
\ENDFOR
\STATE $\bar{\mathbf{q}}_{\ell, \;\theta = \frac{1}{4}} \leftarrow \frac{1}{2^\textit{d}}\bar{\mathbf{q}}_{\ell, \;\theta = \frac{1}{4}}$
\STATE Compute $\bar{\mathbf{q}}_{\ell, \; \theta = \frac{3}{4}}$ using Equation (\ref{ProjRK3equations1}).
\STATE Compute $\bar{\mathbf{q}}^{**}_{\ell}$ using Equation (\ref{ProjRK3equations2}).
\ENDFOR
\ENDFOR
\FOR{$\ell = \ell_{\min} \to L$, $\ell\!=\!\ell\!+\!1$}
\FOR{Every virtual leaf $ \in \Omega^{\ell}$}
\IF{$\ell = \ell_{\min}$}
\STATE Use the solution $\bar{\mathbf{q}}_{\ell-1, \; \theta = \frac{3}{4}}$, at time instant $t^n + \frac{1}{2}\Delta t_{\ell}$, to predict the solution $\bar{\mathbf{q}}_{\ell}^{**}$.
\ELSE
\STATE Use the solution $\bar{\mathbf{q}}_{\ell-1, \; \theta = \frac{1}{4}}$, at time instant $t^n + \frac{1}{2}\Delta t_{\ell}$, to predict the solution $\bar{\mathbf{q}}_{\ell}^{**}$.
\ENDIF
\STATE Compute $\bar{\mathbf{q}}_{\ell, \; \theta = \frac{1}{4}}$ using Equation (\ref{ProjRK3equations2}).
\STATE Compute $\bar{\mathbf{q}}_{\ell, \; \theta = \frac{3}{4}}$ using Equation (\ref{ProjRK3equations1}).
\ENDFOR
\ENDFOR
\end{algorithmic}
\end{algorithm}